\documentclass{amsart}

\usepackage{amssymb}

\allowdisplaybreaks

\newtheorem{theorem}{Theorem}[section]
\newtheorem{lemma}[theorem]{Lemma}
\newtheorem{proposition}[theorem]{Proposition}
\newtheorem{corollary}[theorem]{Corollary}
\newtheorem{remark}[theorem]{Remark}

\newtheorem{maintheo}{Main result}
\newtheorem{theoA}{Theorem}
\newtheorem{theoB}{Theorem}
\newtheorem{theoC}{Theorem}
\newtheorem{theoD}{Theorem}

\theoremstyle{definition}

\newcommand{\N}{\mathbb{N}}

\newcommand{\R}{\mathbb{R}}
\newcommand{\C}{\mathbb{C}}

\newcommand{\oM}{\otimes_{\mathcal{M}}}

\newcommand{\ten}{\otimes}

\newcommand{\dem}{\noindent {\bf Proof. }}
\newcommand{\demI}{\noindent {\bf Part I of the proof of Theorem \ref{Theorem-QS}. }}
\newcommand{\demII}{\noindent {\bf Part II of the proof of Theorem \ref{Theorem-QS}. }}
\newcommand{\demC}{\noindent {\bf Proof of Theorem \ref{QSQSQS}. }}
\newcommand{\fin}{\hspace*{\fill} $\square$ \vskip0.2cm}
\newcommand{\prodd}{\prod\nolimits}
\newcommand{\limm}{\lim\nolimits}

\def\bubl{{\displaystyle \mathop{\mathsf{A}}^{\circ}}}

\newcommand{\summ}{\sum\nolimits}

\begin{document}

\title[Operator space embedding
of $L_q$ into $L_p$] {Operator space embedding of $L_q$ into
$L_p$}

\author[Junge and Parcet]
{Marius Junge$^{\ast}$ and Javier Parcet$^{\dag}$}

\address{University of Illinois, Urbana, IL}

\email{junge@math.uiuc.edu}

\address{Centre de Recerca Matem{\`a}tica, Barcelona, Spain}

\email{javier.parcet@uam.es}

\footnote{$^{\ast}$Partially supported by the NSF DMS-0301116.}
\footnote{$^{\dag}$Partially supported by the Project
MTM2004-00678, Spain.} 


\maketitle

\section*{Introduction}

The idea of replacing functions by linear operators, the process
of quantization, goes back to the foundations of quantum mechanics
and has a great impact in mathematics. This applies for instance
to representation theory, operator algebra, noncommutative
geometry, quantum and free probability or operator space theory.
The quantization of measure theory leads to the theory of $L_p$
spaces defined over general von Neumann algebras, so called
\emph{noncommutative $L_p$ spaces}. This theory was initiated by
Segal, Dixmier and Kunze in the fifties and continued years later
by Haagerup, Fack, Kosaki and many others. We refer to the recent
survey \cite{PX2} for a complete exposition. In this paper we will
investigate noncommutative $L_p$ spaces in the language of
noncommutative Banach spaces, so called \emph{operator spaces}.
The theory of operator spaces took off in 1988 with Ruan's work
\cite{R}. Since then, it has been developed by Blecher/Paulsen,
Effros/Ruan and Pisier as a noncommutative generalization of
Banach space theory, see e.g. \cite{ER,Pau,P3}. In his book
\cite{P2} on vector valued noncommutative $L_p$ spaces, Pisier
considered a distinguished operator space structure on $L_p$. In
fact, the right category when dealing with noncommutative $L_p$ is
in many aspects that of operator spaces. Indeed, this has become
clear in the last years by recent results on noncommutative
martingales and related topics. In this paper, we prove a
fundamental structure theorem of $L_p$ spaces in the category of
operator spaces, solving a problem formulated by Gilles Pisier.

\begin{maintheo} Let $1 \le p < q \le 2$ and let $\mathcal{M}$ be a
von Neumann algebra. Then, there exists a sufficiently large von
Neumann algebra $\mathcal{A}$ and a completely isomorphic
embedding of $L_q(\mathcal{M})$ into $L_p(\mathcal{A})$, where
both spaces are equipped with their respective natural operator
space structures. Moreover, we have
\begin{itemize}
\item[\textbf{(a)}] If $\mathcal{M}$ is $\mathrm{QWEP}$, we can
choose $\mathcal{A}$ to be $\mathrm{QWEP}$.

\item[\textbf{(b)}] If $\mathcal{M}$ is hyperfinite, we can choose
$\mathcal{A}$ to be hyperfinite.
\end{itemize}
\end{maintheo}

In order to put our result in the right context, let us stress the
interaction between harmonic analysis, probability and Banach
space theory carried out mostly in the 70's. Based on previous
results by Beck, Grothendieck, L{\'e}vy, Orlicz, etc... probabilistic
methods in Banach spaces became the heart of the work developed by
Kwapie\'n, Maurey, Pisier, Rosenthal and many others. A
fundamental motivation for this new field relied on the embedding
theory of classical $L_p$ spaces. This theory was born in 1966
with the seminal paper \cite{BDK} of Bretagnolle, Dacunha-Castelle
and Krivine. They constructed an isometric embedding of $L_q$ into
$L_p$ for $1 \le p < q \le 2$, a Banach space version of our main
result. The simplest form of such embedding was known to L\'evy
and is given by
\begin{equation} \label{pstab}
\Big( \sum_{k=1}^\infty |\alpha_k|^q \Big)^{\frac1q} = \Big\|
\sum_{k=1}^\infty \alpha_k \, \theta_k \Big\|_{L_1(\Omega)},
\end{equation}
for scalars $(\alpha_k)_{k \ge 1}$ and where $(\theta_k)_{k \ge
1}$ is a suitable sequence of independent $q$-stable random
variables in $L_1(\Omega)$ for some probability space $(\Omega,
\mu)$. In other words, we have the relation $${\mathbb E} \exp
\big( i \summ_k \alpha_k \theta_k \big) = \exp \big( -c_q \summ_k
|\alpha_k|^q \big).$$

\vskip3pt

More recently, it has been discovered a parallel connection
between operator space theory and quantum probability. The
operator space version of Grothendieck theorem by Pisier and
Shlyakhtenko \cite{PS} and the embedding of $\mathrm{OH}$
\cite{J2} require tools from free probability. In this context we
should replace the $\theta_k$'s by suitable operators so that
\eqref{pstab} holds with matrix-valued coefficients $\alpha_1,
\alpha_2, \ldots$. To that aim, we develop new tools in quantum
probability and construct an operator space version of $q$-stable
random variables. To formulate this quantized form of
\eqref{pstab} we need some basic results of Pisier's theory
\cite{P2}. The most natural operator space structure on
$\ell_\infty$ comes from the diagonal embedding $\ell_\infty
\hookrightarrow \mathcal{B}(\ell_2)$. The natural structure on
$\ell_1$ is given by operator space duality, while the spaces
$\ell_p$ are defined by means of the complex interpolation method
\cite{P1} for operator spaces. Let us denote by $(\delta_k)_{k \ge
1}$ the unit vector basis of $\ell_q$. If $\widehat{\ten}$ denotes
the operator space projective tensor product and $S_p$ stands for
the Schatten $p$-class over $\ell_2$, it was shown in \cite{P2}
that
\[ \Big\| \sum_{k=1}^\infty a_k \ten \delta_k \Big\|_{S_1
\widehat{\ten} \, \ell_q} = \inf_{a_k = \alpha b_k \beta}
\|\alpha\|_{S_{2q'}} \Big( \sum_{k=1}^\infty \|b_k\|_{S_q}^q
\Big)^{\frac1q} \|\beta\|_{S_{2q'}}. \] The answer to Pisier's
problem for $\ell_q$ reads as follows.

\begin{theoA} \label{lp}
If $1 < q \le 2$, there exists a sufficiently large von Neumann
algebra $\mathcal{A}$ and a sequence $(x_k)_{k \ge 1}$ in
$L_1(\mathcal{A})$ such that the equivalence below holds for any
family $(a_k)_{k \ge 1}$ of trace class operators $$\inf_{a_k =
\alpha b_k \beta} \|\alpha\|_{2q'} \Big( \sum_{k=1}^\infty
\|b_k\|_q^q \Big)^{\frac1q} \|\beta\|_{2q'} \sim_c \Big\|
\sum_{k=1}^\infty a_k \ten x_k \Big\|_{L_1( \mathcal{A} \bar\ten
\mathcal{B}(\ell_2))}.$$
\end{theoA}

This gives a completely isomorphic embedding of $\ell_q$ into
$L_1(\mathcal{A})$. Moreover, the sequence $x_1, x_2, \ldots$
provides an operator space version of a $q$-stable process and
motivates a cb-embedding theory of $L_p$ spaces. A particular case
of Theorem A is the recent construction \cite{J2} of a
cb-embedding of Pisier's operator Hilbert space $\mathrm{OH}$ into
a von Neumann algebra predual. In other words, a complete
embedding of $\ell_2$ (with its natural operator space structure)
into a noncommutative $L_1$ space, see also Pisier's paper
\cite{P4} for a shorter proof and Xu's alternative construction
\cite{X2}. Other related results appear in \cite{J3,P5,PS,X3},
while semi-complete embeddings between (vector-valued) $L_p$
spaces can be found in \cite{JP,Pa}. All these papers will play a
role in our analysis.

\vskip5pt

Let us sketch the simplest construction (a more elaborated one is
needed to prove the stability of hyperfiniteness) which leads to
this operator space version of $q$-stable random variables. A key
ingredient in our proof is the notion of the Haagerup tensor
product $\ten_h$. We first note that $\ell_q$ is the diagonal
subspace of the Schatten class $S_q$. According to \cite{P2},
$S_q$ can be written as the Haagerup tensor product of its first
column and first row subspaces $S_q = C_q \otimes_h R_q$.
Moreover, using a simple generalization of \lq\lq Pisier's
exercise\rq\rq${}$ (see Exercise 7.9 in Pisier's book \cite{P3})
we have
\begin{equation} \label{Pexercise}
\begin{array}{rcl} C_q & \hookrightarrow_{cb} & \big( R \oplus
\mathrm{OH} \big) \big/ graph(\Lambda_1)^\perp, \\ [5pt] R_q &
\hookrightarrow_{cb} & \big( C \oplus \mathrm{OH} \big) \big/
graph(\Lambda_2)^\perp,
\end{array}
\end{equation}
with $\Lambda_1: C \to \mathrm{OH}$ and $\Lambda_2: R \to
\mathrm{OH}$ suitable injective, closed, densely-defined operators
with dense range and where $\hookrightarrow_{cb}$ denotes a
cb-embedding, see \cite{J2,P4,X3} for related results. By
\cite{J2} and duality, it suffices to see that $graph(\Lambda_1)
\otimes_h graph(\Lambda_2)$ is cb-isomorphic to a cb-complemented
subspace of $\mathcal{A}(\mathrm{OH})$. By the injectivity of the
Haagerup tensor product, we note that $graph(\Lambda_1) \otimes_h
graph(\Lambda_2)$ is an intersection of four spaces. Let us
explain this in detail. By discretization we may assume that
$\Lambda_j = \mathsf{d}_{\lambda^j} = \sum_k \lambda_k^j e_{kk}$
is a diagonal operator on $\ell_2$. In fact, by a simple
complementation argument, it is no restriction to assume that the
eigenvalues are the same for $j=1,2$. Thus we may write the
Haagerup tensor product above as follows
$$\mathcal{J}_{\infty,2} = graph(\mathsf{d}_\lambda) \otimes_h
graph(\mathsf{d}_\lambda) = \Big( C \cap \ell_2^{oh}(\lambda)
\Big) \otimes_h \Big( R \cap \ell_2^{oh}(\lambda) \Big),$$ where
$\ell_2^{oh}(\lambda)$ is a weighted form of $\mathrm{OH}$
according to the action of $\mathsf{d}_\lambda$. That is
\begin{eqnarray*}
C \cap \ell_2^{oh} (\lambda) & = & \Big\{ ( \hskip0.5pt e_{i1}
\hskip0.5pt, \lambda_i \hskip0.5pt e_{i1} \hskip1pt) \,
\hskip0.5pt \big| \ \hskip0.5pt i \hskip0.5pt \ge 1 \Big\} \subset
C \oplus \mathrm{OH}, \\ R \cap \ell_2^{oh} (\lambda) & = & \Big\{
(e_{1j}, \lambda_j e_{1j}) \, \big| \ j \ge 1 \Big\} \subset R
\hskip0.5pt \oplus \mathrm{OH}.
\end{eqnarray*}
The symbol $\infty$ in $\mathcal{J}_{\infty,2}$ is used because we
shall consider $L_p$ versions of these spaces along the paper. The
number $2$ denotes that this space arises as a \lq middle
point\rq${}$ in the sense of interpolation theory between two
related $\mathcal{J}$-spaces, see Section \ref{NewSect3} below for
more details. Now, regarding $\mathsf{d}_\lambda^4 = \sum_k
\lambda_k^4 e_{kk}$ as the density $d_\psi$ of some normal
strictly semifinite faithful weight $\psi$ on
$\mathcal{B}(\ell_2)$, the space $\mathcal{J}_{\infty,2}$ splits
up into a $4$-term intersection space. In other words, we find
$$\mathcal{J}_{\infty,2}(\psi) = \big( C \ten_h R \big) \cap
\big( C \ten_h \mathrm{OH} \big) d_\psi^{\frac14} \cap
d_\psi^{\frac14} \big( \mathrm{OH} \ten_h R \big) \cap
d_\psi^{\frac14} \big( \mathrm{OH} \ten_h \mathrm{OH} \big)
d_\psi^{\frac14}.$$ The norm of $x$ in
$\mathcal{J}_{\infty,2}(\psi)$ is given by
$$\max \Big\{ \|x\|_{B(\ell_2)}, \big\|x d_\psi^{\frac14}
\big\|_{C \ten_h \mathrm{OH}}, \big\| d_\psi^{\frac14} x
\big\|_{\mathrm{OH} \ten_h R}, \big\| d_\psi^{\frac14} x
d_\psi^{\frac14} \big\|_{\mathrm{OH} \ten_h \mathrm{OH}} \Big\}.$$
The two middle terms are not as unusual as it might seem
\begin{equation} \label{osstruct}
\begin{array}{rcl}
\big\| d_\psi^{\frac14} (x_{ij}) \big\|_{\mathrm{M}_m(\mathrm{OH}
\ten_h R)} & = & \displaystyle \sup_{\|\alpha\|_{S_4^m}\le 1}
\Big\| d_\psi^{\frac14} \Big( \sum_{k=1}^m \alpha_{ik} x_{kj}
\Big) \Big\|_{L_4(\mathrm{M}_m \ten \mathcal{B}(\ell_2))}, \\
\big\| (x_{ij}) d_\psi^{\frac14} \big\|_{\mathrm{M}_m(C \ten_h
\mathrm{OH})} & = & \displaystyle \sup_{\|\beta\|_{S_4^m}\le 1}
\Big\| \Big( \sum_{k=1}^m x_{ik} \beta_{kj} \Big) \hskip1pt
d_\psi^{\frac14} \Big\|_{L_4(\mathrm{M}_m \ten
\mathcal{B}(\ell_2))}.
\end{array}
\end{equation}
Let us now assume that we just try to embed the finite-dimensional
space $S_q^m$. By approximation, it suffices to consider only
finitely many eigenvalues $\lambda_1, \lambda_2, \ldots,
\lambda_n$ and according to the results from \cite{J3}, we can
take $n \sim m \log m$. In this case we rename $\psi$ by $\psi_n$
and we may easily assume that $\mbox{tr} (d_{\psi_n}) = \sum_k
\lambda_k^4 = \mathrm{k}_n$ is an integer. Therefore, we may
consider the following state on $\mathcal{B}(\ell_2^n)$
$$\varphi_n(x) = \frac{1}{\mathrm{k}_n} \sum_{k=1}^n \lambda_k^4
\, x_{kk}.$$ In this particular case, the space
$\mathcal{J}_{\infty,2}(\psi_n)$ can be obtained using free
probability. More precisely, we may identify it as a subspace of
$L_\infty(\mathcal{A}; \mathrm{OH}_{\mathrm{k}_n})$, which in turn
is the space obtained by complex interpolation
$L_\infty(\mathcal{A}; \mathrm{OH}_{\mathrm{k}_n}) =
[C_{\mathrm{k}_n}(\mathcal{A}),
R_{\mathrm{k}_n}(\mathcal{A})]_{1/2}$ where
$R_{\mathrm{k}_n}(\mathcal{A})$ and
$C_{\mathrm{k}_n}(\mathcal{A})$ denote the row and column
subspaces of $\mathrm{M}_{\mathrm{k}_n}(\mathcal{A})$.

\begin{theoB} \label{MAMSth}
Let $\mathcal{A} = \mathsf{A}_1 * \mathsf{A}_2 * \cdots *
\mathsf{A}_{\mathrm{k}_n}$ be the reduced free product of
$\mathrm{k}_n$ copies of $\mathcal{B}(\ell_2^n) \oplus
\mathcal{B}(\ell_2^n)$ equipped with the state $\frac12 (\varphi_n
\oplus \varphi_n)$. If $\pi_k: \mathsf{A}_j \to \mathcal{A}$
denotes the canonical embedding into the $j$-th component of
$\mathcal{A}$, the mapping $$u_n: x \in
\mathcal{J}_{\infty,2}(\psi_n) \mapsto \sum_{j=1}^{\mathrm{k}_n}
\pi_j(x,-x) \ten \delta_j \in L_\infty(\mathcal{A};
\mathrm{OH}_{\mathrm{k}_n})$$ is a cb-embedding with
cb-complemented image and constants independent of $n$.
\end{theoB}

Theorem B and its generalization for arbitrary von Neumann
algebras is a very recent result from \cite{JP2}. However, the
proof given there is rather long and quite technical, see Section
\ref{NewSect3} for further details. In order to be more
self-contained, we provide a second proof only using elementary
tools from free probability. We think this alternative approach is
of independent interest. Now, by duality we obtain a cb-embedding
of $\ell_q^m \hookrightarrow_{cb} S_q^m$ into $L_1(\mathcal{A};
\mathrm{OH}_{\mathrm{k}_n})$. Then we use the cb-embedding
\cite{J2} of $\mathrm{OH}$ into a von Neumann algebra predual.
Using an ultraproduct procedure, the desired cb-embedding of
$\ell_q$ into $L_1(\mathcal{A})$ is obtained. In fact, what we
shall prove is a far reaching generalization of Theorem \ref{lp}.
Namely, the same result holds replacing $C_q$ and $R_q$ by
subspaces of quotients of $R \oplus \mathrm{OH}$ and $C \oplus
\mathrm{OH}$ respectively.

\begin{theoC} \label{QSQSQS}
Let $\mathrm{X}_1$ be a subspace of a quotient $R \oplus
\mathrm{OH}$ and let $\mathrm{X}_2$ be a subspace of a quotient
of $C \oplus \mathrm{OH}$. Then, there exist some $\mathrm{QWEP}$
algebra $\mathcal{A}$ and a cb-embedding
$$\mathrm{X}_1 \otimes_h \mathrm{X}_2 \hookrightarrow_{cb}
L_1(\mathcal{A}).$$
\end{theoC}

Now we may explain how the paper is organized. In Section
\ref{NewSect1} we just prove the complete embeddings
\eqref{Pexercise}. This is a simple consequence of Pisier's
exercise 7.9 in \cite{P3} and Xu's generalization \cite{X3}.
However, we have decided to include the proof since it serves as a
model in our construction of the cb-embedding for general von
Neumann algebras. In Section \ref{NewSect2} we concentrate on the
simplest form of our main result by proving Theorems \ref{lp},
\ref{MAMSth}, \ref{QSQSQS}. This is done in part to motivate a new
class of noncommutative function spaces introduced in Section
\ref{NewSect3}. We will state (in terms of these new spaces) the
main result of \cite{JP2}, a further generalization of Theorem B
and a key tool in the proof of our main result. This is a free
analogue of a generalization of Rosenthal's inequality, as we
shall explain in Section \ref{NewSect3}. The whole Section
\ref{NewSect4} is entirely devoted to the proof of our main
result. We first construct a cb-embedding of $S_q$ into
$L_p(\mathcal{A})$ for some von Neumann algebra $\mathcal{A}$.
This is quite similar to our argument sketched above and provides
an $L_p$ generalization of Theorem \ref{QSQSQS}, but this
construction does not preserve hyperfiniteness. The argument to
fix this is quite involved and requires recent techniques from
\cite{J3} and \cite{J6}. We first apply a transference argument,
via a noncommutative Rosenthal type inequality in $L_1$ for
identically distributed random variables, to replace freeness in
our construction by some sort of noncommutative independence. This
allows to avoid free product von Neumann algebras and use tensor
products instead. Then we combine the algebraic central limit
theorem with the notion of noncommutative Poisson random measure
to eliminate the use of ultraproducts in the process. After these
modifications in our original argument, it is easily seen that
hyperfiniteness is preserved. This more involved construction of
the cb-embedding is the right one to analyze the finite
dimensional case. In other words, we estimate the dimension of
$\mathcal{A}$ in terms of the dimension of $\mathcal{M}$, see
Remark \ref{estfindim} below for details. The proof for general
von Neumann algebras requires to consider a \lq generalized\rq${}$
Haagerup tensor product since we are not in the discrete case
anymore.

\vskip5pt

We conclude with some comments related to our main result. In the
category of Banach spaces, noncommutative versions of a $p$-stable
process were studied in \cite{J0} and further analyzed in
\cite{J6} to construct Banach space embeddings between
noncommutative $L_p$ spaces. If $0 < p < q < 2$ and we write
$\mathcal{R}$ for the hyperfinite $\mathrm{II}_1$ factor, the main
embedding result there asserts that the space $L_q(\mathcal{R}
\bar\ten \mathcal{B}(\ell_2))$ embeds isometrically into
$L_p(\mathcal{R})$. One of the principal techniques in the proof
is a noncommutative version of the notion of Poisson random
measure, which will also be instrumental in this paper. On the
other hand, the cb-embedding theory of $L_p$ spaces presents
significant differences. Indeed, in sharp contrast with the
classical theory, it was proved in \cite{J00} that the operator
space $\mathrm{OH}$ does not embed completely isomorphically into
any $L_p$ space for $2 < p < \infty$. Moreover, after \cite{P5} we
know that there is no possible cb-embedding of $\mathrm{OH}$ into
the predual of a semifinite von Neumann algebra. As it is to be
expected, this also happens in our main result and justifies the
relevance of type III von Neumann algebras in the subject.

\begin{theoD}
If $1 \le p < q \le 2$ and $\ell_q$ cb-embeds into
$L_p(\mathcal{A})$, $\mathcal{A}$ is not semifinite.
\end{theoD}


Unfortunately, the proof of this result is out of the scope of
this paper and will be the subject of a forthcoming publication.
In fact, the proof for the case $p=1$ is much harder than the case
$1 < p < q$ and requires the use of a noncommutative version of
Rosenthal theorem \cite{Ro}, recently obtained in \cite{JP3}. Our
results there are closely related to this paper and complement
Pisier's paper \cite{Pis} and some recent results of
Randrianantoanina \cite{Narcisse}. After a quick look at the main
results in \cite{BDK,J0}, the problem of constructing an
\emph{isometric} cb-embedding of $L_q$ into $L_p$ arises in a
natural way. This remains an open problem.

\vskip0.5cm

\subsection*{Background and notation}

We shall assume that the reader is familiar with those branches of
operator algebra related to the theories of operator spaces and
noncommutative $L_p$ spaces. The recent monographs \cite{ER} and
\cite{P3} on operator spaces contain more than enough information
for our purposes. We shall work over general von Neumann algebras
so that we use Haagerup's definition \cite{H} of $L_p$, see also
Terp's excellent exposition of the subject \cite{T1}. As is well
known, Haagerup $L_p$ spaces have trivial intersection and thereby
do not form an interpolation scale. However, the complex
interpolation method will be a basic tool in this paper. This is
solved using Kosaki's definition \cite{Ko} of $L_p$. We also refer
to Chapter 1 in \cite{JP2} or to the survey \cite{PX2} for a quick
review of Haagerup's and Kosaki's definitions of noncommutative
$L_p$ spaces and the compatibility between them. In particular,
using such compatibility, we shall use in what follows the complex
interpolation method without further details. The basics on von
Neumann algebras and Tomita's modular theory required to work with
these notions can be found in Kadison/Ringrose books \cite{KR}.
There are some other topics related to noncommutative $L_p$ spaces
that will be frequently used in the paper. The main properties of
normal faithful conditional expectations over Haagerup $L_p$
spaces can be found in \cite{JX} and \cite{Ta2}. We shall also
assume certain familiarity with Pisier's theory of vector-valued
non-commutative $L_p$ spaces \cite{P2} and Voiculescu's free
probability theory \cite{VDN}. Along the text we shall find some
other topics such as certain noncommutative function spaces, some
recent inequalities for free random variables, a noncommutative
version of a Poisson process, etc... In all these cases our
exposition intends to be self-contained.

\vskip5pt

We shall follow the standard notation in the subject. Anyway, let
us say a few words on our terminology. The symbols $(\delta_k)$
and $(e_{ij})$ will denote the unit vector basis of $\ell_2$ and
$\mathcal{B}(\ell_2)$ respectively. The letters $\mathcal{A},
\mathcal{M}$ and $\mathcal{N}$ are reserved to denote von Neumann
algebras. Almost all the time, the inclusions $\mathcal{N} \subset
\mathcal{M} \subset \mathcal{A}$ will hold. We shall use $\varphi$
and $\phi$ to denote normal faithful (\emph{n.f.} in short)
states, while the letter $\psi$ will be reserved for normal
strictly semifinite faithful (\emph{n.s.s.f.} in short) weights.
Inner products and duality brackets will be anti-linear on the
first component and linear on the second component. As usual,
given an operator space $\mathrm{X}$ we shall write
$\mathrm{M}_m(\mathrm{X})$ for the space of $m \times m$ matrices
with entries in $\mathrm{X}$ and we shall equip it with the norm
of the minimal tensor product $\mathrm{M}_m \ten_{\mathrm{min}}
\mathrm{X}$. Similarly, the $\mathrm{X}$-valued Schatten $p$-class
over $\mathrm{M}_m$ will be denoted by $S_p^m(\mathrm{X})$.
Accordingly, $L_p(\mathcal{M}; \mathrm{X})$ will stand for the
$\mathrm{X}$-valued $L_p$ space over $\mathcal{M}$. Given $\gamma
> 0$, we shall write $\gamma \mathrm{X}$ to denote the space
$\mathrm{X}$ equipped with the norm $\|x\|_{\gamma \mathrm{X}} =
\gamma \|x\|_{\mathrm{X}}$. In particular, if $\mathcal{M}$ is a
finite von Neumann algebra equipped with a finite weight $\psi$,
we shall usually write $\psi = \mathrm{k} \varphi$ with
$\mathrm{k} = \psi(1)$ so that $\varphi$ becomes a state on
$\mathcal{M}$. In this situation, the associated $L_p$ space will
be denoted as $\mathrm{k}^{1/p} L_p(\mathcal{M})$, so that the
$L_p$ norm is calculated using the state $\varphi$. Any new or
non-standard terminology will be properly introduced in the text.

\numberwithin{equation}{section}

\section{On \lq\lq Pisier's exercise\rq\rq}
\label{NewSect1}

We begin by proving a generalization of Exercise 7.9 in \cite{P3}.
This result became popular after Pisier applied it in \cite{P4} to
obtain a simpler way to cb-embed $\mathrm{OH}$ into the predual of
a von Neumann algebra. In fact, our argument is quite close to the
one given in \cite{X3} for a similar result and might be known to
experts. Nevertheless we include it here for completeness, since
it will be used below and mainly because it will also serve as a
model for our construction in the non-discrete case. Let us set
some notation. Given a Hilbert space $\mathcal{H}$, we shall write
$\mathcal{H}_r = \mathcal{B}(\mathcal{H}, \C)$ and $\mathcal{H}_c
= \mathcal{B}(\C, \mathcal{H})$ for the row and column
quantizations on $\mathcal{H}$. Moreover, given $1 \le p \le
\infty$ we shall use the following terminology
$$\mathcal{H}_{r_p} = \big[ \mathcal{H}_r, \mathcal{H}_c
\big]_{\frac1p} \quad \mbox{and} \quad \mathcal{H}_{c_p} = \big[
\mathcal{H}_c, \mathcal{H}_r\big]_{\frac1p}.$$ There are two
particular cases for which we use another terminology. When
$\mathcal{H} = \ell_2$, we shall use $(R,C,R_p,C_p)$ instead.
Moreover, when the Hilbert space is $L_2(\mathcal{M})$ for some
von Neumann algebra $\mathcal{M}$, we shall write
$L_2^{r_p}(\mathcal{M})$ and $L_2^{c_p}(\mathcal{M})$. In the same
fashion, $\mathcal{H}_{oh}$ and $L_2^{oh}(\mathcal{M})$ stand for
the operator Hilbert space structures. Given two operator spaces
$\mathrm{X}_1$ and $\mathrm{X}_2$, the expression $\mathrm{X}_1
\simeq_{cb} \mathrm{X}_2$ means that there exists a complete
isomorphism between them. We shall write $\mathrm{X}_1 \in
\mathcal{QS}(\mathrm{X}_2)$ to denote that $\mathrm{X}_1$ is
completely isomorphic to a quotient of a subspace of
$\mathrm{X}_2$. Let $\mathcal{S}$ denote the strip
$$\mathcal{S} = \Big\{ z \in \C \, \big| \ 0 < \mbox{Re}(z) < 1
\Big\}$$ and let $\partial \mathcal{S} =
\partial_0 \cup \partial_1$ be the partition of its boundary into
$$\partial_0 = \Big\{ z \in \C \, \big| \ \mbox{Re}(z)=0 \Big\}
\quad \mbox{and} \quad \partial_1 = \Big\{ z \in \C \, \big| \
\mbox{Re}(z)=1 \Big\}.$$ Given $0 < \theta < 1$, let $\mu_\theta$
be the harmonic measure of the point $z = \theta$. This is a
probability measure on $\partial \mathcal{S}$ (with density given
by the Poisson kernel in the strip) that can be written as
$\mu_\theta = (1 - \theta) \mu_0 + \theta \mu_1$, with $\mu_j$
being probability measures supported by $\partial_j$ and such that
\begin{equation} \label{Eq-AnalyticCondition}
f(\theta) = \int_{\partial \mathcal{S}} f d\mu_\theta
\end{equation}
for any bounded analytic function $f: \mathcal{S} \to \C$ extended
non-tangentially to $\partial \mathcal{S}$. Now, before proving
the announced result, we need to set a formula describing the norm
of certain kind of vector-valued noncommutative $L_p$ space. A
more detailed account of these expressions will be given at the
beginning of Section \ref{NewSect3}. Given $2 \le p \le \infty$
and $0 < \theta < 1$, the norm of $x = \sum_k x_k \ten \delta_k$
in $\big[ S_p(C_p), S_p(R_p) \big]_{\theta}$ is given by
\begin{equation} \label{fbasexu}
\sup \Big\{ \Big( \summ_k \big\| \alpha x_k \beta \big\|_{S_2}^2
\Big)^{\frac12} \big| \ \|\alpha\|_{S_u}, \|\beta\|_{S_v} \le 1
\Big\}
\end{equation}
where the indices $(u,v)$ are determined by $$(1/u,1/v) =
(\theta/q, (1-\theta)/q) \quad \mbox{with} \quad 1/2 = 1/p +
1/q.$$ Of course, this formula trivially generalizes for the norm
in the space $$\big[ S_p(\mathcal{H}_{c_p}),
S_p(\mathcal{H}_{r_p}) \big]_{\theta}.$$

\begin{lemma} \label{Lemma-Motivation}
If $1 \le p < q \le 2$, we have
$$R_q \in \mathcal{QS} \big( R_p \oplus_2
\mathrm{OH} \big) \quad \mbox{and} \quad C_q \in \mathcal{QS}
\big( C_p \oplus_2 \mathrm{OH} \big).$$
\end{lemma}

\dem We only prove the first assertion since the arguments for
both are the same. In what follows we fix $0 < \theta < 1$
determined by the relation $R_q = [R_p, \mathrm{OH}]_\theta$. In
other words, we have $1/q = (1-\theta)/p + \theta/2$. According to
the complex interpolation method and its operator space extension
\cite{P1}, given a compatible couple $(\mathrm{X}_0,
\mathrm{X}_1)$ of operator spaces we define
$\mathcal{F}(\mathrm{X}_0, \mathrm{X}_1)$ as the space of bounded
analytic functions $f: \mathcal{S} \to \mathrm{X}_0 +
\mathrm{X}_1$ and we equip it with the following norm
$$\|f\|_{\mathcal{F}(\mathrm{X}_0, \mathrm{X}_1)} = \Big(
(1-\theta) \,
\|{f_{\mid_{\partial_0}}}\|_{L_2(\partial_0;\mathrm{X}_0)}^2 +
\theta \,
\|{f_{\mid_{\partial_1}}}\|_{L_2(\partial_1;\mathrm{X}_1)}^2
\Big)^{\frac12}.$$ Then, the complex interpolation space
$\mathrm{X}_\theta = [\mathrm{X}_0, \mathrm{X}_1]_{\theta}$ can be
defined as the space of all $x \in \mathrm{X}_0 + \mathrm{X}_1$
such that there exists a function $f \in \mathcal{F}(\mathrm{X}_0,
\mathrm{X}_1)$ with $f(\theta) = x$. We equip $\mathrm{X}_\theta$
with the norm $$\|x\|_{\mathrm{X}_\theta} = \inf \Big\{
\|f\|_{\mathcal{F}(\mathrm{X}_0, \mathrm{X}_1)} \, \big| \ f \in
\mathcal{F}(\mathrm{X}_0, \mathrm{X}_1) \ \mbox{and} \ f(\theta) =
x \Big\}.$$ In our case we set $\mathrm{X}_0 = R_p$ and
$\mathrm{X}_1 = \mathrm{OH}$. If we define $$\mathcal{H} =
(1-\theta)^{\frac12} L_2(\partial_0; \ell_2) \quad \mbox{and}
\quad \mathcal{K} = \theta^\frac12 L_2(\partial_1; \ell_2),$$ it
turns out that $\mathcal{F}(R_p,\mathrm{OH})$ can be regarded (via
Poisson integration) as a subspace of $\mathcal{H} \oplus_2
\mathcal{K}$. Moreover, we equip $\mathcal{F}(R_p,\mathrm{OH})$
with the operator space structure inherited from
$\mathcal{H}_{r_p} \oplus_2 \mathcal{K}_{oh}$. Then, we define the
mapping $$\mathcal{Q}: f \in \mathcal{F}(R_p,\mathrm{OH}) \mapsto
f(\theta) \in R_q.$$ The assertion follows from the fact that
$\mathcal{Q}$ is a complete metric surjection. Indeed, in that
case we have $R_q \simeq_{cb} \mathcal{F}(R_p,\mathrm{OH}) / \ker
\mathcal{Q}$, which is a quotient of a subspace of $R_p \oplus_2
\mathrm{OH}$. In order to see that $\mathcal{Q}$ is a complete
surjection, it suffices to see that the map $id_{S_{p'}} \otimes
\mathcal{Q}: S_{p'}(\mathcal{F}(R_p,\mathrm{OH})) \to S_{p'}(R_q)$
is a metric surjection. We begin by showing that $id_{S_{p'}}
\otimes \mathcal{Q}$ is contractive. Let $f \in
S_{p'}(\mathcal{F}(R_p,\mathrm{OH}))$ be of norm $< 1$ and let us
write $$f(\theta) = \summ_k f_k(\theta) \otimes \delta_k \in
S_{p'}(R_q).$$ To compute the norm of $f(\theta)$ we note that
$$S_{p'}(R_q) = \big[ S_{p'}(C_{p'}), S_{p'}(R_{p'}) \big]_\eta
\quad \mbox{with} \quad 1/q = (1-\eta)/p + \eta/p'.$$ Then it
follows from \eqref{fbasexu} that
\begin{equation} \label{Eq-int-f(theta)}
\|f(\theta)\|_{S_{p'}(R_q)} = \sup \Big\{ \Big( \summ_k \big\|
\alpha f_k(\theta) \beta \big\|_{S_2}^2 \Big)^{\frac12} \, \big| \
\|\alpha\|_{S_{2r/\eta}}, \|\beta\|_{S_{2r/(1-\eta)}} \le 1 \Big\}
\end{equation}
where $1/2r = 1/2 - 1/p' = 1/p - 1/2$. Moreover, it is clear that
we can restrict the supremum above to all $\alpha$ and $\beta$ in
the positive parts of their respective unit balls. Taking this
restriction in consideration, we define
$$g(z) = \summ_k g_k(z) \otimes \delta_k \quad \mbox{with} \quad
g_k(z) = \alpha^{\frac{z}{\theta}} f_k(z)
\beta^{\frac{2-z}{2-\theta}}.$$ The $g_k$'s are analytic in
$\mathcal{S}$ and take values in $S_2$. Thus, we have the identity
\begin{eqnarray} \label{Eq-id-harm}
\lefteqn{\Big( \summ_k \big\| \alpha f_k(\theta) \beta
\big\|_{S_2}^2 \Big)^{\frac12}} \\ \nonumber & = & \Big( \summ_k
\|g_k(\theta)\|_{S_2}^2 \Big)^{\frac12} \\ \nonumber & = & \Big(
(1-\theta) \int_{\partial_0} \summ_k \|g_k(z)\|_{S_2}^2 d\mu_0 +
\theta \int_{\partial_1} \summ_k \|g_k(z)\|_{S_2}^2 d\mu_1
\Big)^{\frac12}.
\end{eqnarray}
The contractivity of $id_{S_{p'}} \otimes \mathcal{Q}$ will follow
from
\begin{eqnarray}
\label{Est000} \int_{\partial_0} \summ_k \|g_k(z)\|_{S_2}^2 d\mu_0
& \le & \|f_{\mid_{\partial_0}}\|_{S_{p'}
(L_2^{r_p}(\partial_0;\ell_2))}^2, \\ \label{Est111}
\int_{\partial_1} \summ_k \|g_k(z)\|_{S_2}^2 d\mu_1 & \le &
\|f_{\mid_{\partial_1}}\|_{S_{p'}(L_2^{oh}(\partial_1;\ell_2))}^2.
\end{eqnarray}
Indeed, if we combine \eqref{Eq-int-f(theta)} and
\eqref{Eq-id-harm} with the operator space structure defined on
$\mathcal{F}(R_p, \mathrm{OH})$, it turns out that inequalities
\eqref{Est000} and \eqref{Est111} are exactly what we need. To
prove \eqref{Est000} we observe that $2 \eta = \theta$ and
$r'=p'/2$, so that
\begin{eqnarray*}
\int_{\partial_0} \summ_k \|g_k(z)\|_{S_2}^2 d\mu_0 & = &
\int_{\partial_0} \summ_k \|f_k(z)
\beta^{\frac{2}{2-\theta}}\|_{S_2}^2 d\mu_0 \\ & = &
\int_{\partial_0} \summ_k \mbox{tr} \big( f_k(z)^* f_k(z)
\beta^{\frac{1}{1-\eta}} (\beta^{\frac{1}{1-\eta}})^* \big) d\mu_0
\\ & \le & \|\beta^{\frac{2}{1-\eta}}\|_{S_r} \Big\|
\int_{\partial_0} \summ_k f_k(z)^* f_k(z) d\mu_0 \Big\|_{S_{r'}}.
\end{eqnarray*}
This gives $$\int_{\partial_0} \summ_k \|g_k(z)\|_{S_2}^2 d\mu_0
\le \|\beta\|_{S_{2r/(1-\eta)}}^{\frac{2}{1-\eta}} \Big\| \Big(
\int_{\partial_0} \summ_k f_k(z)^* f_k(z) d\mu_0 \Big)^{\frac12}
\Big\|_{S_{p'}}^2.$$ The first term on the right is $\le 1$.
Hence, we have
\begin{eqnarray*}
\int_{\partial_0} \summ_k \|g_k(z)\|_{S_2}^2 d\mu_0 & \le &
\|f_{\mid_{\partial_0}}\|_{S_{p'}(L_2^{c_{p'}}(\partial_0;
\ell_2))}^2 = \|f_{\mid_{\partial_0}}\|_{S_{p'}
(L_2^{r_p}(\partial_0;\ell_2))}^2.
\end{eqnarray*}
This proves \eqref{Est000} while for \eqref{Est111} we proceed in
a similar way
\begin{eqnarray*}
\lefteqn{\int_{\partial_1} \summ_k \|g_k(z)\|_{S_2}^2 d\mu_1} \\ &
= & \int_{\partial_1} \summ_k \| \alpha^{\frac{1}{\theta}} f_k(z)
\beta^{\frac{1}{2-\theta}}\|_{S_2}^2 d\mu_1 \\ & = &
\int_{\partial_1} \summ_k \| \alpha^{\frac{1}{2 \eta}} f_k(z)
\beta^{\frac{1}{2- 2 \eta}}\|_{S_2}^2 d\mu_1 \\ & \le & \sup
\left\{ \int_{\partial_1} \summ_k \| a f_k(z) b \|_{S_2}^2 d\mu_1
\, \big| \ \|a\|_{S_{4r}}, \|b\|_{S_{4r}} \le 1 \right\}
\\ & = & \sup \left\{ \Big\| a \Big(
\summ_k {f_k}_{\mid_{\partial_1}} \otimes \delta_k \Big) b
\Big\|_{S_2(L_2^{oh}(\partial_1; \ell_2))}^2 \, \big| \
\|a\|_{S_{4r}}, \|b\|_{S_{4r}} \le 1 \right\}.
\end{eqnarray*}
According to \eqref{fbasexu}, this proves \eqref{Est111} and we
have a contraction. Reciprocally, given $x \in S_{p'}(R_q)$ of
norm $< 1$, we are now interested on finding $f \in
S_{p'}(\mathcal{F}(R_p, \mathrm{OH}))$ such that $f(\theta) = x$
and $\|f\|_{S_{p'}(\mathcal{F}(R_p, \mathrm{OH}))} \le 1$. Since
$$[S_{p'}(R_p), S_{p'}(\mathrm{OH})]_\theta = S_{p'}(R_q)$$ there
must exists $f \in \mathcal{F}(S_{p'}(R_p), S_{p'}(\mathrm{OH}))$
such that $f(\theta) = x$ and
$$\|f\|_{\mathcal{F}(S_{p'}(R_p),
S_{p'}(\mathrm{OH}))} = \Big( (1-\theta)
\|f_{\mid_{\partial_0}}\|_{L_2(\partial_0;S_{p'}(R_p))}^2 + \theta
\|f_{\mid_{\partial_1}}\|_{L_2(\partial_1;S_{p'}(\mathrm{OH}))}^2
\Big)^{\frac12} \le 1.$$ Therefore, it remains to see that
\begin{eqnarray}
\label{Eq-222}
\|f_{\mid_{\partial_0}}\|_{S_{p'}(L_2^{r_p}(\partial_0;\ell_2))} &
\le & \|f_{\mid_{\partial_0}}\|_{L_2(\partial_0;S_{p'}(R_p))}, \\
\label{Eq-333}
\|f_{\mid_{\partial_1}}\|_{S_{p'}(L_2^{oh}(\partial_1;\ell_2))} &
\le &
\|f_{\mid_{\partial_1}}\|_{L_2(\partial_1;S_{p'}(\mathrm{OH}))}.
\end{eqnarray}
However, the identities below are clear by now
\begin{eqnarray*}
\|f_{\mid_{\partial_0}}\|_{S_{p'}(L_2^{r_p}(\partial_0;\ell_2))} &
= & \Big\| \Big( \int_{\partial_0} \summ_k f_k(z)^* f_k(z) d\mu_0
\Big)^{\frac12}\Big\|_{S_{p'}}, \\
\|f_{\mid_{\partial_0}}\|_{L_2 (\partial_0; S_{p'}(R_p))}
\hskip1.5pt & = & \Big( \int_{\partial_0} \Big\| \Big( \summ_k
f_k(z)^* f_k(z) \Big)^{\frac12} \Big\|_{S_{p'}}^2 d\mu_0
\Big)^{\frac12}.
\end{eqnarray*}
In particular, \eqref{Eq-222} follows automatically. On the other
hand, \eqref{fbasexu} gives
\begin{eqnarray*}
\|f_{\mid_{\partial_1}}\|_{S_{p'}(L_2^{oh}(\partial_1;\ell_2))} &
= & \sup_{\|\alpha\|_{4r}, \|\beta\|_{4r} \le 1} \Big(
\int_{\partial_1} \summ_k \|\alpha f_k(z) \beta\|_{S_2}^2 d\mu_1
\Big)^{\frac12}, \\
\|f_{\mid_{\partial_1}}\|_{L_2 (\partial_1; S_{p'}(\mathrm{OH}))}
& = & \Big( \int_{\partial_1} \sup_{\|\alpha\|_{4r},
\|\beta\|_{4r} \le 1} \summ_k \|\alpha f_k(z) \beta\|_{S_2}^2
d\mu_1 \Big)^{\frac12}.
\end{eqnarray*}
Thus, inequality \eqref{Eq-333} also follows easily and
$\mathcal{Q}$ is a complete metric surjection. \fin

\begin{remark}
\emph{If $1 \le p_0 \le p \le p_1 \le \infty$, it is also true
that $$R_p \in \mathcal{QS} \big( R_{p_0} \oplus_2 R_{p_1} \big)
\quad \mbox{and} \quad C_p \in \mathcal{QS} \big( C_{p_0} \oplus_2
C_{p_1} \big).$$ The arguments to prove it are the same. However,
we have preferred to state and prove the particular case with $p_1
= 2$ for the sake of clarity, since we think of Lemma
\ref{Lemma-Motivation} as a model to follow when dealing with
non-discrete algebras.}
\end{remark}

\begin{remark} \label{Remark-Graph_Sq}
\emph{In Lemma \ref{Lemma-Motivation} we have obtained $$R_q
\simeq_{cb} \mathcal{F}(R_p, \mathrm{OH}) / \ker \mathcal{Q} \in
\mathcal{QS}(R_p \oplus_2 \mathrm{OH}),$$ $$C_q \simeq_{cb}
\mathcal{F}(C_p, \mathrm{OH}) / \ker \mathcal{Q} \in
\mathcal{QS}(C_p \oplus_2 \mathrm{OH}).$$ However, it will be
convenient in the sequel to observe that $\ker \mathcal{Q}$ can be
regarded in both cases as the annihilator of the graph of certain
linear operator, see \eqref{Pexercise}. Recall that for a given a
linear map between Hilbert spaces $\Lambda: \mathcal{K}_1 \to
\mathcal{K}_2$ with domain $\mathsf{dom}(\Lambda)$, its graph is
defined by
$$graph (\Lambda) = \Big\{ (x_1,x_2) \in \mathcal{K}_1 \oplus_2
\mathcal{K}_2 \, \big| \ x_1 \in \mathsf{dom}(\Lambda) \
\mbox{and} \ x_2 = \Lambda(x_1) \Big\}.$$ Let us consider for
instance the situation with $R_q$. We first observe that
$\mathcal{F}(R_p,\mathrm{OH})$ is the graph of an injective,
closed, densely-defined operator $\Lambda$ with dense range. This
operator is given $$\Lambda(f_{\mid_{\partial_0}}) =
f_{\mid_{\partial_1}}$$ for all $f \in \mathcal{F}(R_p,
\mathrm{OH})$. $\Lambda$ is well defined and injective by Poisson
integration due to the analyticity of elements in
$\mathcal{F}(R_p, \mathrm{OH})$. On the other hand, $\ker
\mathcal{Q}$ is the subspace of $\mathcal{F}(R_p, \mathrm{OH})$
composed of functions $f$ vanishing at $z = \theta$. Then, it
easily follows from \eqref{Eq-AnalyticCondition} that $\ker
\mathcal{Q}$ is the annihilator of $\mathcal{F}(R_{p'},
\mathrm{OH}) = graph (\Lambda)$ regarded as a subspace of
$$(1-\theta)^{\frac12} L_2^{r_{p'}}(\partial_0; \ell_2) \oplus_2
\theta^{\frac12} L_2^{oh} (\partial_1; \ell_2).$$}
\end{remark}

\section{The simplest case}
\label{NewSect2}

Given $1 < q \le 2$, we construct a completely isomorphic
embedding of $S_q$ into the predual of a $\mathrm{QWEP}$ (not yet
hyperfinite) von Neumann algebra. This is the simplest case of our
main result and will serve as a motivation for the general case.
We start with the proof of a generalized form of Theorem
\ref{MAMSth} for arbitrary von Neumann algebras, although we just
use for the moment the discrete version as stated in the
Introduction. The general formulation will be instrumental when
dealing with non-discrete algebras. In the second part of this
section, we prove Theorem \ref{QSQSQS} and deduce our cb-embedding
via Lemma \ref{Lemma-Motivation}. Theorem \ref{lp} and the
subsequent family of operator space $q$-stable random variables
arise by injecting the space $\ell_q$ into the diagonal of the
Schatten class $S_q$.

\subsection{Free harmonic analysis}
\label{P2.1}

Our starting point is a von Neumann algebra $\mathcal{M}$ equipped
with a \emph{n.f.} state $\varphi$ and associated density
$d_{\varphi}$. Let $\mathcal{N}$ be a von Neumann subalgebra of
$\mathcal{M}$. According to Takesaki \cite{Ta2}, the existence and
uniqueness of a \emph{n.f.} conditional expectation $\mathsf{E}:
\mathcal{M} \rightarrow \mathcal{N}$ is equivalent to the
invariance of $\mathcal{N}$ under the action of the modular
automorphism group $\sigma_t^{\varphi}$ associated to
$(\mathcal{M},\varphi)$. Moreover, in that case $\mathsf{E}$ is
$\varphi$-invariant and following Connes \cite{Co} we have
$\mathsf{E} \circ \sigma_t^\varphi = \sigma_t^\varphi \circ
\mathsf{E}$. In what follows, we assume these properties in all
subalgebras considered. Now we set $\mathsf{A}_k = \mathcal{M}
\oplus \mathcal{M}$ and define $\mathcal{A}$ to be the reduced
amalgamated free product $*_{\mathcal{N}} \mathsf{A}_k$ of the
family $\mathsf{A}_1, \mathsf{A}_2, \ldots, \mathsf{A}_n$ over the
subalgebra $\mathcal{N}$. Note that our notation $*_\mathcal{N}
\mathsf{A}_k$ for reduced amalgamated free products does not make
explicit the dependence on the conditional expectations
$\mathsf{E}_k: \mathsf{A}_k \to \mathcal{N}$, given by
$\mathsf{E}_k(a,b) = \frac12 \mathsf{E}(a) + \frac12
\mathsf{E}(b)$. The following is the operator-valued version
\cite{J2,JPX} of Voiculescu inequality \cite{V2}, for which we
need to introduce the mean-zero subspaces $$\bubl_k = \Big\{ x \in
\mathsf{A}_k \, \big| \ \mathsf{E}_k(a_k) = 0 \Big\}.$$

\begin{lemma} \label{Lemma-Voiculescu}
If $a_k \in \bubl_k$ for $1 \le k \le n$ and
$\mathsf{E}_\mathcal{N}: \mathcal{A} \to \mathcal{N}$ stands for
the conditional expectation of $\mathcal{A}$ onto $\mathcal{N}$,
the following equivalence of norms holds with constants
independent of $n$ $$\Big\| \sum_{k=1}^n a_k \Big\|_{\mathcal{A}}
\sim \max_{1 \le k \le n} \|x_k\|_{\mathsf{A}_k} + \Big\| \Big(
\sum_{k=1}^n \mathsf{E}_{\mathcal{N}}(a_k a_k^*) \Big)^{\frac12}
\Big\|_{\mathcal{N}} + \Big\| \Big( \sum_{k=1}^n
\mathsf{E}_{\mathcal{N}}(a_k^* a_k) \Big)^{\frac12}
\Big\|_{\mathcal{N}}.$$ Moreover, we also have
$$\Big\| \big( \sum_{k=1}^n a_k a_k^* \big)^{\frac12}
\Big\|_{\mathcal{A}} \sim \max_{1 \le k \le n}
\|a_k\|_{\mathsf{A}_k} + \Big\| \Big( \sum_{k=1}^n
\mathsf{E}_{\mathcal{N}}(a_k a_k^*) \Big)^{\frac12}
\Big\|_{\mathcal{N}},$$ $$\Big\| \big( \sum_{k=1}^n a_k^*a_k
\big)^{\frac12} \Big\|_{\mathcal{A}} \sim \max_{1 \le k \le n}
\|a_k\|_{\mathsf{A}_k} + \Big\| \Big( \sum_{k=1}^n
\mathsf{E}_{\mathcal{N}}(a_k^* a_k) \Big)^{\frac12}
\Big\|_{\mathcal{N}}.$$
\end{lemma}

\dem For the first inequality we refer to \cite{J2}. The others
can be proved in a similar way. Alternatively, both can be deduced
from the first one. Indeed, using the identity $$\Big\| \big(
\sum_{k=1}^n a_k a_k^* \big)^{\frac12} \Big\|_\mathcal{A} = \Big\|
\sum_{k=1}^n a_k \otimes e_{1k}
\Big\|_{\mathrm{M}_n(\mathcal{A})}$$ and recalling the isometric
isomorphism $$\mathrm{M}_n \big( *_{\mathcal{N}} \mathsf{A}_k
\big) = *_{\mathrm{M}_n(\mathcal{N})}
\mathrm{M}_n(\mathsf{A}_k),$$ we may apply Voiculescu's inequality
over the triple $$\Big( \mathrm{M}_n(\mathcal{A}),
\mathrm{M}_n(\mathsf{A}_k), \mathrm{M}_n(\mathcal{N}) \Big).$$
Taking $\widetilde{\mathsf{E}}_{\mathcal{N}} = id_{\mathrm{M}_n}
\ten \mathsf{E}_{\mathcal{N}}$, the last term disappears because
$$\Big\| \Big( \sum_{k=1}^n \widetilde{\mathsf{E}}_{\mathcal{N}}
\big( (a_k \ten e_{1k})^* (a_k \ten e_{1k}) \big) \Big)^{\frac12}
\Big\|_{\mathrm{M}_n(\mathcal{N})} \! = \! \sup_{1 \le k \le n}
\big\| \mathsf{E}_\mathcal{N}(a_k^* a_k)^{\frac12}
\big\|_{\mathcal{N}} \le \sup_{1 \le k \le n}
\|a_k\|_{\mathsf{A}_k}.$$ The third equivalence follows by taking
adjoints. The proof is complete. \fin

Let $\pi_k: \mathsf{A}_k \to \mathcal{A}$ denote the embedding of
$\mathsf{A}_k$ into the $k$-th component of $\mathcal{A}$. Given
$x \in \mathcal{M}$, we shall write $x_k$ as an abbreviation of
$\pi_k(x,-x)$. Note that $x_k$ is mean-zero. In the following we
shall use with no further comment the identities
$\mathsf{E}_{\mathcal{N}}(x_k x_k^*) = \mathsf{E}(x x^*)$ and
$\mathsf{E}_{\mathcal{N}}(x_k^* x_k) = \mathsf{E}(x^* x)$. We will
mostly work with identically distributed variables. In other
words, given $x \in \mathcal{M}$ we shall work with the sequence
$x_k = \pi_k(x,-x)$ for $1 \le k \le n$. In terms of the last
equivalences in Lemma \ref{Lemma-Voiculescu}, we may consider the
following norms
\begin{eqnarray*}
\|x\|_{\mathcal{R}_{\infty,1}^n} & = & \max \Big\{
\|x\|_{\mathcal{M}}, \sqrt{n} \, \big\| \mathsf{E}
(xx^*)^{\frac12} \big\|_{\mathcal{N}} \Big\}, \\
\|x\|_{\hskip1pt \mathcal{C}_{\infty,1}^n} \hskip2pt & = & \max
\Big\{ \|x\|_{\mathcal{M}}, \sqrt{n} \, \big\| \mathsf{E}
(x^*x)^{\frac12} \big\|_{\mathcal{N}} \Big\}.
\end{eqnarray*}
Here the letters $\mathcal{R}$ and $\mathcal{C}$ stand for row and
column according to Lemma \ref{Lemma-Voiculescu}. The symbol
$\infty$ is motivated because in the following section we shall
consider $L_p$ versions of these spaces. The number $1$ arises
from interpolation theory, because we think of these spaces as
endpoints in an interpolation scale. Finally, the norms on the
right induce to introduce the spaces $L_\infty^r(\mathcal{M},
\mathsf{E})$ and $L_\infty^c(\mathcal{M}, \mathsf{E})$ as the
closure of $\mathcal{M}$ with respect to the norms $$\big\|
\mathsf{E} (xx^*)^{\frac12} \big\|_{\mathcal{N}} \quad \mbox{and}
\quad \big\| \mathsf{E} (x^*x)^{\frac12} \big\|_{\mathcal{N}}.$$
In this way, we obtain the spaces
$$\begin{array}{rclcl} \mathcal{R}_{\infty,1}^n (\mathcal{M},
\mathsf{E}) & = &
\mathcal{M} \cap \sqrt{n} \, L_{\infty}^r(\mathcal{M}, \mathsf{E}), \\
[5pt] \mathcal{C}_{\infty,1}^n \, (\mathcal{M}, \mathsf{E}) & = &
\mathcal{M} \cap \sqrt{n} \, L_{\infty}^c(\mathcal{M},
\mathsf{E}). \end{array}$$

\begin{remark}
\emph{It is easily seen that
\begin{eqnarray*}
\big\| \mathsf{E} (xx^*)^{\frac12} \big\|_{\mathcal{N}} & = & \sup
\Big\{ \|\alpha x\|_{L_2(\mathcal{M})} \, \big| \
\|\alpha\|_{L_2(\mathcal{N})} \le 1 \Big\}, \\ \big\| \mathsf{E}
(x^*x)^{\frac12} \big\|_{\mathcal{N}} & = & \sup \Big\{ \|x
\beta\|_{L_2(\mathcal{M})} \, \big| \ \|\beta\|_{L_2(\mathcal{N})}
\le 1 \Big\}
\end{eqnarray*}
This relation will be crucial in this paper and will be assumed in
what follows.}
\end{remark}

The state $\varphi$ induces the \emph{n.f.} state $\phi = \varphi
\circ \mathsf{E}_\mathcal{N}$ on $\mathcal{A}$. If
$\mathcal{A}_{\oplus n}$ denotes the $n$-fold direct sum
$\mathcal{A} \oplus \mathcal{A} \oplus \ldots \oplus \mathcal{A}$,
we consider the \emph{n.f.} state $\phi_n: \mathcal{A}_{\oplus n}
\rightarrow \C$ and the conditional expectation $\mathcal{E}_n:
\mathcal{A}_{\oplus n} \rightarrow \mathcal{A}$ given by $$\phi_n
\big( \sum_{k=1}^n a_k \ten \delta_k \big) = \frac{1}{n}
\sum_{k=1}^n \phi(a_k) \qquad \mbox{and} \qquad \mathcal{E}_n
\big( \sum_{k=1}^n a_k \ten \delta_k \big) = \frac{1}{n}
\sum_{k=1}^n a_k.$$ Let us consider the map
\begin{equation} \label{Eq-Map-u}
u: x \in \mathcal{M} \mapsto \sum_{k=1}^n x_k \ten \delta_k \in
\mathcal{A}_{\oplus n} \quad \mbox{with} \quad x_k = \pi_k(x,-x).
\end{equation}

\begin{lemma} \label{Lemma-Complemented-Isomorphism}
The mappings $$\begin{array}{rrcl} u_r: & x \in
\mathcal{R}_{\infty,1}^n (\mathcal{M}, \mathsf{E}) & \mapsto &
\displaystyle \sum_{k=1}^n x_k \ten e_{1k} \in R_n(\mathcal{A}), \\
[5pt] u_c: & x \in \mathcal{C}_{\infty,1}^n \, (\mathcal{M},
\mathsf{E}) & \mapsto & \displaystyle \sum_{k=1}^n x_k \ten e_{k1}
\in C_n(\mathcal{A}),
\end{array}$$ are isomorphisms onto complemented subspaces with
constants independent of $n$.
\end{lemma}

\dem Given $x \in \mathcal{R}_{\infty,1}^n (\mathcal{M},
\mathsf{E})$, Lemma \ref{Lemma-Voiculescu} gives
$$\big\| u_r(x) \big\|_{R_n(\mathcal{A})} = \Big\| \big(
\sum_{k=1}^n x_k x_k^* \big)^{\frac12} \Big\|_{\mathcal{A}} \sim
\max_{1 \le k \le n} \|x_k\|_{\mathsf{A}_k} + \Big\| \Big(
\sum_{k=1}^n \mathsf{E}_{\mathcal{N}} \big( x_k x_k^* \big)
\Big)^{\frac12} \Big\|_{\mathcal{N}}.$$ In other words, we have
$$\big\| u_r(x) \big\|_{R_n(\mathcal{A})} \sim
\|x\|_{\mathcal{M}} + \sqrt{n} \, \|x\|_{L_{\infty}^r(\mathcal{M},
\mathsf{E})} \sim \|x\|_{\mathcal{R}_{\infty,1}^n(\mathcal{M},
\mathsf{E})}.$$ Thus $u_r$ is an isomorphism onto its image with
constants independent of $n$. The same argument yields to the same
conclusion for $u_c$. Let $d_\varphi$ and $d_\phi$ be the
densities associated to the states $\varphi$ and $\phi$. To prove
the complementation, we consider the map
$$\omega_r: x \in L_1(\mathcal{M}) + \frac{1}{\sqrt{n}} \,
L_1^r(\mathcal{M}, \mathsf{E}) \longmapsto \frac{1}{n}
\sum_{k=1}^n x_k \otimes e_{1k} \in R_1^n(L_1(\mathcal{A})),$$
where $L_1^r(\mathcal{M}, \mathsf{E})$ is the closure of
$\mathcal{N} d_\varphi \mathcal{M}$ with respect to the norm
$\|\mathsf{E}(xx^*)^{\frac12}\|_1$. Now we use approximation and
assume that $x = \alpha d_\varphi a$ for some $(\alpha,a) \in
\mathcal{N} \times \mathcal{M}$. Then, taking $a_k = \pi_k(a,-a)$
it follows from Theorem 7.1 in \cite{JX} that
\begin{eqnarray*}
\big\| \omega_r(x) \big\|_{R_1^n(L_1(\mathcal{A}))} & = &
\frac{1}{n} \, \Big\| \alpha d_\phi \Big( \sum_{k=1}^n a_k a_k^*
\Big) d_\phi \alpha^* \Big\|_{L_{1/2}(\mathcal{A})}^{1/2}
\\ & \le & \frac{1}{n} \, \Big\| \alpha d_\varphi \Big(
\sum_{k=1}^n \mathsf{E}_{\mathcal{N}} \big( a_k a_k^* \big) \Big)
d_{\varphi} \alpha^* \Big\|_{L_{1/2}(\mathcal{N})}^{1/2}.
\end{eqnarray*}
This gives $$\big\| \omega_r(x) \big\|_{R_1^n(L_1(\mathcal{A}))}
\le \frac{1}{\sqrt{n}} \, \|x\|_{L_1^r(\mathcal{M},
\mathsf{E})}.$$ On the other hand, by the triangle inequality
$$\big\| \omega_r(x) \big\|_{R_1^n(L_1(\mathcal{A}))} = \frac{1}{n}
\, \Big\| \sum_{k=1}^n x_k \otimes e_{1k}
\Big\|_{R_1^n(L_1(\mathcal{A}))} \le \|x\|_{L_1(\mathcal{M})}.$$
These estimates show that $\omega_r$ is a contraction. Note also
that $$\big\langle u_r(x), \omega_r(y) \big\rangle = \frac{1}{n}
\sum_{k=1}^n \mbox{tr}_{\mathcal{A}} \big( x_{k}^* y_{k}^{} \big)
= \frac{1}{n} \sum_{k=1}^n \mbox{tr}_{\mathcal{M}} (x^* y) =
\langle x, y \rangle.$$ In particular, since it follows from
Corollary 2.12 of \cite{J1} that
$$\mathcal{R}_{\infty,1}^n (\mathcal{M}; \mathsf{E}) = \Big(
L_1(\mathcal{M}) + \frac{1}{\sqrt{n}} \, L_1^r(\mathcal{M},
\mathsf{E}) \Big)^*,$$ it turns out that the map $\omega_r^* u_r$
is the identity on $\mathcal{R}_{\infty,1}^n (\mathcal{M};
\mathsf{E})$ and $u_r \omega_r^*$ is a bounded projection onto the
image of $u_r$ with constants independent of $n$. This completes
the proof in the row case. The column case follows in the same
way. \fin

In what follows we shall use the vector-valued space
$L_\infty(\mathcal{A}; \mathrm{OH}_n)$. This space is defined in
\cite{P2} for $\mathcal{A}$ hyperfinite, but we shall work in this
paper with $\mathcal{A}$ being a reduced free product von Neumann
algebra as defined above, which is no longer hyperfinite. There is
however a natural definition mentioned in the Introduction and
motivated by Pisier's formula $\mathrm{OH}_n = [C_n, R_n]_{1/2}$.
Indeed, recalling that the spaces $R_n(\mathcal{A})$ and
$C_n(\mathcal{A})$ are defined for any von Neumann algebra
$\mathcal{A}$, we may define
\begin{equation} \label{LpOH}
L_\infty(\mathcal{A}; \mathrm{OH}_n) = \big[ C_n(\mathcal{A}),
R_n(\mathcal{A}) \big]_{\frac12}.
\end{equation}
Pisier showed in \cite{P0} that
\begin{equation} \label{pisint}
\Big\| \sum_{k=1}^n a_k \ten \delta_k
\Big\|_{L_\infty(\mathcal{A}; \mathrm{OH}_n)} = \sup \left\{
\Big\| \sum_{k=1}^n a_k^* \alpha a_k
\Big\|_{L_2(\mathcal{A})}^{\frac12} \, \big| \ \alpha \ge 0, \,
\|\alpha\|_2 \le 1 \right\}.
\end{equation}
More general results can be found in \cite{JP2, X} or Section
\ref{NewSect3} below.

\begin{remark}
\emph{We know from \cite{J1} that the spaces
$$L_\infty(\mathcal{A}; \ell_1) \quad \mbox{and} \quad
L_\infty(\mathcal{A}; \ell_\infty)$$ are also defined for every
von Neumann algebra $\mathcal{A}$. In particular, we might wonder
whether or not our definition \eqref{LpOH} of
$L_\infty(\mathcal{A}; \mathrm{OH}_n)$ satisfies the following
complete isometry $$L_\infty(\mathcal{A}; \mathrm{OH}_n) = \big[
L_\infty(\mathcal{A}; \ell_\infty^n), L_\infty(\mathcal{A};
\ell_1^n) \big]_{\frac12}.$$ Fortunately this is the case. A
similar remark holds for $L_p(\mathcal{A}; \mathrm{OH}_n)$, see
\cite{JP2}. }
\end{remark}

The careful reader will have observed that the projection maps
$u_r w_r^*$ and $u_c w_c^*$ are the same, modulo the
identification of $R_n(\mathcal{A})$ and $C_n(\mathcal{A})$ with
$\mathcal{A}_{\oplus n}$. This is the same identification as in
Pisier's result \eqref{pisint}. In particular, this allows us to
identify via Lemma \ref{Lemma-Complemented-Isomorphism} the
interpolation space $$\mathrm{X}_{\frac12} = \big[
\mathcal{C}_{\infty,1}^n (\mathcal{M}, \mathsf{E}),
\mathcal{R}_{\infty,1}^n (\mathcal{M}, \mathsf{E})
\big]_{\frac12}$$ with a complemented subspace of
$L_\infty(\mathcal{A}; \mathrm{OH}_n)$. However, the difficult
part in proving Theorem \ref{MAMSth} is to identify the norm of
the space $\mathrm{X}_{\frac12}$. Of course, according to the fact
that we are interpolating $2$-term intersection spaces, we expect
a $4$-term maximum. This is the case and we define
$\mathcal{J}_{\infty,2}^n(\mathcal{M}, \mathsf{E})$ as the space
of elements $x$ in $\mathcal{M}$ equipped with the norm
$$\max_{u, v \in \{4,\infty\}} \left\{ n^{\frac{1}{\xi(u,v)}} \,
\sup \big\{ \|\alpha x \beta\|_{L_{\xi(u,v)}(\mathcal{M})} \, | \,
\|\alpha\|_{L_u(\mathcal{N})}, \|\beta\|_{L_v(\mathcal{N})} \le 1
\big\} \right\},$$ where $\xi(u,v)$ is given by
$\frac{1}{\xi(u,v)} = \frac1u + \frac1v$. Obviously, multiplying
by elements $\alpha, \beta$ in the unit ball of
$L_\infty(\mathcal{N})$ and taking suprema does not contribute to
the corresponding $L_{\xi(u,v)}(\mathcal{M})$ term. In other
words, we may rewrite the norm of $x$ in $\mathcal{J}_{\infty,
2}^n(\mathcal{M}, \mathsf{E})$ as $$\|x\|_{\mathcal{J}_{\infty,
2}^n(\mathcal{M}, \mathsf{E})} = \max \Big\{
\|x\|_{\Lambda_{(u,v)}^n} \, \big| \ u,v \in \{4,\infty\} \Big\}$$
where the $\Lambda_{(u,v)}^n$ norms are given by
$$\begin{array}{lcl} \|x\|_{\Lambda_{(\infty, \infty)}^n} & = &
\displaystyle \|x\|_{\mathcal{M}}, \\ [7pt]
\|x\|_{\Lambda_{(\infty, 4)}^n} & = & n^{\frac14} \, \displaystyle
\sup \Big\{ \|x \beta\|_{L_4(\mathcal{M})} \, \big|
\ \|\beta\|_{L_4(\mathcal{N})} \le 1 \Big\}, \\
[7pt] \|x\|_{\Lambda_{(4,\infty)}^n} & = & n^{\frac14} \,
\displaystyle \sup \Big\{ \|\alpha x\|_{L_4(\mathcal{M})} \, \big|
\ \|\alpha\|_{L_4(\mathcal{N})} \le 1 \Big\}, \\
[7pt]\|x\|_{\Lambda_{(4,4)}^n} & = & n^{\frac12} \, \displaystyle
\sup \Big\{ \|\alpha x \beta\|_{L_2(\mathcal{M})} \, \big| \
\|\alpha\|_{L_4(\mathcal{N})}, \|\beta\|_{L_4(\mathcal{N})} \le 1
\Big\}. \end{array}$$ These norms arise as particular cases of the
so-called conditional $L_p$ spaces, to be analyzed in Section
\ref{NewSect3}. Before identifying the norm of
$\mathrm{X}_{\frac12}$, we need some information on interpolation
spaces.

\begin{lemma} \label{Lemintinfty}
If $(1/u,1/v) = (\theta/2, (1-\theta)/2)$, we have for $x \in
\mathcal{M}$
\begin{eqnarray*}
\|x\|_{[\mathcal{M}, L_{\infty}^r(\mathcal{M},
\mathsf{E})]_{\theta}} & = & \sup \Big\{ \|\alpha x
\|_{L_u(\mathcal{M})} \, \big| \
\|\alpha\|_{L_u(\mathcal{N})} \le 1 \Big\}, \\
\|x\|_{[L_{\infty}^c(\mathcal{M}, \mathsf{E}),
\mathcal{M}]_{\theta}} & = & \sup \Big\{ \|x
\beta\|_{L_v(\mathcal{M})} \, \big| \
\|\beta\|_{L_v(\mathcal{N})} \le 1 \Big\}, \\
\|x\|_{[L_{\infty}^c(\mathcal{M}, \mathsf{E})]_{\theta},
L_{\infty}^r(\mathcal{M}, \mathsf{E})]_{\theta}} & = & \sup \Big\{
\|\alpha x \beta\|_{L_2(\mathcal{M})} \, \big| \
\|\alpha\|_{L_u(\mathcal{N})}, \|\beta\|_{L_v(\mathcal{N})} \le 1
\Big\}.
\end{eqnarray*}
\end{lemma}

The proof can be found in \cite{JP2}. In the finite setting, this
result follows from a well-known application of
Helson/Lowdenslager, Wiener/Masani type results on the existence
of operator-valued analytic functions. This kind of applications
has been used extensively by Pisier in his theory of vector-valued
$L_p$ spaces.

\begin{theorem} \label{Theorem-Intersection}
We have isomorphically
\begin{eqnarray*}
\big[ \mathcal{C}_{\infty,1}^n (\mathcal{M}, \mathsf{E}),
\mathcal{R}_{\infty,1}^n (\mathcal{M}, \mathsf{E}) \big]_{\frac12}
& \simeq & \mathcal{J}_{\infty,2}^n (\mathcal{M}, \mathsf{E}).
\end{eqnarray*}
Moreover, the constants in these isomorphisms are uniformly
bounded on $n$.
\end{theorem}

\dem We have a contractive inclusion
$$\mathrm{X}_{\frac12} \subset \mathcal{J}_{\infty,2}^n
(\mathcal{M}, \mathsf{E}).$$ Indeed, by elementary properties of
interpolation spaces we find $$\mathrm{X}_{\frac12} \subset
[\mathcal{M}, \mathcal{M}]_{\frac12} \cap [\sqrt{n} L_\infty^c,
\mathcal{M}]_{\frac12} \cap [\mathcal{M}, \sqrt{n}
L_\infty^r]_{\frac12} \cap [\sqrt{n} L_\infty^c, \sqrt{n}
L_\infty^r]_{\frac12},$$ where $L_\infty^r$ and $L_\infty^c$ are
abbreviations for $L_\infty^r(\mathcal{M}, \mathsf{E})$ and
$L_\infty^c(\mathcal{M}, \mathsf{E})$ respectively. Using the
obvious identity $[\lambda_0 \mathrm{X}_0, \lambda_1
\mathrm{X}_1]_\theta = \lambda_0^{1-\theta} \lambda_1^\theta \,
\mathrm{X}_\theta$ and applying Lemma \ref{Lemintinfty} we
rediscover the norm of the space $\mathcal{J}_{\infty,2}^n
(\mathcal{M}, \mathsf{E})$ on the right hand side. Therefore the
lower estimate holds with constant $1$.

\vskip5pt

\noindent To prove the upper estimate, we note from Lemma
\ref{Lemma-Complemented-Isomorphism} and \eqref{pisint} that
\begin{eqnarray*}
\|x\|_{\mathrm{X}_{\frac12}} & \sim & \big\| u(x)
\big\|_{[C_n(\mathcal{A}), R_n(\mathcal{A})]_{\frac{1}{2}}} \\ & =
& \Big\| \sum_{k=1}^n x_k \ten \delta_k
\Big\|_{L_\infty(\mathcal{A}; \mathrm{OH}_n)} \\ & = & \sup
\left\{ \Big\| \sum_{k=1}^n x_k^* a x_k
\Big\|_{L_2(\mathcal{A})}^{\frac12} \, \big| \ a \ge 0, \, \|a\|_2
\le 1 \right\} = \mathrm{A}.
\end{eqnarray*}
Thus, it remains to see that
\begin{equation} \label{Equation-Remaining-Estimate2}
\mathrm{A} \ \lesssim \max \Big\{ \|x\|_{\Lambda_{(u,v)}^n} \,
\big| \ u,v \in \{4,\infty\} \Big\} =
\|x\|_{\mathcal{J}_{\infty,2}^n(\mathcal{M}, \mathsf{E})}.
\end{equation}
In order to justify \eqref{Equation-Remaining-Estimate2}, we
introduce the orthogonal projections $\mathsf{L}_k$ and
$\mathsf{R}_k$ in $L_2(\mathcal{A})$ defined as follows. Given $a
\in L_2(\mathcal{A})$, the vector $\mathsf{L}_k(a)$ (resp.
$\mathsf{R}_k(a)$) collects the reduced words in $a$ starting
(resp. ending) with a letter in $\mathsf{A}_k$. In other words,
following standard terminology in free probability, we have
\begin{eqnarray*}
\mathsf{L}_k: L_2(\mathcal{A}) & \longrightarrow & L_2 \Big( \big[
\bigoplus_{m \ge 1} \bigoplus_{j_1 =k \neq j_2 \neq \cdots
\neq j_m} \bubl_{j_1} \bubl_{j_2} \cdots \bubl_{j_m} \big]'' \Big), \\
\mathsf{R}_k: L_2(\mathcal{A}) & \longrightarrow & L_2 \Big( \big[
\bigoplus_{m \ge 1} \bigoplus_{j_1 \neq j_2 \neq \cdots \neq k =
j_m} \bubl_{j_1} \bubl_{j_2} \cdots \bubl_{j_m} \big]'' \Big).
\end{eqnarray*}
Now, given a positive operator $a$ in $L_2(\mathcal{A})$ and a
fixed integer $1 \le k \le n$, we consider the following way to
decompose $a$ in terms of the projections $\mathsf{L}_k$ and
$\mathsf{R}_k$ and the conditional expectation
$\mathsf{E}_\mathcal{N}: \mathcal{A} \to \mathcal{N}$
\begin{equation} \label{Equation-Decomposition}
a = \mathsf{E}_{\mathcal{N}}(a) + \mathsf{L}_k(a) +
\mathsf{R}_k(a) - \mathsf{R}_k \mathsf{L}_k(a) + \gamma_k(a),
\end{equation} where the term $\gamma_k(a)$ has the
following form $$\gamma_k(a) = a - \mathsf{E}_{\mathcal{N}}(a) -
\mathsf{L}_k(a) - \mathsf{R}_k(a - \mathsf{L}_k(a)).$$ The
triangle inequality gives $\mathrm{A}^2 \le \sum_{k=1}^5
\mathrm{A}_j^2$, where the terms $\mathrm{A}_j$ are the result of
replacing $a$ in $\mathrm{A}$ by the $j$-th term in the
decomposition \eqref{Equation-Decomposition}. Let us estimate
these terms separately. For the first term
$\mathsf{E}_{\mathcal{N}}(a)$ we use
\begin{eqnarray*}
\mathrm{A}_1^2 & = & \Big\| \sum_{k=1}^n x_k^*
\mathsf{E}_{\mathcal{N}}(a) x_k \Big\|_{L_2(\mathcal{A})} \\
& \le & \Big\| \sum_{k=1}^n \mathsf{E}_{\mathcal{N}} \big( x_k^*
\mathsf{E}_{\mathcal{N}}(a) x_k \big) \Big\|_{L_2(\mathcal{N})} \\
& + & \Big\| \sum_{k=1}^n x_k^* \mathsf{E}_{\mathcal{N}}(a) x_k -
\mathsf{E}_{\mathcal{N}} \big( x_k^* \mathsf{E}_{\mathcal{N}}(a)
x_k \big) \Big\|_{L_2(\mathcal{A})} = \mathrm{A}_{11}^2 +
\mathrm{A}_{12}^2.
\end{eqnarray*}
Since $\mathsf{E}_{\mathcal{N}} \big( x_k^*
\mathsf{E}_{\mathcal{N}}(a) x_k \big) = \mathsf{E} \big( x^*
\mathsf{E}_{\mathcal{N}}(a) x \big)$ and $a \in
\mathsf{B}_{L_2(\mathcal{A})}^+$, we obtain
\begin{eqnarray*}
\mathrm{A}_{11} & = & n^{\frac12} \, \sup \Big\{
\mbox{tr}_{\mathcal{N}} \Big( \beta^* \mathsf{E} \big( x^*
\mathsf{E}_{\mathcal{N}}(a) x \big) \beta \Big)^{\frac12} \, \big|
\ \|\beta\|_{L_4(\mathcal{N})} \le 1 \Big\} \\ & \le & n^{\frac12}
\, \sup \Big\{ \mbox{tr}_{\mathcal{M}} \big( \beta^* x^* \alpha^*
\alpha x \beta \big)^{\frac12} \, \big| \
\|\alpha\|_{L_4(\mathcal{N})}, \|\beta\|_{L_4(\mathcal{N})} \le 1
\Big\}.
\end{eqnarray*}
This gives $\mathrm{A}_{11} \le \|x\|_{\Lambda_{(4,4)}^n} \le
\|x\|_{\mathcal{J}_{\infty,2}^n(\mathcal{M}, \mathsf{E})}$. On the
other hand, by freeness
\begin{eqnarray*}
\mathrm{A}_{12}^2 & = & \Big( \sum_{k=1}^n \big\| x_k^*
\mathsf{E}_{\mathcal{N}}(a) x_k - \mathsf{E}_{\mathcal{N}} \big(
x_k^* \mathsf{E}_{\mathcal{N}}(a) x_k \big)
\big\|_{L_2(\mathcal{A})}^2 \Big)^{\frac12} \\ & \le & 2 \, \Big(
\sum_{k=1}^n \big\| x_k^* \mathsf{E}_{\mathcal{N}}(a) x_k
\big\|_{L_2(\mathcal{A})}^2 \Big)^{\frac12} = 2 \, n^{\frac12} \,
\big\| x^* \mathsf{E}_{\mathcal{N}}(a) x
\big\|_{L_2(\mathcal{M})}.
\end{eqnarray*}
Then positivity gives
$$\mathrm{A}_{12} \le \sqrt{2} \, \|x\|_{\Lambda_{(4,\infty)}^n}
\le \sqrt{2} \, \|x\|_{\mathcal{J}_{\infty,2}^n(\mathcal{M},
\mathsf{E})}.$$ The second, third and fourth terms in
\eqref{Equation-Decomposition} satisfy
\begin{eqnarray*}
\lefteqn{\Big\| \sum_{k=1}^n x_k^* \big( \mathsf{L}_k(a) +
\mathsf{R}_k(a) + \mathsf{R}_k \mathsf{L}_k(a) \big) x_k
\Big\|_{L_2(\mathcal{A})}}
\\ & = & \sup \left\{ \sum_{k=1}^n \mbox{tr}_{\mathcal{A}} \Big( b
x_k^* \big( \mathsf{L}_k(a) + \mathsf{R}_k(a) + \mathsf{R}_k
\mathsf{L}_k(a) \big) x_k \Big) \, \big| \
\|b\|_{L_2(\mathcal{A})} \le 1 \right\} \\ & = & \sup \left\{
\sum_{k=1}^n \mbox{tr}_{\mathcal{A}} \Big( x_k b x_k^* \big(
\mathsf{L}_k(a) + \mathsf{R}_k(a) + \mathsf{R}_k \mathsf{L}_k(a)
\big) \Big) \, \big| \ \|b\|_{L_2(\mathcal{A})} \le 1 \right\}
\\ & \le & \sup_{\|b\|_{L_2(\mathcal{A})} \le 1} \Big( \sum_{k=1}^n
\big\| x_k b x_k^* \big\|_{L_2(\mathcal{A})}^2 \Big)^{\frac12}
\Big( \sum_{k=1}^n \big\| \mathsf{L}_k(a) + \mathsf{R}_k(a) +
\mathsf{R}_k \mathsf{L}_k(a) \big\|_{L_2(\mathcal{A})}^2
\Big)^{\frac12}.
\end{eqnarray*}
The second factor on the right is estimated by orthogonality
$$\Big( \sum_{k=1}^n \big\| \mathsf{L}_k(a) + \mathsf{R}_k(a) +
\mathsf{R}_k \mathsf{L}_k(a) \big\|_{L_2(\mathcal{A})}^2
\Big)^{\frac12} \le 3 \, \|a\|_{L_2(\mathcal{A})} \le 3.$$ For the
first factor, we write $b$ as a linear combination $(b_1 - b_2) +
i (b_3 - b_4)$ of four positive operators. Therefore, all these
terms are covered by the following estimate, to be proved below.

\vskip5pt

\noindent \textbf{Claim.} Given $a \in
\mathsf{B}_{L_2(\mathcal{A})}^+$, we have
\begin{equation} \label{claim}
\Big( \sum_{k=1}^n \big\| x_k a x_k^* \big\|_{L_2(\mathcal{A})}^2
\Big)^{\frac14} \lesssim \max \Big\{ \|x\|_\mathcal{M}, \,
n^{\frac14} \sup_{\|\beta\|_{L_4(\mathcal{N})} \le 1} \|x
\beta\|_{L_4(\mathcal{M})} \Big\}.
\end{equation}
Before justifying our claim, we complete the proof. It remains to
estimate the term $\mathrm{A}_5$ associated to $\gamma_k(a)$. We
first observe that $\gamma_k(a)$ is a mean-zero element of
$L_2(\mathcal{A})$ made up of reduced words not starting nor
ending with a letter in $\mathsf{A}_k$. Indeed, note that
$\mathsf{E}_{\mathcal{N}}(\gamma_k(a)) = 0$ and that we first
eliminate the words starting with a letter in $\mathsf{A}_k$ by
subtracting $\mathsf{L}_k(a)$ and, after that, we eliminate the
remaining words which end with a letter in $\mathsf{A}_k$ by
subtracting $\mathsf{R}_k(a - \mathsf{L}_k(a))$. Therefore, it
turns out that $x_1^* \gamma_1(a) x_1, x_2^* \gamma_2(a) x_2,
\ldots, x_n^* \gamma_n(a) x_n$ is a free family of random
variables. In particular, by orthogonality
\begin{eqnarray*}
\Big\| \sum_{k=1}^n x_k^* \gamma_k(a) x_k
\Big\|_{L_2(\mathcal{A})}^{\frac12} & = & \Big( \sum_{k=1}^n
\big\| x_k^* \gamma_k(a) x_k \big\|_{L_2(\mathcal{A})}^2
\Big)^{\frac14}.
\end{eqnarray*}
However, recalling that $\gamma_k(a)$ is a mean-zero element made
up of words not starting nor ending with a letter in
$\mathsf{A}_k$, the following identities hold for the conditional
expectation $\mathcal{E}_{\mathsf{A}_k}: L_2(\mathcal{A}) \to
L_2(\mathsf{A}_k)$
\begin{equation} \label{Expectation-Only}
\begin{array}{c}
\mathcal{E}_{\mathsf{A}_k} \Big( \gamma_k(a)^* \big( x_kx_k^* -
\mathsf{E}_{\mathcal{N}}(x_kx_k^*) \big) \gamma_k(a) \Big) = 0, \\
[5pt] \mathcal{E}_{\mathsf{A}_k} \Big( \gamma_k(a) \big( x_kx_k^*
- \mathsf{E}_{\mathcal{N}}(x_kx_k^*) \big) \gamma_k(a)^* \Big) =
0.
\end{array}
\end{equation}
Using property (\ref{Expectation-Only}) we find
\begin{eqnarray*}
\big\|x_k^* \gamma_k(a) x_k \big\|_{L_2(\mathcal{A})}^2 & = &
\mbox{tr}_{\mathcal{A}} \big( x_k^* \gamma_k(a)^* x_k x_k^*
\gamma_k(a) x_k \big) \\ & = & \mbox{tr}_{\mathcal{A}} \big( x_k^*
\mathcal{E}_{\mathsf{A}_k} \big( \gamma_k(a)^* x_k x_k^*
\gamma_k(a) \big) x_k \big) \\ & = & \mbox{tr}_{\mathcal{A}} \big(
x_k^* \mathcal{E}_{\mathsf{A}_k} \big( \gamma_k(a)^*
\mathsf{E}_{\mathcal{N}}(x_k x_k^*) \gamma_k(a) \big) x_k \big) \\
& = & \mbox{tr}_{\mathcal{A}} \big( \gamma_k(a) x_k x_k^*
\gamma_k(a)^* \mathsf{E}_{\mathcal{N}}(x_k x_k^*) \big) \\ & = &
\mbox{tr}_{\mathcal{A}} \big( \mathcal{E}_{\mathsf{A}_k} \big(
\gamma_k(a) x_k x_k^* \gamma_k(a)^* \big)
\mathsf{E}_{\mathcal{N}}(x_k x_k^*) \big) \\ & = &
\mbox{tr}_{\mathcal{A}} \big( \mathcal{E}_{\mathsf{A}_k} \big(
\gamma_k(a) \mathsf{E}_{\mathcal{N}}(x_k x_k^*) \gamma_k(a)^*
\big) \mathsf{E}_{\mathcal{N}}(x_k x_k^*) \big) \\ & = &
\mbox{tr}_{\mathcal{A}} \big( \mathsf{E}_{\mathcal{N}}(x_k
x_k^*)^{\frac12} \gamma_k(a) \mathsf{E}_{\mathcal{N}}(x_k x_k^*)
\gamma_k(a)^* \mathsf{E}_{\mathcal{N}}(x_k x_k^*)^{\frac12} \big).
\end{eqnarray*}
In combination with $\|\gamma_k(a)\|_2 \le 5 \, \|a\|_2$ and
H{\"o}lder's inequality, this yields
$$\big\|x_k^* \gamma_k(a) x_k \big\|_{L_2(\mathcal{A})}^2 = \big\|
\mathsf{E}_{\mathcal{N}}(x x^*)^{\frac12} \gamma_k(a)
\mathsf{E}_{\mathcal{N}}(x x^*)^{\frac12}
\big\|_{L_2(\mathcal{A})}^2 \le 25 \,
\|x\|_{L_{\infty}^r(\mathcal{M},\mathsf{E})}^4.$$ We refer to
Lemma \ref{Lemintinfty} or \cite{JP2} for the fact that
$$\|x\|_{L_\infty^r(\mathcal{M}, \mathsf{E})} \le \sup \Big\{
\|\alpha x\|_{L_4(\mathcal{M})} \, \big| \
\|\alpha\|_{L_4(\mathcal{N})} \le 1 \Big\}.$$ The inequalities
proved so far give rise to the following estimate $$\Big\|
\sum_{k=1}^n x_k^* \gamma_k(a) x_k
\Big\|_{L_2(\mathcal{A})}^{\frac12} \le \sqrt{5} \,
\|x\|_{\Lambda_{(4,\infty)}^n} \le \sqrt{5} \,
\|x\|_{\mathcal{J}_{\infty,2}^n(\mathcal{M}, \mathsf{E})}.$$
Therefore, it remains to prove the claim. We proceed in a similar
way. According to the decomposition
\eqref{Equation-Decomposition}, we may use the triangle inequality
and decompose the left hand side of \eqref{claim} into five terms
$\mathrm{B}_1, \mathrm{B}_2, \ldots, \mathrm{B}_5$. For the first
term, we deduce from positivity that $$\Big( \sum_{k=1}^n \big\|
x_k \mathsf{E}_\mathcal{N}(a) x_k^* \big\|_{L_2(\mathcal{A})}^2
\Big)^{\frac14} = n^{\frac14} \big\| x \mathsf{E}_\mathcal{N}(a)
x^* \big\|_{L_2(\mathcal{M})}^{\frac12} \le n^{\frac14}
\sup_{\|\beta\|_{L_4(\mathcal{N})} \le 1} \|x
\beta\|_{L_4(\mathcal{M})}.$$ The terms $\mathrm{B}_2,
\mathrm{B}_3$ and $\mathrm{B}_4$ satisfy
\begin{eqnarray*}
\Big( \sum_{k=1}^n \big\| x_k^* \mathsf{L}_k(a) x_k
\big\|_{L_2(\mathcal{A})}^2 \Big)^{\frac14} & \le &
\|x\|_\mathcal{M} \Big( \sum_{k=1}^n
\|\mathsf{L}_k(a)\|_{L_2(\mathcal{A})}^2 \Big)^{\frac14}, \\ \Big(
\sum_{k=1}^n \big\| x_k^* \mathsf{R}_k(a) x_k
\big\|_{L_2(\mathcal{A})}^2 \Big)^{\frac14} & \le &
\|x\|_\mathcal{M} \Big( \sum_{k=1}^n
\|\mathsf{R}_k(a)\|_{L_2(\mathcal{A})}^2 \Big)^{\frac14},
\\ \Big( \sum_{k=1}^n \big\| x_k^* \mathsf{R}_k \mathsf{L}_k (a)
x_k \big\|_{L_2(\mathcal{A})}^2 \Big)^{\frac14} & \le &
\|x\|_\mathcal{M} \Big( \sum_{k=1}^n \|\mathsf{R}_k
\mathsf{L}_k(a)\|_{L_2(\mathcal{A})}^2 \Big)^{\frac14}.
\end{eqnarray*}
Therefore, by orthogonality we have
$$\Big( \sum_{k=1}^n \big\| x_k^* \big( \mathsf{L}_k(a) +
\mathsf{R}_k(a) + \mathsf{R}_k \mathsf{L}_k(a) \big) x_k
\big\|_{L_2(\mathcal{A})}^2 \Big)^{\frac14} \le 3 \,
\|x\|_\mathcal{M}.$$ This leaves us with the term $\mathrm{B}_5$.
Arguing as above $$\Big( \sum_{k=1}^n \big\| x_k \gamma_k(a) x_k^*
\big\|_{L_2(\mathcal{A})}^2 \Big)^{\frac14} \le \sqrt{5} \
n^{\frac14} \sup_{\|\beta\|_{L_4(\mathcal{N})} \le 1} \|x
\beta\|_{L_4(\mathcal{M})}.$$ Therefore, the claim holds and the
proof is complete. \fin

\begin{remark}
\emph{The arguments in Theorem \ref{Theorem-Intersection} also
give
\begin{eqnarray*}
\|x\|_{[\mathcal{M}, \mathcal{R}_{\infty,1}^n (\mathcal{M},
\mathsf{E})]_{\frac12}} & \sim & \max \Big\{ \|x\|_\mathcal{M},
\|x\|_{\Lambda_{(4,\infty)}^n} \Big\}, \\
\|x\|_{[\mathcal{C}_{\infty,1}^n (\mathcal{M}, \mathsf{E}) \, , \,
\mathcal{M}]_{\frac12}} & \sim & \max \Big\{ \|x\|_{\mathcal{M}},
\|x\|_{\Lambda_{(\infty,4)}^n} \Big\}.
\end{eqnarray*}}
\end{remark}

Now we show how the space $\mathrm{X}_{1/2}$ is related to Theorem
\ref{MAMSth}. The idea follows from a well-known argument in which
complete boundedness arise as a particular case of amalgamation.
More precisely, if $L_2^r(\mathcal{M})$/$L_2^c(\mathcal{M})$
denote the row/column quantizations of $L_2(\mathcal{M})$ and $2
\le q \le \infty$, the row/column operator space structures on
$L_q(\mathcal{M})$ are defined as follows
\begin{equation} \label{Eq-RCLp}
\begin{array}{rcl}
L_q^r(\mathcal{M}) & = & \big[ \mathcal{M},
L_2^r(\mathcal{M})\big]_{\frac{2}{q}}, \\ [5pt] L_q^c(\mathcal{M})
& = & \big[ \mathcal{M}, L_2^c(\mathcal{M})\big]_{\frac{2}{q}}.
\end{array}
\end{equation}
The following result from \cite{JP2} is a generalized form of
\eqref{osstruct}.

\begin{lemma} \label{Lemma-Isometry-Rho}
If $\mathcal{M}_m = \mathrm{M}_m(\mathcal{M})$, we have
\begin{eqnarray*}
\big\| d_{\varphi}^{\frac{1}{4}} \big( x_{ij} \big)
\big\|_{\mathrm{M}_m(L_4^r(\mathcal{M}))} & = &
\sup_{\|\alpha\|_{S_4^m} \le 1} \Big\| d_\varphi^{\frac{1}{4}}
\Big( \sum_{k=1}^m \alpha_{ik} x_{kj} \Big)
\Big\|_{L_4(\mathcal{M}_m)}, \\ \big\| \big( x_{ij} \big)
d_{\varphi}^{\frac{1}{4}}
\big\|_{\mathrm{M}_m(L_4^c(\mathcal{M}))} & = &
\sup_{\|\beta\|_{S_4^m} \le 1} \Big\| \Big( \sum_{k=1}^m x_{ik}
\beta_{kj} \Big) d_\varphi^{\frac{1}{4}}
\Big\|_{L_4(\mathcal{M}_m)}.
\end{eqnarray*}
\end{lemma}

\noindent The proof follows from
\begin{equation} \label{Eq-Schatten-RC}
\begin{array}{rcl}
\big\| d_{\varphi}^{\frac{1}{2}} \big( x_{ij} \big)
\big\|_{\mathrm{M}_m(L_2^r(\mathcal{M}))} & = & \displaystyle
\sup_{\|\alpha\|_{S_2^m} \le 1} \Big\| d_\varphi^{\frac{1}{2}}
\Big( \sum_{k=1}^m \alpha_{ik} x_{kj} \Big)
\Big\|_{L_2(\mathcal{M}_m)}, \\ \big\| \big( x_{ij} \big)
d_{\varphi}^{\frac{1}{2}}
\big\|_{\mathrm{M}_m(L_2^c(\mathcal{M}))} & = & \displaystyle
\sup_{\|\beta\|_{S_2^m} \le 1} \Big\| \Big( \sum_{k=1}^m x_{ik}
\beta_{kj} \Big) d_\varphi^{\frac{1}{2}}
\Big\|_{L_2(\mathcal{M}_m)},
\end{array}
\end{equation}
and some complex interpolation formulas developed in \cite{JP2}.
The identity \eqref{Eq-Schatten-RC} from which we interpolate is a
well-known expression in operator space theory (see e.g. p.56 in
\cite{ER}) and will play a crucial role in the last section of
this paper. Now we define the space
$\mathcal{J}_{\infty,2}^n(\mathcal{M})$ as follows
$$\mathcal{J}_{\infty,2}^n(\mathcal{M}) = \mathcal{M} \cap
n^{\frac{1}{4}} L_4^c(\mathcal{M}) \cap n^{\frac{1}{4}}
L_4^r(\mathcal{M}) \cap n^{\frac{1}{2}} L_2(\mathcal{M}).$$ Lemma
\ref{Lemma-Isometry-Rho} determines the operator space structure
of the cross terms in $\mathcal{J}_{\infty,2}^n(\mathcal{M})$. On
the other hand, according to Pisier's fundamental identities
\eqref{fbasexu} or \eqref{pisint}, it is easily seen that we have
$$\big\| d_{\varphi}^{\frac{1}{4}} \big( x_{ij} \big)
d_{\varphi}^{\frac{1}{4}} \big\|_{\mathrm{M}_m(L_2(\mathcal{M}))}
= \sup_{\|\alpha\|_{S_4^m}, \|\beta\|_{S_4^m} \le 1} \Big\|
d_\varphi^{\frac{1}{4}} \Big( \sum_{k,l=1}^m \alpha_{ik} x_{kl}
\beta_{lj} \Big) d_\varphi^{\frac{1}{4}}
\Big\|_{L_2(\mathcal{M}_m)}.$$ In other words, the o.s.s. of
$\mathcal{J}_{\infty,2}^n(\mathcal{M})$ is described by the
isometry
\begin{equation} \label{ossJ}
\mathrm{M}_m \big( \mathcal{J}_{\infty,2}^n(\mathcal{M}) \big) =
\mathcal{J}_{\infty,2}^n(\mathcal{M}_m, \mathsf{E}_m),
\end{equation}
where $\mathcal{M}_m = \mathrm{M}_m(\mathcal{M})$ and
$\mathsf{E}_m = id_{\mathrm{M}_m} \otimes \varphi: \mathcal{M}_m
\to \mathrm{M}_m$ for $m \ge 1$. This means that the
\emph{vector-valued} spaces $\mathcal{J}_{\infty,2}^n(\mathcal{M},
\mathsf{E})$ describe the o.s.s. of the \emph{scalar-valued}
spaces $\mathcal{J}_{\infty,2}^n(\mathcal{M})$. In the result
below we prove the operator space/free analogue of a form of
Rosenthal's inequality in the limit case $p \to \infty$, see
Section \ref{NewSect3} or \cite{JP2} for more details. This result
does not have a commutative counterpart. The particular case for
$\mathcal{M} = \mathcal{B}(\ell_2^n)$ recovers Theorem
\ref{MAMSth}. Given a von Neumann algebra $\mathcal{M}$, we set as
above $\mathsf{A}_k = \mathcal{M} \oplus \mathcal{M}$.

\renewcommand{\theequation}{$\Sigma_{\infty 2}$}
\addtocounter{equation}{-1}

\begin{corollary} \label{Corollary-Sigma-Infty-2}
If $\mathcal{A}_\mathcal{N} = *_{\mathcal{N}} \mathsf{A}_k$, the
map $$u: x \in \mathcal{J}_{\infty,2}^n(\mathcal{M}, \mathsf{E})
\mapsto \sum_{k=1}^n x_k \otimes \delta_k \in
L_\infty(\mathcal{A}_\mathcal{N}; \mathrm{OH}_n)$$ is an
isomorphism with complemented image and constants independent of
$n$. In particular, replacing as usual $(\mathcal{M}, \mathcal{N},
\mathsf{E})$ by $(\mathcal{M}_m, \mathrm{M}_m, \mathsf{E}_m)$ and
replacing $\mathcal{A}_\mathcal{N}$ by the non-amalgamated algebra
$\mathcal{A}_\C = \mathsf{A}_1
* \mathsf{A}_2 * \cdots * \mathsf{A}_n$, we obtain a
cb-isomorphism with cb-complemented image and constants
independent of $n$
\begin{equation}
\sigma: x \in \mathcal{J}_{\infty,2}^n(\mathcal{M}) \mapsto
\sum_{k=1}^n x_k \otimes \delta_k \in L_\infty(\mathcal{A}_\C;
\mathrm{OH}_n).
\end{equation}
\end{corollary}

\renewcommand{\theequation}{\arabic{equation}}
\numberwithin{equation}{section}

\dem The first assertion follows from Lemma
\ref{Lemma-Complemented-Isomorphism} and Theorem
\ref{Theorem-Intersection}. To prove the second assertion we
choose the triple $(\mathcal{M}_m, \mathrm{M}_m, \mathsf{E}_m)$
and apply \eqref{ossJ}. This provides us with an isomorphic
embedding $$\sigma_m: x \in \mathrm{M}_m \big(
\mathcal{J}_{\infty,2}^n (\mathcal{M}) \big) \mapsto \sum_{k=1}^n
x_k \otimes \delta_k \in L_\infty (\mathcal{A}_m;
\mathrm{OH}_n),$$ where the von Neumann algebra $\mathcal{A}_m$ is
given by $$\mathcal{A}_m = \mathrm{M}_m(\mathcal{A}_\C) =
\mathrm{M}_m(\mathsf{A}_1) *_{\mathrm{M}_m}
\mathrm{M}_m(\mathsf{A}_2) *_{\mathrm{M}_m} \cdots
*_{\mathrm{M}_m} \mathrm{M}_m(\mathsf{A}_n).$$ The last isometry
is well-known, see e.g. \cite{J2}. In particular
$$L_\infty(\mathcal{A}_m; \mathrm{OH}_n) = \mathrm{M}_m \big(
L_\infty(\mathcal{A}_\C; \mathrm{OH}_n) \big)$$ and it turns out
that $\sigma_m = id_{\mathrm{M}_m} \otimes \sigma$. This completes
the proof. \fin

\begin{remark} \label{Remark-HalfWay}
\emph{A quick look at Corollary \ref{Corollary-Sigma-Infty-2}
shows that our formulation of $(\Sigma_{\infty 2})$ is the
half-way result (in the sense of complex interpolation) between
the stated isomorphisms in Lemma \ref{Lemma-Voiculescu}. In the
same way, as we shall see in Section \ref{NewSect3}, the free
analogue of $(\Sigma_{p2})$ is the half-way result between the row
and column formulations of the free Rosenthal inequality
\cite{JPX} for positive random variables. However, this nice
property is no longer true for $(\Sigma_{pq})$ with $q \neq 2$,
see below for details.}
\end{remark}


\subsection{Embedding $S_q$ into $L_1(\mathcal{A})$}

The tools developed so far allow us to prove Theorem \ref{QSQSQS}
and thereby obtain a complete embedding of the Schatten class
$S_q$ into $L_1(\mathcal{A})$ for some $\mathrm{QWEP}$ von Neumann
algebra $\mathcal{A}$. Our main concern here is to set a model
from which we may motivate/justify the forthcoming definitions and
arguments. We shall use some well-known facts from the theory of
operator spaces which we do not state here to simplify the
exposition. All these results will be properly stated in Section
\ref{NewSect4} and we shall refer to them. We fix $\mathcal{M} =
\mathcal{B}(\ell_2)$ and consider a family $\gamma_1, \gamma_2,
\ldots \in \R_+$ of strictly positive numbers. Then we define
$\mathsf{d}_\gamma$ to be the diagonal operator on $\ell_2$
defined by $\mathsf{d}_\gamma = \sum_k \gamma_k e_{kk}$. This
operator can be regarded as the density $d_{\psi}$ associated to a
normal strictly semifinite faithful (\emph{n.s.s.f.} in short)
weight $\psi$ on $\mathcal{B}(\ell_2)$. Let us set $q_n$ to be the
projection $\sum_{k \le n} e_{kk}$ and let us consider the
restriction of $\psi$ to the subalgebra $q_n \mathcal{B}(\ell_2)
q_n$ $$\psi_n \Big( q_n \big( \summ_{ij} x_{ij} e_{ij} \big) q_n
\Big) = \sum_{k=1}^n \gamma_k x_{kk}.$$ Note that if we set
$\mathrm{k}_n = \psi_n(q_n)$, we obtain $\psi_n = \mathrm{k}_n
\varphi_n$ for some state $\varphi_n$ on $q_n \mathcal{B}(\ell_2)
q_n$. If $d_{\psi_n}$ denotes the density on $q_n
\mathcal{B}(\ell_2) q_n$ associated to the weight $\psi_n$, we
define the space $\mathcal{J}_{\infty,2}(\psi_n)$ as the subspace
$$\Big\{ \big( z, z d_{\psi_n}^{\frac14}, d_{\psi_n}^{\frac14} z,
d_{\psi_n}^{\frac14} z d_{\psi_n}^{\frac14} \big) \, \big| \ z \in
q_n \mathcal{B}(\ell_2) q_n \Big\}$$ of the direct sum
$$\mathcal{L}_{\infty}^n = \big( C_n \otimes_h R_n \big)
\oplus_2 \big( C_n \otimes_h \mathrm{OH}_n \big) \oplus_2 \big(
\mathrm{OH}_n \otimes_h R_n \big) \oplus_2 \big( \mathrm{OH}_n
\otimes_h \mathrm{OH}_n \big).$$ In other words, we obtain the
intersection space considered in the Introduction
$$(C_n \ten_h R_n) \cap (C_n \ten_h \mathrm{OH}_n)
d_{\psi_n}^{\frac14} \cap d_{\psi_n}^{\frac14} (\mathrm{OH}_n
\ten_h R_n) \cap d_{\psi_n}^{\frac14} (\mathrm{OH}_n \ten_h
\mathrm{OH}_n) d_{\psi_n}^{\frac14}.$$

\begin{lemma} \label{Lemma-PredualMAMS11}
Let us consider $$\mathcal{K}_{1,2}(\psi_n) =
\mathcal{J}_{\infty,2}(\psi_n)^*.$$ Assume that $\mathrm{k}_n =
\sum_{k=1}^n \gamma_k$ is an integer and define $\mathcal{A}_n$ to
be the $\mathrm{k}_n$-fold reduced free product of $q_n
\mathcal{B}(\ell_2) q_n \oplus q_n \mathcal{B}(\ell_2) q_n$. If
$\pi_j: q_n \mathcal{B}(\ell_2) q_n \oplus q_n \mathcal{B}(\ell_2)
q_n \to \mathcal{A}_n$ is the natural embedding into the $j$-th
component of $\mathcal{A}_n$ and we set $x_j = \pi_j(x,-x)$, the
mapping $$\omega: x \in \mathcal{K}_{1,2}(\psi_n) \mapsto
\frac{1}{\ \mathrm{k}_n} \sum_{j=1}^{\mathrm{k}_n} x_j \otimes
\delta_j \in L_1(\mathcal{A}_n; \mathrm{OH}_{\mathrm{k}_n})$$ is a
cb-embedding with cb-complemented image and constants independent
of $n$.
\end{lemma}

\dem We claim that $$\mathcal{J}_{\infty,2}(\psi_n) =
\mathcal{J}_{\infty,2}^{\mathrm{k}_n} (q_n \mathcal{B}(\ell_2)
q_n)$$ completely isometrically. Indeed, by \eqref{Eq-RCLp}
\begin{eqnarray*}
\mathrm{k}_n^{\frac14} L_4^r(q_n \mathcal{B}(\ell_2) q_n,
\varphi_n) & = & \mathrm{k}_n^{\frac14} \, \big[
\mathcal{B}(\ell_2^n), L_2^r (\mathcal{B}(\ell_2^n), \varphi_n)
\big]_{\frac12} \\ & = & \mathrm{k}_n^{\frac14} \, \big[
\mathcal{B}(\ell_2^n), d_{\varphi_n}^{\frac12} L_2^r
(\mathcal{B}(\ell_2^n), \mathrm{tr}_n) \big]_{\frac12} \\ & = &
\mathrm{k}_n^{\frac14} d_{\varphi_n}^{\frac14} \big[ C_n \ten_h
R_n, R_n \ten_h R_n \big]_{\frac12} = d_{\psi_n}^{\frac14}
(\mathrm{OH}_n \ten_h R_n).
\end{eqnarray*}
Similarly, we can treat the other terms and obtain
\begin{eqnarray*}
\mathrm{k}_n^{\frac14} L_4^c(q_n \mathcal{B}(\ell_2) q_n,
\varphi_n) & = & (C_n \ten_h \mathrm{OH}_n) d_{\psi_n}^{\frac14},
\\ \mathrm{k}_n^{\frac12} L_2(q_n \mathcal{B}(\ell_2) q_n,
\varphi_n) & = & d_{\psi_n}^{\frac14} (\mathrm{OH}_n \ten_h
\mathrm{OH}_n) d_{\psi_n}^{\frac14}.
\end{eqnarray*}
In particular, Corollary \ref{Corollary-Sigma-Infty-2} provides a
cb-isomorphism
$$\sigma: x \in \mathcal{J}_{\infty,2}(\psi_n) \mapsto
\sum_{j=1}^{\mathrm{k}_n} x_j \otimes \delta_j \in
L_{\infty}(\mathcal{A}_n; \mathrm{OH}_{\mathrm{k}_n})$$ onto a
cb-complemented subspace with constants independent of $n$ and
$$\big\langle \sigma(x), \omega(y) \big\rangle = \frac{1}{\ \mathrm{k}_n}
\sum_{j=1}^{\mathrm{k}_n} \mbox{tr}_{\mathcal{A}_n} (x_j^*y_j) =
\mbox{tr}_n(x^*y) = \langle x,y \rangle.$$ In particular, the
stated properties of $\omega$ follow from those of the mapping
$\sigma$. \fin

Now we give a more explicit description of
$\mathcal{K}_{1,2}(\psi_n)$. Using the terminology introduced
before Lemma \ref{Lemma-PredualMAMS11}, the dual of the space
$\mathcal{L}_\infty^n$ is given by the following direct sum
$$\mathcal{L}_1^n = \big( R_n \otimes_h C_n \big)
\oplus_2 \big( R_n \otimes_h \mathrm{OH}_n \big) \oplus_2 \big(
\mathrm{OH}_n \otimes_h C_n \big) \oplus_2 \big( \mathrm{OH}_n
\otimes_h \mathrm{OH}_n \big).$$ Thus, we may consider the map
$$\Psi_n: \mathcal{L}_1^n \to L_1(q_n \mathcal{B}(\ell_2) q_n)$$
given by
$$\Psi_n (x_1, x_2, x_3,x_4) = x_1 + x_2
d_{\psi_n}^{\frac{1}{4}} + d_{\psi_n}^{\frac{1}{4}} x_3 +
d_{\psi_n}^{\frac{1}{4}} x_4 d_{\psi_n}^{\frac{1}{4}}.$$ Then it
is easily checked that $\ker \Psi_n =
\mathcal{J}_{\infty,2}(\psi_n)^{\perp}$ with respect to the
anti-linear duality bracket and we deduce
$\mathcal{K}_{1,2}(\psi_n) = \mathcal{L}_1^n / \ker \Psi_n$. The
finite-dimensional spaces defined so far allow us to take direct
limits
$$\mathcal{J}_{\infty,2} (\psi) = \overline{\bigcup_{n \ge 1}
\mathcal{J}_{\infty,2}(\psi_n)} \quad \mbox{and} \quad
\mathcal{K}_{1,2}(\psi) = \overline{\bigcup_{n \ge 1}
\mathcal{K}_{1,2}(\psi_n)}.$$

\begin{lemma} \label{Lemma-Direct-Sum11}
Let $\lambda_1, \lambda_2, \ldots \in \R_+$ be a sequence of
strictly positive numbers and define the diagonal operator
$\mathsf{d}_\lambda = \sum_k \lambda_k e_{kk}$ on $\ell_2$. Let us
equip the space $graph(\mathsf{d}_\lambda)$ with the following
operator space structures $$\begin{array}{rclcl} R \cap
\ell_2^{oh}(\lambda) & = & graph(\mathsf{d}_\lambda) & \subset & R
\oplus_2 \mathrm{OH}, \\ [3pt] C \cap \ell_2^{oh}(\lambda) & = &
graph(\mathsf{d}_\lambda) & \subset & C \hskip0.5pt \oplus_2
\mathrm{OH}. \end{array}$$ Then, if we consider the dual spaces
$$\begin{array}{rclcl} C + \ell_2^{oh}(\lambda) & = & \big( C
\oplus_2 \mathrm{OH} \big) \big/ \big( R \cap \ell_2^{oh}(\lambda)
\big)^\perp,
\\ [3pt] R + \ell_2^{oh}(\lambda) & = & \big( R \oplus_2
\mathrm{OH} \big) \big/ \big( C \cap \ell_2^{oh}(\lambda)
\big)^\perp, \end{array}$$ there exists a n.s.s.f. weight $\psi$
on $\mathcal{B}(\ell_2)$ such that
$$\big( R + \ell_2^{oh}(\lambda) \big) \otimes_h
\big( C + \ell_2^{oh}(\lambda) \big) = \mathcal{K}_{1,2}(\psi).$$
\end{lemma}

\dem If we set $$q_n = \sum_{k \le n} e_{kk},$$ then we define
$$\begin{array}{lclcl} q_n \big( C + \ell_2^{oh}(\lambda) \big)
\!\!\!\! & = & \!\!\!\! \Big\{ q_n(a,b) + \big( R \cap
\ell_2^{oh}(\lambda) \big)^\perp \, \big| \ (a,b) \in C \oplus_2
\mathrm{OH} \Big\} \!\!\!\! & \subset & \!\!\!\! C +
\ell_2^{oh}(\lambda), \\ [5pt] q_n \big( R + \ell_2^{oh}(\lambda)
\big) \!\!\!\! & = & \!\!\!\! \Big\{ q_n(a,b) + \big( C \cap
\ell_2^{oh}(\lambda) \big)^\perp \, \big| \ (a,b) \in R \oplus_2
\mathrm{OH} \Big\} \!\!\!\! & \subset & \!\!\!\! R +
\ell_2^{oh}(\lambda).
\end{array}$$ Note that, since the corresponding annihilators are
$q_n$-invariant, these are quotients of $C_n \oplus_2
\mathrm{OH}_n$ and $R_n \oplus_n \mathrm{OH}_n$ respectively.
Moreover, recalling that $q_n(x) \to x$ as $n \to \infty$ in the
norms of $R, \mathrm{OH}, C$, it is not difficult to see that we
may write the Haagerup tensor product $\big( R +
\ell_2^{oh}(\lambda) \big) \otimes_h \big( C +
\ell_2^{oh}(\lambda) \big)$ as the direct limit
$$\overline{\bigcup_{n \ge 1} q_n(R + \ell_2^{oh}(\lambda))
\otimes_h q_n(C + \ell_2^{oh}(\lambda))}.$$ Therefore, it suffices
to show that $$q_n(R + \ell_2^{oh}(\lambda)) \otimes_h q_n(C +
\ell_2^{oh}(\lambda)) = \mathcal{K}_{1,2}(\psi_n),$$ where
$\psi_n$ denotes the restriction to $q_n \mathcal{B}(\ell_2) q_n$
of some \emph{n.s.s.f.} weight $\psi$. However, by duality this is
equivalent to see that $q_n(C \cap \ell_2^{oh}(\lambda)) \otimes_h
q_n(R \cap \ell_2^{oh}(\lambda)) = \mathcal{J}_{\infty,2}(\psi_n)$
where the spaces $q_n(R \cap \ell_2^{oh}(\lambda))$/$q_n(C \cap
\ell_2^{oh}(\lambda))$ are the span of $$\Big\{ (\delta_k,
\lambda_k \delta_k) \, \big| \ 1 \le k \le n \Big\}$$ in $R_n
\oplus_2 \mathrm{OH}_n$/$C_n \oplus_2 \mathrm{OH}_n$ respectively.
Indeed, we have
\begin{eqnarray*}
q_n \big( C + \ell_2^{oh}(\lambda) \big) & = & \big( C_n \oplus_2
\mathrm{OH}_n \big) / q_n(R \cap \ell_2^{oh}(\lambda))^{\perp},
\\ q_n \big( R + \ell_2^{oh}(\lambda) \big) & = & \big(
R_n \oplus_2 \mathrm{OH}_n \big) / q_n(C \cap
\ell_2^{oh}(\lambda))^{\perp},
\end{eqnarray*}
completely isometrically. Using row/column terminology in terms of
matrix units
\begin{eqnarray*}
q_n \big( C \cap \ell_2^{oh}(\lambda) \big) & = & \mbox{span}
\Big\{ (e_{i1} \hskip0.5pt , \hskip0.5pt \lambda_i \hskip0.5pt
e_{i1} \hskip0.5pt ) \in C_n \oplus_2
\hskip1pt \mathrm{OH}_n \Big\}, \\
q_n \big( R \cap \ell_2^{oh}(\lambda) \big) & = & \mbox{span}
\Big\{ (e_{1j}, \lambda_j e_{1j}) \in R_n \oplus_2 \mathrm{OH}_n
\Big\}.
\end{eqnarray*}
Therefore, the space $q_n(C \cap \ell_2^{oh}(\lambda)) \otimes_h
q_n(R \cap \ell_2^{oh}(\lambda))$ is the subspace
$$\mbox{span} \Big\{ (e_{ij}, \lambda_j e_{ij}, \lambda_i e_{ij},
\lambda_i \lambda_j e_{ij} ) \Big\} = \Big\{ (z, z
\mathsf{d}_{\lambda}, \mathsf{d}_\lambda z, \mathsf{d}_\lambda z
\mathsf{d}_\lambda) \, \big| \ z \in q_n \mathcal{B}(\ell_2) q_n
\Big\}$$ of the space $\mathcal{L}_\infty^n$ defined above. Then,
we define $\gamma_k \in \R_+$ by the relation $\lambda_k =
\gamma_k^{\frac14}$ and consider the \emph{n.s.s.f.} weight $\psi$
on $\mathcal{B}(\ell_2)$ induced by $\mathsf{d}_\gamma$. In
particular, we immediate obtain
$$q_n \big( C \cap \ell_2^{oh}(\lambda) \big) \otimes_h
q_n \big( R \cap \ell_2^{oh}(\lambda) \big) = \Big\{ \big( z, z
d_{\psi_n}^{\frac{1}{4}}, d_{\psi_n}^{\frac{1}{4}} z,
d_{\psi_n}^{\frac{1}{4}} z d_{\psi_n}^{\frac{1}{4}} \big)
\Big\}.$$ The space on the right is by definition
$\mathcal{J}_{\infty,2}(\psi_n)$. This completes the proof. \fin

\demC By injectivity of the Haagerup tensor product, we may assume
that $(\mathrm{X}_1, \mathrm{X}_2) \in \mathcal{Q}(R \oplus_2
\mathrm{OH}) \times \mathcal{Q}(C \oplus_2 \mathrm{OH})$. In
particular, the duals $\mathrm{X}_1^*$ and $\mathrm{X}_2^*$ are
subspaces of $C \oplus_2 \mathrm{OH}$ and $R \oplus_2 \mathrm{OH}$
respectively. Therefore, using a well-known result (see Lemma
\ref{Lemma-Xu-Partition} below), we may find Hilbert spaces
$\mathcal{H}_{ij}$ and $\mathcal{K}_{ij}$ for $i,j=1,2$ such that
\begin{eqnarray*}
\mathrm{X}_1^* & \simeq_{cb} & \mathcal{H}_{11,c} \oplus_2
\mathcal{H}_{12,oh} \oplus_2 graph(\Lambda_1), \\
\mathrm{X}_2^* & \simeq_{cb} & \mathcal{H}_{21,r} \oplus_2
\mathcal{H}_{22,oh} \oplus_2 graph(\Lambda_2),
\end{eqnarray*}
where the operators $\Lambda_1: \mathcal{K}_{11,c} \to
\mathcal{K}_{12,oh}$ and $\Lambda_2: \mathcal{K}_{21,r} \to
\mathcal{K}_{22,oh}$ are injective, closed, densely-defined with
dense range. On the other hand, using the complete isometries
$\mathcal{H}_r^* = \mathcal{H}_c$ and $\mathcal{H}_c^* =
\mathcal{H}_r$, we easily obtain the cb-isomorphisms
\begin{eqnarray*}
\mathrm{X}_1 & \simeq_{cb} & \mathcal{H}_{11,r} \oplus_2
\mathcal{H}_{12,oh} \oplus_2 \Big( \big( \mathcal{K}_{11,r}
\oplus_2 \mathcal{K}_{12,oh} \big) \big/ graph(\Lambda_1)^\perp
\Big), \\ \mathrm{X}_2 & \simeq_{cb} & \mathcal{H}_{21,c} \oplus_2
\mathcal{H}_{22,oh} \oplus_2 \Big( \big( \mathcal{K}_{21,c}
\oplus_2 \mathcal{K}_{22,oh} \big) \big/ graph(\Lambda_2)^\perp
\Big).
\end{eqnarray*}
Let us set for the sequel
\begin{eqnarray*}
\mathcal{Z}_1 & = & \big( \mathcal{K}_{11,r} \oplus_2
\mathcal{K}_{12,oh} \big) \big/ graph(\Lambda_1)^\perp, \\
\mathcal{Z}_2 & = & \big( \mathcal{K}_{21,c} \oplus_2
\mathcal{K}_{22,oh} \big) \big/ graph(\Lambda_2)^\perp.
\end{eqnarray*}
Then, we have the following cb-isometric inclusion
\begin{eqnarray} \label{Eq-6terms11}
\mathrm{X}_1 \otimes_h \mathrm{X}_2 \subset \mathrm{A}_1 \oplus_2
\mathrm{A}_2 \oplus_2 \mathrm{A}_3 \oplus_2 \mathrm{A}_4 \oplus_2
\mathrm{A}_5 \oplus_2 \mathrm{A}_6,
\end{eqnarray}
where the $\mathrm{A}_j$'s are given by
\begin{eqnarray*}
\mathrm{A}_1 & = & \mathcal{Z}_1 \otimes_h \mathcal{Z}_2 \\
\mathrm{A}_2 & = & \mathcal{H}_{11,r} \otimes_h \mathrm{X}_2
\\ \mathrm{A}_3 & = & \mathrm{X}_1 \otimes_h \mathcal{H}_{21,c} \\
\mathrm{A}_4 & = & \mathcal{H}_{12,oh} \otimes_h \mathcal{Z}_2
\\ \mathrm{A}_5 & = & \mathcal{Z}_1 \otimes_h
\mathcal{H}_{22,oh} \\ \mathrm{A}_6 & = & \mathcal{H}_{12,oh}
\otimes_h \mathcal{H}_{22,oh}.
\end{eqnarray*}
Let us show that the proof can be reduced to the construction of a
cb-embedding $\mathcal{Z}_1 \otimes_h \mathcal{Z}_2 \to
L_1(\mathcal{A})$ for some $\mathrm{QWEP}$ von Neumann algebra
$\mathcal{A}$. Indeed, according to \cite{J2} we know that
$\mathrm{OH}$ cb-embeds in $L_1(\mathcal{A})$ for some
$\mathrm{QWEP}$ type $\mathrm{III}$ factor $\mathcal{A}$. Hence,
the last term on the right of \eqref{Eq-6terms11} automatically
satisfies the assertion. A similar argument works for the second
and third terms. Indeed, they clearly embed into
$S_1(\mathrm{X}_1)$ and $S_1(\mathrm{X}_2)$ completely
isometrically. On the other hand, since $\mathrm{OH} \in
\mathcal{QS}(C \oplus R)$ by \lq\lq Pisier's exercise\rq\rq${}$
and we have by hypothesis
$$\mathrm{X}_1 \in \mathcal{QS}(R \oplus_2 \mathrm{OH}) \quad
\mbox{and} \quad \mathrm{X}_2 \in \mathcal{QS}(C \oplus_2
\mathrm{OH}),$$ both $\mathrm{X}_1$ and $\mathrm{X}_2$ are
cb-isomorphic to an element in $\mathcal{QS}(C \oplus R)$.
According to \cite{J2} one more time, we know that any operator
space in $\mathcal{QS}(C \oplus R)$ cb-embeds into
$L_1(\mathcal{A})$ for some $\mathrm{QWEP}$ von Neumann algebra
$\mathcal{A}$. Thus, the spaces $S_1(\mathrm{X}_1)$ and
$S_1(\mathrm{X}_2)$ also satisfy the assertion. Finally, for the
fourth and fifth terms on the right of \eqref{Eq-6terms11}, we may
write $\mathrm{OH}$ as the graph of a diagonal operator on
$\ell_2$, see Lemma \ref{Lemma-OH-graph} below for further
details. In particular, by the self-duality of $\mathrm{OH}$ we
conclude that these terms can be regarded as particular cases of
the first term $\mathcal{Z}_1 \otimes_h \mathcal{Z}_2$. It remains
to see that the term $\mathcal{Z}_1 \otimes_h \mathcal{Z}_2$
satisfies the assertion. By discretization (see Lemma
\ref{Lemma-Diagonalization}) we may assume that the graphs
appearing in the terms $\mathcal{Z}_1$ and $\mathcal{Z}_2$ above
are graphs of diagonal operators $\mathsf{d}_{\lambda_1}$ and
$\mathsf{d}_{\lambda_2}$. Moreover, using polar decomposition we
may also assume that both diagonal operators are positive, see the
proof of Lemma \ref{Lemma-Diagonalization} one more time. In fact,
by adding a perturbation term we can take the eigenvalues
$\lambda_{1k}, \lambda_{2k} \in \mathbb{R}_+$ strictly positive.
Indeed, if we replace $\lambda_{jk}$ by $\xi_{jk} = \lambda_{jk} +
\varepsilon_k$ for $j=1,2$, the new diagonal operators
$\mathsf{d}_{\xi_1}$ and $\mathsf{d}_{\xi_2}$ satisfy the
cb-isomorphisms
$$graph(\mathsf{d}_{\lambda_j}) \simeq_{cb}
graph(\mathsf{d}_{\xi_j}) \quad \mbox{for} \quad j=1,2$$ where
(arguing as in Lemma \ref{Lemma-OH-graph} below) the cb-norms are
controlled by
$$\Big( \summ_k |\varepsilon_k|^4 \Big)^{\frac14}.$$
Therefore, taking the $\varepsilon_k$'s small enough, we may write
\begin{eqnarray*}
\mathcal{Z}_1 = \big( R \oplus_2 \mathrm{OH} \big) \big/ \big( C
\cap \ell_2^{oh}(\lambda_1) \big)^{\perp} = R +
\ell_2^{oh}(\lambda_1) & \mbox{with} &
\mathsf{d}_{\lambda_1}: C \to \mathrm{OH}, \\
\mathcal{Z}_2 = \big( C \oplus_2 \mathrm{OH} \big) \big/ \big( R
\cap \ell_2^{oh}(\lambda_2) \big)^{\perp} = C +
\ell_2^{oh}(\lambda_2) & \mbox{with} & \mathsf{d}_{\lambda_2}: R
\to \mathrm{OH},
\end{eqnarray*}
where the diagonal operators above are positive and invertible.
Now we define $$\lambda_k =
\begin{cases} \lambda_{1, \frac{k+1}{2}} & \mbox{if $k$ is odd},
\\ \lambda_{2, \, \frac{k}{2}} & \mbox{if $k$ is even}. \end{cases}$$
This defines a positive invertible operator $\mathsf{d}_{\lambda}$
such that $$\mathcal{Z}_1 \otimes_h \mathcal{Z}_2 \subset \big( R
+ \ell_2^{oh}(\lambda) \big) \ten_h \big( C + \ell_2^{oh}(\lambda)
\big),$$ where the former is clearly cb-complemented in the
latter. According to Lemma \ref{Lemma-Direct-Sum11}, we conclude
that $\mathcal{Z}_1 \otimes_h \mathcal{Z}_2$ can be regarded as a
completely complemented subspace of the direct limit
$$\mathcal{K}_{1,2}(\psi) = \overline{\bigcup_{n \ge 1}
\mathcal{K}_{1,2}(\psi_n)}^{\null},$$ for some \emph{n.s.s.f.}
weight $\psi$ on $\mathcal{B}(\ell_2)$. It remains to construct a
completely isomorphic embedding from $\mathcal{K}_{1,2}(\psi)$
into $L_1(\mathcal{A})$ for some $\mathrm{QWEP}$ algebra
$\mathcal{A}$. To that aim, letting the constants is such
cb-embedding a little perturbation, we may assume without loss of
generality that the numbers $\mathrm{k}_n = \psi_n(q_n)$ are
non-decreasing positive integers since we may approximate each
$\mathrm{k}_n$ to its closest integer. This will allow us to apply
Lemma \ref{Lemma-PredualMAMS11} below. Now, in order to cb-embed
$\mathcal{K}_{1,2}(\psi)$ into $L_1(\mathcal{A})$, it suffices to
construct a cb-embedding of $\mathcal{K}_{1,2}(\psi_n)$ into
$L_1(\mathcal{A}_n')$ for some $\mathcal{A}_n'$ being
$\mathrm{QWEP}$ and with relevant constants independent of $n$.
Indeed, if so we may consider an ultrafilter $\mathcal{U}$
containing all the intervals $(n,\infty)$, so that we have a
completely isometric embedding
$$\mathcal{K}_{1,2}(\psi) = \overline{\bigcup_{n \ge 1}
\mathcal{K}_{1,2}(\psi_n)}^{\null} \to \prodd_{n, \mathcal{U}}
\mathcal{K}_{1,2}(\psi_n).$$ Then, according to \cite{Ra}, our
assumption provides a cb-embedding
$$\mathcal{K}_{1,2}(\psi) \to L_1(\mathcal{A}) \quad
\mbox{with} \quad \mathcal{A} = \Big( \prodd_{n, \mathcal{U}}
{\mathcal{A}'_n}_* \Big)^*.$$ Moreover, we know from \cite{J5}
that $\mathcal{A}$ is $\mathrm{QWEP}$ provided the
$\mathcal{A}_n'$'s are. Therefore, it remains to construct the
cb-embeddings $\mathcal{K}_{1,2}(\psi_n) \to L_1(\mathcal{A}_n')$.
This follows from the cb-embedding \cite{J2} of $\mathrm{OH}$ into
$L_1(\mathcal{B})$ for some $\mathrm{QWEP}$ type $\mathrm{III}$
factor $\mathcal{B}$ and from Lemma \ref{Lemma-PredualMAMS11}
$$\mathcal{K}_{1,2}(\psi_n) \to L_1(\mathcal{A}_n;
\mathrm{OH}_{\mathrm{k}_n}) \to L_1(\mathcal{A}_n \bar\otimes
\mathcal{B}) = L_p(\mathcal{A}_n').$$ Let us show that
$\mathcal{A}_n'$ is $\mathrm{QWEP}$. The algebra $\mathcal{A}_n$
is the free product of $\mathrm{k}_n$ copies of $\mathrm{M}_n
\oplus \mathrm{M}_n$. Therefore, since we know after \cite{J2} and
\cite{J5} that the $\mathrm{QWEP}$ is stable under free products
and tensor products, $\mathcal{A}_n' = \mathcal{A}_n \otimes
\mathcal{B}$ is $\mathrm{QWEP}$. \fin

\begin{corollary} \label{SqQWEP11}
$S_q$ cb-embeds into $L_1(\mathcal{A})$ for some $\mathrm{QWEP}$
algebra $\mathcal{A}$.
\end{corollary}

\dem Using the complete isometry $$S_q = C_q \otimes_h R_q,$$ the
assertion follows combining Lemma \ref{Lemma-Motivation} and
Theorem \ref{QSQSQS}. \fin

\section{Mixed norms of free variables} \label{NewSect3}

In this section we present a variation of the free Rosenthal
inequality \cite{JPX}. This is the main result of \cite{JP2} and
will be a key point to prove the complete embedding of $L_q$ into
$L_p$ in the general case. Its statement forces us to introduce
some new classes of noncommutative function spaces. The motivation
comes from our construction in the previous section and some
classical probabilistic estimates.

\subsection{Conditional $L_p$ spaces}

Inspired by Pisier's theory \cite{P2} of noncommutative
vector-valued $L_p$ spaces, several noncommutative function spaces
have been recently introduced in quantum probability. The first
insight came from some of Pisier's fundamental equalities, which
we briefly review. Let $\mathcal{N}_1$ and $\mathcal{N}_2$ be two
hyperfinite von Neumann algebras. Given $1 \le p,q \le \infty$, we
define $1/r = |1/p - 1/q|$. If $p \le q$, the norm of $x$ in $L_p
(\mathcal{N}_1; L_q(\mathcal{N}_2))$ is given by
\begin{equation} \label{AEq-p<q}
\inf \Big\{ \|\alpha\|_{L_{2r}(\mathcal{N}_1)}
\|y\|_{L_q(\mathcal{N}_1 \bar\otimes \mathcal{N}_2)}
\|\beta\|_{L_{2r}(\mathcal{N}_1)} \, \big| \ x = \alpha y \beta
\Big\}.
\end{equation}
If $p \ge q$, the norm of $x \in L_p (\mathcal{N}_1;
L_q(\mathcal{N}_2))$ is given by
\begin{equation} \label{AEq-p>q}
\sup \Big\{ \|\alpha x \beta \|_{L_q(\mathcal{N}_1 \bar\otimes
\mathcal{N}_2)} \, \big| \ \alpha, \beta \in
\mathsf{B}_{L_{2r}(\mathcal{N}_1)} \Big\}.
\end{equation}

The hyperfiniteness is an essential assumption in \cite{P2}.
However, when dealing with mixed $L_p(L_q)$ norms, Pisier's
identities remain true for general von Neumann algebras, see
\cite{JX4}. On the other hand, the \emph{row} and \emph{column}
subspaces of $L_p$ are defined as follows
$$R_p^n(L_p(\mathcal{M})) = \Big\{ \sum_{k=1}^n x_k \ten e_{1k} \,
\big| \ x_k \in L_p(\mathcal{M}) \Big\} \subset L_p \big(
\mathcal{M} \bar\ten \mathcal{B}(\ell_2) \big),$$
$$C_p^n(L_p(\mathcal{M})) = \Big\{ \sum_{k=1}^n x_k \ten
e_{k1} \, \big| \ x_k \in L_p(\mathcal{M}) \Big\} \subset L_p
\big( \mathcal{M} \bar\ten \mathcal{B}(\ell_2) \big),$$ where
$(e_{ij})$ denotes the unit vector basis of $\mathcal{B}(\ell_2)$.
These spaces are crucial in the noncommutative
Khintchine/Rosenthal type inequalities \cite{JPX,LuP,PP} and in
noncommutative martingale inequalities \cite{JX,PR,PX1}, where the
row and column spaces are traditionally denoted by
$L_p(\mathcal{M}; \ell_2^r)$ and $L_p(\mathcal{M}; \ell_2^c)$. The
norm in these spaces is given by $$\begin{array}{c} \displaystyle
\Big\| \sum_{k=1}^n x_k \otimes e_{1k}
\Big\|_{R_p^n(L_p(\mathcal{M}))} = \Big\| \Big( \sum_{k=1}^n x_k
x_k^* \Big)^{\frac12} \Big\|_{L_p(\mathcal{M})}, \\ [8pt]
\displaystyle \Big\| \sum_{k=1}^n x_k \otimes e_{k1}
\Big\|_{C_p^n(L_p(\mathcal{M}))} = \Big\| \Big( \sum_{k=1}^n x_k^*
x_k \Big)^{\frac12} \Big\|_{L_p(\mathcal{M})}.
\end{array}$$

\begin{remark}
\emph{In what follows we shall write
$$R_p^n(L_p(\mathcal{M})) = L_p(\mathcal{M}; R_p^n) \quad
\mbox{and} \quad C_p^n(L_p(\mathcal{M})) = L_p(\mathcal{M};
C_p^n).$$}
\end{remark}

\vskip5pt

Now, let us assume that $\mathcal{N}$ is a von Neumann subalgebra
of $\mathcal{M}$ and that there exists a \emph{n.f.} conditional
expectation $\mathsf{E}: \mathcal{M} \to \mathcal{N}$. Then we may
define $L_p$ norms of the \emph{conditional square functions}
$$\Big( \sum_{k=1}^n \mathsf{E}(x_k x_k^*) \Big)^{\frac12} \quad
\mbox{and} \quad \Big( \sum_{k=1}^n \mathsf{E}(x_k^* x_k)
\Big)^{\frac12}.$$ The expressions $\mathsf{E}(x_k x_k^*)$ and
$\mathsf{E}(x_k^* x_k)$ have to be defined properly for $1 \le p
\le 2$, see \cite{J1} or Chapter 1 of \cite{JP2}. Note that the
resulting spaces coincide with the row and column spaces defined
above when $\mathcal{N}$ is $\mathcal{M}$ itself. When $n=1$ we
recover the spaces $L_p^r(\mathcal{M}, \mathsf{E})$ and
$L_p^c(\mathcal{M}, \mathsf{E})$, which have been instrumental in
proving Doob's inequality \cite{J1}, see also \cite{JX2} for more
applications. In particular, taking $\mathcal{M}_{\oplus n}$ to be
the $n$-fold direct sum $\mathcal{M} \oplus \mathcal{M} \oplus
\cdots \oplus \mathcal{M}$ and considering the conditional
expectation
$$\mathcal{E}_n: \sum_{k=1}^n x_k \ten \delta_k \in
\mathcal{M}_{\oplus n} \mapsto \frac1n \sum_{k=1}^n x_k \in
\mathcal{M},$$ we easily obtain the following isometric
isomorphisms
\begin{equation} \label{AEq-RpCp-Cond}
\begin{array}{c}
L_p(\mathcal{M}; R_p^n) = \sqrt{n} \,
L_p^r \big( \mathcal{M}_{\oplus n}, \mathcal{E}_n \big), \\
[8pt] L_p(\mathcal{M}; C_p^n) = \sqrt{n} \, L_p^c \big(
\mathcal{M}_{\oplus n}, \mathcal{E}_n \big).
\end{array}
\end{equation}

We have already introduced $L_p(L_q)$ spaces, row and column
subspaces of $L_p$ and some variations associated to a given
conditional expectation. The careful reader may have noticed that
many of these norms have been used in Section \ref{NewSect2}. Now
we present a unified approach for this kind of spaces. All the
norms considered so far fit into more general families of
noncommutative function spaces, which we now define. Let us
consider the solid $\mathsf{K}$ in $\R^3$ defined by
$$\mathsf{K} = \Big\{(1/u,1/v,1/q) \, \big| \ 2 \le u,v \le
\infty, \ 1 \le q \le \infty, \ 1/u + 1/q + 1/v \le 1 \Big\}.$$
Let $\mathcal{M}$ be a von Neumann algebra equipped with a
\emph{n.f.} state $\varphi$ and let $\mathcal{N}$ be a given von
Neumann subalgebra with corresponding conditional expectation
$\mathsf{E}$. The amalgamated and conditional $L_p$ spaces are
defined as follows.

\vskip5pt

\begin{itemize}
\item[\textbf{(i)}] Let $(1/u,1/v,1/q) \in \mathsf{K}$ and take
$1/p = 1/u + 1/q + 1/v$. Then we define the \emph{amalgamated
$L_p$ space} associated to the indices $(u,q,v)$ as the subspace
$L_u(\mathcal{N}) L_q(\mathcal{M}) L_v(\mathcal{N})$ of
$L_p(\mathcal{M})$. The norm is given by
$$\inf \Big\{ \|\alpha\|_{L_u(\mathcal{N})}
\|y\|_{L_q(\mathcal{M})} \|\beta\|_{L_v(\mathcal{N})} \, \big| \ x
= \alpha y \beta \Big\}.$$

\item[\textbf{(ii)}] Let $(1/u,1/v,1/p) \in \mathsf{K}$ and take
$1/s = 1/u+1/p+1/v$. Then we define the \emph{conditional $L_p$
space} associated to the indices $(u,v)$ as the completion of
$L_p(\mathcal{M})$ with respect to the following norm
$$\sup \Big\{ \|axb\|_{L_s(\mathcal{M})} \, \big| \
\|a\|_{L_u(\mathcal{N})}, \|b\|_{L_v(\mathcal{N})} \le 1 \Big\}.$$
This space will be denoted by $$L_{(u,v)}^p(\mathcal{M},
\mathsf{E}).$$
\end{itemize}
The reader is referred to \cite{JP2} for a much more detailed
exposition of these spaces. In the following, it will also be
useful to recognize some important spaces in the terminology just
introduced. Here are the basic examples of amalgamated and
conditional $L_p$ spaces. The non-trivial isometric identities
below, which will be used in the following with no further
reference, are proved in \cite{JP2}.

\begin{itemize}
\item[\textbf{(a)}] The spaces $L_p(\mathcal{M})$ satisfy
$$\quad L_p(\mathcal{M}) = L_{\infty}(\mathcal{N}) L_p(\mathcal{M})
L_{\infty}(\mathcal{N}) \quad \mbox{and} \quad L_p(\mathcal{M}) =
L_{(\infty,\infty)}^p(\mathcal{M}, \mathsf{E}).$$

\item[\textbf{(b)}] The spaces $L_p \big( \mathcal{N}_1;
L_q(\mathcal{N}_2) \big)$:

\vskip3pt

\begin{itemize} \item[$\bullet$] If $p \le q$ and $1/r =
1/p - 1/q$, we have $$L_p \big( \mathcal{N}_1; L_q(\mathcal{N}_2)
\big) = L_{2r}(\mathcal{N}_1) L_q(\mathcal{N}_1 \bar\otimes
\mathcal{N}_2) L_{2r}(\mathcal{N}_1).$$

\item[$\bullet$] If $p \ge q$ and $1/r = 1/q - 1/p$, we have
$$L_p \big( \mathcal{N}_1;
L_q(\mathcal{N}_2) \big) = L_{(2r,2r)}^p(\mathcal{N}_1 \bar\otimes
\mathcal{N}_2, \mathsf{E}),$$ where $\mathsf{E}: \mathcal{N}_1
\bar\otimes \mathcal{N}_2 \to \mathcal{N}_1$ is given by
$\mathsf{E} = id_{\mathcal{N}_1} \otimes \varphi_{\mathcal{N}_2}$.
\end{itemize}

\vskip3pt

\item[\textbf{(c)}] The spaces $L_p^r(\mathcal{M}, \mathsf{E})$
and $L_p^c(\mathcal{M}, \mathsf{E})$:

\vskip3pt

\begin{itemize} \item[$\bullet$] If $1 \le p \le 2$ and
$1/p = 1/2 + 1/s$, we have
\begin{eqnarray*}
L_p^r(\mathcal{M}, \mathsf{E}) & = & L_s(\mathcal{N})
L_2(\mathcal{M}) L_{\infty}(\mathcal{N}), \\ L_p^c(\mathcal{M},
\mathsf{E}) & = & L_{\infty}(\mathcal{N}) L_2(\mathcal{M})
L_s(\mathcal{N}).
\end{eqnarray*}

\item[$\bullet$] If $2 \le p \le \infty$ and $1/p + 1/s = 1/2$, we
have
\begin{eqnarray*} L_p^r(\mathcal{M},
\mathsf{E}) & = & L_{(s,\infty)}^p(\mathcal{M}, \mathsf{E}), \\
L_p^c(\mathcal{M}, \mathsf{E}) & = & L_{(\infty,s)}^p(\mathcal{M},
\mathsf{E}).
\end{eqnarray*}
\end{itemize}
By \eqref{AEq-RpCp-Cond}, we have also representations for
$L_p(\mathcal{M}, R_p^n)$ and $L_p(\mathcal{M}, C_p^n)$.

\vskip5pt

\item[\textbf{(d)}] Along the paper, we shall also find
representations of asymmetric spaces $L_{(u,v)} (\mathcal{M})$ (a
non-standard operator space structure on $L_p$ defined below which
will be crucial in this paper) in terms of either amalgamated or
conditional $L_p$ spaces. This will be a key point since we need
to handle spaces of the form $S_p (L_{(u,v)} (\mathcal{M}))$. The
use of amalgamated $L_p$ spaces or conditional $L_p$ spaces
depends on the sign of the term $1/u + 1/v - 1/p$.
\end{itemize}

\vskip5pt

Now we collect the complex interpolation and duality properties of
amalgamated and conditional $L_p$ spaces from \cite{JP2}. Our
interpolation identities generalize some previous results by
Pisier \cite{P0} and very recently by Xu \cite{X}. We need to
consider the following subset of the solid $\mathsf{K}$
$$\mathsf{K}_0 = \Big\{ (1/u,1/v,1/q) \in \mathsf{K} \, \big| \ 2
< u,v \le \infty, \ 1 < q < \infty, \ 1/u + 1/q + 1/v < 1
\Big\}.$$

\begin{theorem} \label{TheoApp1}
The following properties hold$\, :$
\begin{itemize}
\item[\textbf{a)}] If $(1/u,1/v,1/q) \in \mathsf{K}$,
$L_u(\mathcal{N}) L_q(\mathcal{M}) L_v(\mathcal{N})$ is a Banach
space.

\item[\textbf{b)}] If $(1/u_j,1/v_j,1/q_j) \in \mathsf{K}$ for
$j=0,1$ and $$(1/u_\theta,1/v_\theta,1/q_\theta) = \sum_{j=0,1} |1
- j - \theta| (1/u_j,1/v_j,1/q_j),$$ the space
$L_{u_{\theta}}(\mathcal{N}) L_{q_{\theta}}(\mathcal{M})
L_{v_{\theta}}(\mathcal{N})$ is isometrically isomorphic to
$$\Big[L_{u_0}(\mathcal{N}) L_{q_0}(\mathcal{M})
L_{v_0}(\mathcal{N}), L_{u_1}(\mathcal{N}) L_{q_1}(\mathcal{M})
L_{v_1}(\mathcal{N}) \Big]_{\theta}^{\null}.$$

\item[\textbf{c)}] If $(1/u,1/v,1/q) \in \mathsf{K}_0$ and $1 -
1/p = 1/u + 1/q + 1/v$, we have $$\big( L_u(\mathcal{N})
L_q(\mathcal{M}) L_v(\mathcal{N}) \big)^* =
L_{(u,v)}^p(\mathcal{M}, \mathsf{E}),$$
$$\big( L_{(u,v)}^p(\mathcal{M}, \mathsf{E}) \big)^* =
L_u(\mathcal{N}) L_q(\mathcal{M}) L_v(\mathcal{N}).$$

\item[\textbf{d)}] In particular, we obtain the following
isometric isomorphisms
$$\Big[ L_{(u_0,v_0)}^{p_0}(\mathcal{M}, \mathsf{E}),
L_{(u_1,v_1)}^{p_1}(\mathcal{M}, \mathsf{E}) \Big]_{\theta} =
L_{(u_\theta,v_\theta)}^{p_\theta}(\mathcal{M}, \mathsf{E}).$$
\end{itemize}
\end{theorem}

\vskip5pt

In the following result we list some particular cases of Theorem
\ref{TheoApp1} under the restriction $p_0 = p_1$, since these are
the main interpolation identities used in this paper. The case
where both $p_0$ and $p_1$ are $\infty$ is excluded in Theorem
\ref{TheoApp1}. That case was stated in Lemma \ref{Lemintinfty}
above and the proof was also given in \cite{JP2}.

\begin{corollary} \label{CorApp2}
If $2 \le p < \infty$, we set $$(1/u,1/v) = (\theta/q,
(1-\theta)/q) \quad \mbox{for \ $0 < \theta < 1$ and \ $q$ given
by} \quad 1/2 = 1/p + 1/q.$$ Then, the following isometric
isomorphisms hold
\begin{eqnarray*}
\big[ L_p(\mathcal{M}), L_p^r(\mathcal{M}, \mathsf{E})
\big]_{\theta} & = & L_{(u, \infty)}^p(\mathcal{M}, \mathsf{E}),
\\ \big[ L_p^c(\mathcal{M}, \mathsf{E}), L_p(\mathcal{M})
\big]_{\theta} & = & L_{(\infty,v)}^p(\mathcal{M}, \mathsf{E}), \\
\big[ L_p^c(\mathcal{M}, \mathsf{E}), L_p^r(\mathcal{M},
\mathsf{E}) \big]_{\theta} & = & L_{(u,v)}^p \ (\mathcal{M},
\mathsf{E}).
\end{eqnarray*}
\end{corollary}

\subsection{A variant of free Rosenthal's inequality}

In this paragraph we study a variation of the free Rosenthal
inequality \cite{JPX}, which will be applied in the sequel.
Intersection of $L_p$ spaces appear naturally in the theory of
noncommutative Hardy spaces. These spaces are also natural
byproducts of Rosenthal's inequality for sums of independent
random variables. Let us first illustrate this point in the
commutative setting and then provide the link to the spaces
defined above. Let $g_1, g_2, \ldots, g_n$ be a finite collection
of independent random variables on a probability space $(\Omega,
\mu)$. If $\varepsilon_1, \varepsilon_2, \ldots, \varepsilon_n$ is
an independent family of Bernoullis equidistributed on $\pm 1$,
the Khintchine inequality implies for $0 < s < \infty$ that \[
\Big( \int_\Omega \Big[ \sum_{k=1}^n |g_k|^2 \Big]^{\frac{s}{2}}
d\mu \Big)^{\frac{1}{s}} \sim_{c_s} \mathbb{E} \, \Big\|
\sum_{k=1}^n \varepsilon_k g_k \Big\|_s \, .\] Therefore,
Rosenthal's inequality \cite{Ro0} gives for $2 \le s < \infty$
\begin{equation} \label{Sp2}
\Big( \int_\Omega \Big[ \sum_{k=1}^n |g_k|^2 \Big]^{\frac{s}{2}}
d\mu \Big)^{\frac{1}{s}} \sim_{c_s} \Big( \sum_{k=1}^n \|g_k\|_s^s
\Big)^{\frac{1}{s}} + \Big( \sum_{k=1}^n \|g_k\|_2^2
\Big)^{\frac{1}{2}}.
\end{equation}
Now, given $1 \le q \le p < \infty$ and an independent family
$f_1, f_2, \ldots f_n$ of $p$-integrable random variables, we
define $g_k = |f_k|^{q/2}$ for $1 \le k \le n$. Then we have the
following identity for the index $s = 2p/q$
\begin{equation} \label{ChangeSp2}
\Big( \int_\Omega \Big[ \sum_{k=1}^n |f_k|^q \Big]^{\frac{p}{q}}
d\mu \Big)^{\frac1p} = \Big( \int_\Omega \Big[ \sum_{k=1}^n
|g_k|^2 \Big]^{\frac{s}{2}} d\mu \Big)^{\frac{2}{qs}}.
\end{equation}
Since the $g_k$'s are independent and $2 \le s < \infty$,
\eqref{Sp2} and \eqref{ChangeSp2} give
\renewcommand{\theequation}{$\Sigma_{pq}$}
\addtocounter{equation}{-1}
\begin{equation}
\Big( \int_\Omega \Big( \sum_{k=1}^n |f_k|^q \Big)^{\frac{p}{q}}
d\mu \Big)^{\frac1p} \sim_{c_p} \Big( \sum_{k=1}^n \|f_k\|_p^p
\Big)^{\frac{1}{p}} + \Big( \sum_{k=1}^n \|f_k\|_q^q
\Big)^{\frac{1}{q}}.
\end{equation}
\renewcommand{\theequation}{\arabic{equation}}
\numberwithin{equation}{section}

\vskip-5pt

\noindent In particular, ($\Sigma_{pq}$) provides a natural way to
realize the space
$$\mathcal{J}_{p,q}^n(\Omega) = n^{\frac1p} L_p(\Omega)
\cap n^{\frac1q} L_q(\Omega)$$ as an isomorph of a subspace of
$L_p(\Omega;\ell_q^n)$. More precisely, if $f_1, f_2, \ldots, f_n$
are taken to be independent copies of a given random variable $f$,
the right hand side of $(\Sigma_{pq})$ is the norm of $f$ in the
intersection space $\mathcal{J}_{p,q}^n(\Omega)$ and inequality
($\Sigma_{pq}$) provides an isomorphic embedding
$$f \in \mathcal{J}_{p,q}^n(\Omega) \mapsto (f_1,f_2, \ldots, f_n)
\in L_p(\Omega; \ell_q^n).$$

\vskip5pt

Quite surprisingly, replacing independent variables by matrices of
independent variables in $(\Sigma_{pq})$ requires to intersect
\emph{four} spaces using the so-called \emph{asymmetric} $L_p$
spaces. In other words, the natural operator space structure of
$\mathcal{J}_{p,q}^n$ comes from a $4$-term intersection space.
This phenomenon was discovered for the first time in \cite{JP} and
we have already met it in Section \ref{NewSect2}. To justify this
point, instead of giving precise definitions we note that
H\"{o}lder inequality gives $L_p = L_{2p} L_{2p}$, meaning that
the $p$-norm of $f$ is the infimum of $\|g\|_{2p} \|h\|_{2p}$ over
all possible factorizations $f = g h$. If $L_{2p}^r$ and
$L_{2p}^c$ denote the row and column quantizations \eqref{Eq-RCLp}
of $L_{2p}$, the operator space analogue of this isometry is given
by the complete isometry
$$L_p = L_{2p}^r L_{2p}^c.$$ This will be further explained below.
In particular, according to the algebraic definition of
$L_p(\ell_q)$, the intersection space $\mathcal{J}_{p,q}^n$ has to
be redefined as the product
$$\mathcal{J}_{p,q}^n = \Big( n^{\frac{1}{2p}} L_{2p}^r \cap
n^{\frac{1}{2q}} L_{2q}^r \Big) \Big( n^{\frac{1}{2p}} L_{2p}^c
\cap n^{\frac{1}{2q}} L_{2q}^c \Big).$$ According to \cite{JP2},
we find
\begin{equation} \label{Eq-Aspect}
\mathcal{J}_{p,q}^n = n^{\frac{1}{p}} L_{2p}^r L_{2p}^c \cap
n^{\frac{1}{2p}+\frac{1}{2q}} L_{2p}^r L_{2q}^c \cap
n^{\frac{1}{2q} + \frac{1}{2p}} L_{2q}^r L_{2p}^c \cap
n^{\frac{1}{q}} L_{2q}^rL_{2q}^c.
\end{equation}
Of course, these notions are not rigorously defined and will be
analyzed in more detail below. Our only aim here is to motivate
the forthcoming results. Let us now see how the space in
\eqref{Eq-Aspect} generalizes our first definition of
$\mathcal{J}_{p,q}^n(\Omega)$. On the Banach space level we have
the isometries
$$L_{2p}^r L_{2q}^c = L_s = L_{2q}^r L_{2p}^c \quad \mbox{with}
\quad 1/s = 1/2p + 1/2q.$$ Moreover, again by H\"{o}lder
inequality it is clear that
$$n^{\frac1s} \|f\|_s \le \max \Big\{ n^{\frac1p}
\|f\|_p, n^{\frac1q} \|f\|_q \Big\}.$$ Therefore, the two cross
terms in the middle of \eqref{Eq-Aspect} disappear in the Banach
space level. However, as we shall see in this section, in the
category of operator spaces all the four terms have a significant
contribution. The operator space/free version of $(\Sigma_{pq})$
is the main result in \cite{JP2}. It is worthy of mention that
this result goes a bit further than its commutative counterpart.
More precisely, in contrast with the classical case, we find a
\emph{right} formulation for ($\Sigma_{\infty q}$). Indeed, as it
happens with the Khintchine and Rosenthal inequalities, the limit
case as $p \to \infty$ holds when replacing independence by
Voiculescu's concept of freeness \cite{VDN}. Unfortunately, the
techniques for $(p,q) = (\infty,2)$ used in Section \ref{NewSect2}
do not apply in the general case and the arguments in \cite{JP2}
become more involved. This is mainly because a concrete Fock space
representation does not seem available for $L_p({\mathcal A})$
with $\mathcal{A}$ a free product algebra and $p<\infty$.
Therefore, a \emph{purely free} proof of $(\Sigma_{pq})$ seems out
of the scope by now. Nevertheless, since we shall need to be
familiar with the results in \cite{JP2}, we summarize them here.
We observe in advance that all the spaces and results presented in
this paragraph for $1 \le q \le p \le \infty$ are consistent with
their corresponding versions for $(p,q) = (\infty,2)$ used in
Section \ref{NewSect2}.

\vskip5pt

Now, if $2 \le u,v \le \infty$ and $1/p = 1/u + 1/v$ for some $1
\le p \le \infty$, we define the \emph{asymmetric $L_p$ space}
\label{LpAsimetrico} associated to the pair $(u,v)$ as the
$\mathcal{M}$-amalgamated Haagerup tensor product
\begin{equation} \label{Eq-Asymmetric}
L_{(u,v)}(\mathcal{M}) = L_u^r(\mathcal{M})
\otimes_{\mathcal{M},h} L_v^c(\mathcal{M}).
\end{equation}
That is, we consider the quotient of $L_u^r(\mathcal{M}) \otimes_h
L_v^c(\mathcal{M})$ by the closed subspace $\mathcal{I}$ generated
by the differences $x_1 \gamma \otimes x_2 - x_1 \otimes \gamma
x_2$ with $\gamma \in \mathcal{M}$. Recall that the row and column
operator space structures on $L_u(\mathcal{M})$ and
$L_v(\mathcal{M})$ have been already defined in \eqref{Eq-RCLp}.
By a well-known factorization argument, see e.g. Lemma 3.5 in
\cite{P2}, the norm of an element $x$ in $L_{(u,v)}(\mathcal{M})$
is given by
$$\|x\|_{(u,v)}^{\null} = \inf_{x = \alpha \beta}
\|\alpha\|_{L_u(\mathcal{M})} \|\beta\|_{L_v(\mathcal{M})}.$$

\begin{remark} \label{Rem-Lp-2p2p}
\emph{We have a complete isometry $L_p(\mathcal{M}) =
L_{(2p,2p)}(\mathcal{M})$.}
\end{remark}

\begin{remark} \label{Remark-Asymm-Schatten}
\emph{Asymmetric $L_p$ spaces were introduced in \cite{JP} for
matrix algebras. In fact, if $\mathcal{M}$ is the algebra
$\mathrm{M}_m$ of $m \times m$ matrices, we define the
\emph{asymmetric Schatten $p$-class} as follows
$$S_{(u,v)}^m = C_{u/2}^m \otimes_h R_{v/2}^m.$$ As observed in
\cite{JP2}, this definition is consistent with our definition
\eqref{Eq-Asymmetric}.}
\end{remark}

According to the discussion which led to \eqref{Eq-Aspect}, we
know how the general aspect of $\mathcal{J}_{p,q}^n(\mathcal{M})$
should be. Now, equipped with asymmetric $L_p$ spaces we know how
to factorize noncommutative $L_p$ spaces in the right way and
define $$\mathcal{J}_{p,q}^n(\mathcal{M}) = \bigcap_{u,v \in
\{2p,2q\}} n^{\frac{1}{2p} + \frac{1}{2q}} \,
L_{(u,v)}(\mathcal{M}).$$ The following result generalizes
\eqref{ossJ}, see \cite{JP2} for the proof.

\begin{lemma} \label{Lem-oss-Jpq}
If we take $$\mathcal{M}_m = \mathrm{M}_m(\mathcal{M}) \quad
\mbox{and} \quad \mathsf{E}_m = id_{\mathrm{M}_m} \otimes \varphi:
\mathcal{M}_m \to \mathrm{M}_m$$ for $m \ge 1$ and consider the
index $1/r = 1/q - 1/p$, we have an isometry $$S_p^m \big(
\mathcal{J}_{p,q}^n(\mathcal{M}) \big) = \bigcap_{u,v \in
\{2r,\infty\}}^{\null} n^{\frac{1}{u} + \frac{1}{p} + \frac{1}{v}}
\, L_{(u,v)}^p(\mathcal{M}_m, \mathsf{E}_m).$$
\end{lemma}

\begin{remark}
\emph{According to Lemma \ref{Lem-oss-Jpq}, we set
$$\mathcal{J}_{p,q}^n(\mathcal{M}, \mathsf{E}) = \bigcap_{u,v \in
\{2r,\infty\}}^{\null} n^{\frac{1}{u} + \frac{1}{p} + \frac{1}{v}}
\, L_{(u,v)}^p(\mathcal{M}, \mathsf{E}).$$}
\end{remark}

Lemma \ref{Lem-oss-Jpq} shows us the way to work in what follows.
Indeed, instead of working with the o.s.s. of the spaces
$\mathcal{J}_{p,q}^n(\mathcal{M})$, it suffices to argue with the
Banach space structure of the more general spaces
$\mathcal{J}_{p,q}^n(\mathcal{M}, \mathsf{E})$. In this spirit,
for $1 \le q \le p \le \infty$ we set $1/r = 1/q - 1/p$ and
introduce the spaces
\begin{eqnarray*}
\mathcal{R}_{2p,q}^n(\mathcal{M}, \mathsf{E}) & = &
n^{\frac{1}{2p}} \, L_{2p}(\mathcal{M}) \, \cap \,
n^{\frac{1}{2q}} \, L_{(2r,\infty)}^{2p} (\mathcal{M},
\mathsf{E}), \\ \mathcal{C}_{2p,q}^n \, (\mathcal{M}, \mathsf{E})
& = & n^{\frac{1}{2p}} \, L_{2p}(\mathcal{M}) \, \cap \,
n^{\frac{1}{2q}} \, L_{(\infty,2r)}^{2p} (\mathcal{M},
\mathsf{E}).
\end{eqnarray*}

\begin{remark} \label{Re-oM}
\emph{Let $(\mathrm{X}_1, \mathrm{X}_2)$ be a pair of operator
spaces containing a von Neumann algebra $\mathcal{M}$ as a common
two-sided ideal. We define the \emph{amalgamated} Haagerup tensor
product $\mathrm{X}_1 \otimes_{\mathcal{M}, h} \mathrm{X}_2$ as
the quotient of $\mathrm{X}_1 \otimes_h \mathrm{X}_2$
\label{AmalHaa} by the closed subspace $\mathcal{I}$ generated by
the differences $x_1 \gamma \otimes x_2 - x_1 \otimes \gamma x_2$
with $\gamma \in \mathcal{M}$. This notion has already been used
above in the definition of asymmetric $L_p$ spaces. Let us write
$\mathrm{X}_1 \oM \mathrm{X}_2$ to denote the underlying Banach
space of $\mathrm{X}_1 \otimes_{\mathcal{M},h} \mathrm{X}_2$. Our
definition uses the operator space structure of the
$\mathrm{X}_j$'s since the row (resp. column) square functions are
not necessarily closed operations in $\mathrm{X}_1$ (resp.
$\mathrm{X}_2$). However, in the sequel it will be important to
note that much less structure on $(\mathrm{X}_1, \mathrm{X}_2)$ is
needed to define the norm in $\mathrm{X}_1 \otimes_{\mathcal{M}}
\mathrm{X}_2$. Indeed, we just need to impose conditions under
which the row and column square functions become closed operations
in $\mathrm{X}_1$ and $\mathrm{X}_2$ respectively. In particular,
this is guaranteed if $\mathrm{X}_1$ is a right
$\mathcal{M}$-module and $\mathrm{X}_2$ is a left
$\mathcal{M}$-module, see Chapter 6 of \cite{JP2} for further
details. Therefore, we may define the Banach space
$$\mathcal{R}_{2p,q}^n(\mathcal{M}, \mathsf{E}) \oM
\mathcal{C}_{2p,q}^n(\mathcal{M}, \mathsf{E}).$$}
\end{remark}

\vskip5pt

\noindent The theorem below collects the key results in
\cite{JP2}.

\begin{theorem} \label{Theorem-MAMS}
The following isomorphisms hold$\, :$

\vskip5pt

\begin{itemize}
\item[\textbf{a)}] If $1 \le p \le \infty$ and $1/q = 1-\theta +
\theta/p$, we have
\begin{eqnarray*}
\big[ \mathcal{R}_{2p,1}^n (\mathcal{M}, \mathsf{E}),
\mathcal{R}_{2p,p}^n (\mathcal{M}, \mathsf{E}) \big]_{\theta} &
\simeq & \mathcal{R}_{2p,q}^n (\mathcal{M}, \mathsf{E}), \\ \big[
\hskip1.5pt \mathcal{C}_{2p,1}^n \, (\mathcal{M}, \mathsf{E}), \,
\mathcal{C}_{2p,p}^n \, (\mathcal{M}, \mathsf{E}) \big]_{\theta} &
\simeq & \, \mathcal{C}_{2p,q}^n \, (\mathcal{M}, \mathsf{E}).
\end{eqnarray*}

\item[\textbf{b)}] If $1 \le p \le \infty$, we have $$\qquad \big[
\mathcal{R}_{2p,1}^n(\mathcal{M},\mathsf{E}),
\mathcal{C}_{2p,1}^n(\mathcal{M},\mathsf{E}) \big]_{\theta} \
\simeq \bigcap_{u,v \in \{ 2p', \infty \}}^{\null}
n^{\frac{1-\theta}{u} + \frac{1}{2p} + \frac{\theta}{v}} \,
L_{(\frac{u}{1-\theta},\frac{v}{\theta})}^{2p} (\mathcal{M},
\mathsf{E}).$$

\item[\textbf{c)}] If $1 \le q \le p \le \infty$, we have
$$\mathcal{J}_{p,q}^n(\mathcal{M}, \mathsf{E}) \simeq
\mathcal{R}_{2p,q}^n(\mathcal{M}, \mathsf{E}) \oM
\mathcal{C}_{2p,q}^n(\mathcal{M}, \mathsf{E}).$$

\item[\textbf{d)}] If $1 \le p \le \infty$ and $1/q = 1-\theta +
\theta/p$, we have
$$\mathcal{J}_{p,q}^n(\mathcal{M}, \mathsf{E}) \simeq \big[
\mathcal{J}_{p,1}^n(\mathcal{M}, \mathsf{E}),
\mathcal{J}_{p,p}^n(\mathcal{M}, \mathsf{E}) \big]_{\theta}.$$
\end{itemize}

\vskip3pt

\noindent Moreover, the involved relevant constants are in all
cases independent of $n$.
\end{theorem}

Of course, the main result in the proof of the free analogue of
$(\Sigma_{pq})$ is the last interpolation isomorphism in (d). In
contrast with the case $(p,q) = (\infty,2)$, its proof does not
follow from (b) but from the combination of (a) and (c). Indeed,
the right hand side of (b) only gives a $\mathcal{J}$-space when
$\theta = 1/2$. This is related to Remark \ref{Remark-HalfWay}
above. The following corollary gives the free analogue of
$(\Sigma_{pq})$ in the operator space case
($\mathcal{J}_{p,q}^n(\mathcal{M})$ spaces) and in the amalgamated
case ($\mathcal{J}_{p,q}^n(\mathcal{M}, \mathsf{E})$ spaces). To
that aim we take again $\mathsf{A}_k = \mathcal{M} \oplus
\mathcal{M}$ for $1 \le k \le n$.

\renewcommand{\theequation}{$\Sigma_{pq}$}
\addtocounter{equation}{-1}

\begin{corollary} \label{Corollary-Sigma-pq}
If $\mathcal{A}_\mathcal{N} = *_{\mathcal{N}} \mathsf{A}_k$, the
map $$u: x \in \mathcal{J}_{p,q}^n(\mathcal{M}, \mathsf{E})
\mapsto \sum_{k=1}^n x_k \otimes \delta_k \in
L_p(\mathcal{A}_\mathcal{N}; \ell_q^n)$$ is an isomorphism with
complemented image and constants independent of $n$. In
particular, replacing as usual $(\mathcal{M}, \mathcal{N},
\mathsf{E})$ by $(\mathcal{M}_m, \mathrm{M}_m, \mathsf{E}_m)$ and
replacing $\mathcal{A}_\mathcal{N}$ by the non-amalgamated algebra
$\mathcal{A}_\C = \mathsf{A}_1
* \mathsf{A}_2 * \cdots * \mathsf{A}_n$, we obtain a
cb-isomorphism with cb-complemented image and constants
independent of $n$
\begin{equation}
\sigma: x \in \mathcal{J}_{p,q}^n(\mathcal{M}) \mapsto
\sum_{k=1}^n x_k \otimes \delta_k \in L_p(\mathcal{A}_\C;
\ell_q^n).
\end{equation}
\end{corollary}

\renewcommand{\theequation}{\arabic{equation}}
\numberwithin{equation}{section}

\section{Construction of the main embedding}
\label{NewSect4}

Now we have all the tools to prove our main result. In the first
paragraph we embed the Schatten class $S_q$ into
$L_p(\mathcal{A})$ for some $\mathrm{QWEP}$ von Neumann algebra
$\mathcal{A}$. Roughly speaking, the proof is almost identical to
our original argument after replacing Corollary
\ref{Corollary-Sigma-Infty-2} by Corollary
\ref{Corollary-Sigma-pq}. In particular, we shall omit some
details in our construction and include some others, such as the
proofs of Lemmas \ref{Lemma-OH-graph} and
\ref{Lemma-Diagonalization} which we applied in our proof of
Theorem \ref{QSQSQS}. The second paragraph is devoted to the
stability of hyperfiniteness and there we will present the
transference argument mentioned in the Introduction. Finally, the
last paragraph contains our construction for general von Neumann
algebras.

\subsection{Embedding Schatten classes}
\label{SSS3.1}

We being by embedding the Schatten class $S_q$ into
$L_p(\mathcal{A})$ for some $\mathrm{QWEP}$ von Neumann algebra.
In fact, we shall prove a more general statement (a generalization
of Theorem \ref{QSQSQS}) for which we need some preliminaries.
Although the following results might be well-known, we state them
in detail since they will be key tools in our construction. The
next lemma has been known to Xu and the first-named author for
quite some time. We refer to Xu's paper \cite{X3} for an even more
general statement than the result presented below.

\begin{lemma} \label{Lemma-Xu-Partition}
Given $1 \le p \le \infty$ and a closed subspace $\mathrm{X}$ of
$R_p \oplus_2 \mathrm{OH}$, there exist closed subspaces
$\mathcal{H}_1, \mathcal{H}_2, \mathcal{K}_1, \mathcal{K}_2$ of
$\ell_2$ and an injective closed densely-defined operator
$\Lambda: \mathcal{K}_1 \to \mathcal{K}_2$ with dense range such
that $$\mathrm{X} \simeq_{cb} \mathcal{H}_{1,r_p} \oplus_2
\mathcal{H}_{2,oh} \oplus_2 graph(\Lambda),$$ where the graph of
$\Lambda$ is regarded as a subspace of $\mathcal{K}_{1,r_p}
\oplus_2 \mathcal{K}_{2,oh}$ and the relevant constants in the
complete isomorphism above do not depend on the subspace
$\mathrm{X}$. Moreover, since $R_p = C_{p'}$ the same result can
be written in terms of column spaces.
\end{lemma}

In the following, we shall also need to recognize Pisier's
operator Hilbert space $\mathrm{OH}$ as the graph of certain
diagonal operator on $\ell_2$. More precisely, the following
result will be used below.

\begin{lemma} \label{Lemma-OH-graph}
Given $1 \le p \le \infty$, there exists a sequence $\lambda_1,
\lambda_2, \ldots$ in $\mathbb{R}_+$ for which the associated
diagonal map $\mathsf{d}_\lambda = \sum_k \lambda_k e_{kk}: R_p
\to \mathrm{OH}$ satisfies the following complete isomorphism
$$\mathrm{OH} \simeq_{cb} graph (\mathsf{d}_\lambda).$$
\end{lemma}

\dem Let us define
$$u: \delta_k \in \mathrm{OH} \mapsto (\lambda_k^{-1} \delta_k,
\delta_k) \in graph(\mathsf{d}_\lambda).$$ The mapping $u$
establishes a linear isomorphism between $\mathrm{OH}$ and
$graph(\mathsf{d}_\lambda)$. The inverse map of $u$ is the
coordinate projection into the second component, which is clearly
a complete contraction. Regarding the cb-norm of $u$, since
$graph(\mathsf{d}_\lambda)$ is equipped with the o.s.s. of $R_p
\oplus_2 \mathrm{OH}$, we have
$$\|u\|_{cb} = \sqrt{1 + \xi^2}$$ with $\xi$ standing for the
cb-norm of $\mathsf{d}_{\lambda^{-1}}: \mathrm{OH} \to R_p$. We
claim that $$\xi \le \Big( \summ_k |\lambda_k^{-1}|^4
\Big)^{\frac14},$$ so that it suffices to take $\lambda_1,
\lambda_2, \ldots$ large enough to deduce the assertion. Indeed,
it is well-known that the inequality above holds for the map
$\mathsf{d}_{\lambda^{-1}}: \mathrm{OH} \to R$ and also for
$\mathsf{d}_{\lambda^{-1}}: \mathrm{OH} \to C$. Therefore, our
claim follows by complex interpolation. \fin

\begin{remark}
\emph{The constants in Lemma \ref{Lemma-OH-graph} are uniformly
bounded on $p$.}
\end{remark}

The main embedding result in Xu's paper \cite{X3} claims that any
quotient of a subspace of $C_p \oplus_p R_p$ cb-embeds in
$L_p(\mathcal{A})$ for some sufficiently large von Neumann algebra
$\mathcal{A}$ whenever $1 \le p < 2$. In particular, if $1 \le p <
q \le 2$, both $R_q$ and $C_q$ embed completely isomorphically in
$L_p(\mathcal{A})$ since both are in $\mathcal{QS}(C_p \oplus_p
R_p)$. The last assertion follows as in Lemma
\ref{Lemma-Motivation}. More precisely, Xu's construction can be
done either with $\mathcal{A}$ being the Araki-Woods quasi-free
CAR factor and also with Shlyakhtenko's generalization of it in
the free setting \cite{S}. In any case, $\mathcal{A}$ can be
chosen to be a $\mathrm{QWEP}$ type $\mathrm{III}_\lambda$ factor,
$0 < \lambda \le 1$. In our first embedding result of this
section, we generalize Xu's embedding.

\begin{theorem} \label{Theorem-QS}
If for some $1 \le p \le 2$ $$(\mathrm{X}_1, \mathrm{X}_2) \in
\mathcal{QS}(C_p \oplus_2 \mathrm{OH}) \times \mathcal{QS}(R_p
\oplus_2 \mathrm{OH}),$$ there exist a cb-embedding $\mathrm{X}_1
\otimes_h \mathrm{X}_2 \to L_p(\mathcal{A})$, for some
$\mathrm{QWEP}$ algebra $\mathcal{A}$.
\end{theorem}


The rest of this paragraph is devoted to the proof of Theorem
\ref{Theorem-QS}, which is formally identical to our proof of
Theorem \ref{QSQSQS}. In the first part of the proof, we reduce
the problem to the particular case where both $\mathrm{X}_1$ and
$\mathrm{X}_2$ are quotients over certain (annihilators of)
graphs.

\vskip5pt

\demI By injectivity of the Haagerup tensor product, we may assume
that $(\mathrm{X}_1, \mathrm{X}_2) \in \mathcal{Q}(C_p \oplus_2
\mathrm{OH}) \times \mathcal{Q}(R_p \oplus_2 \mathrm{OH})$. On the
other hand, recalling that $C_p = R_p^* = R_{p'}$ for $1 \le p \le
\infty$, we easily obtain from Lemma \ref{Lemma-Xu-Partition} and
duality the following cb-isomorphisms
\begin{eqnarray*}
\mathrm{X}_1 & \simeq_{cb} & \mathcal{H}_{11,c_p} \oplus_2
\mathcal{H}_{12,oh} \oplus_2 \Big( \big( \mathcal{K}_{11,c_p}
\oplus_2 \mathcal{K}_{12,oh} \big) \big/ graph(\Lambda_1)^\perp
\Big), \\ \mathrm{X}_2 & \simeq_{cb} & \mathcal{H}_{21,r_p}
\oplus_2 \mathcal{H}_{22,oh} \oplus_2 \Big( \big(
\mathcal{K}_{21,r_p} \oplus_2 \mathcal{K}_{22,oh} \big) \big/
graph(\Lambda_2)^\perp \Big),
\end{eqnarray*}
for certain subspaces $\mathcal{H}_{ij}, \mathcal{K}_{ij}$ $(1 \le
i,j \le 2)$ of $\ell_2$ and
$$\begin{array}{rrcl} \Lambda_1: & \mathcal{K}_{11,c_{p'}} & \to &
\mathcal{K}_{12,oh},
\\ \Lambda_2: & \mathcal{K}_{21,r_{p'}} & \to &
\mathcal{K}_{22,oh},
\end{array}$$
satisfying the properties stated in Lemma
\ref{Lemma-Xu-Partition}. Let us set
\begin{eqnarray*}
\mathcal{Z}_1 & = & \big( \mathcal{K}_{11,c_p} \oplus_2
\mathcal{K}_{12,oh} \big) \big/ graph(\Lambda_1)^\perp, \\
\mathcal{Z}_2 & = & \big( \mathcal{K}_{21,r_p} \oplus_2
\mathcal{K}_{22,oh} \big) \big/ graph(\Lambda_2)^\perp.
\end{eqnarray*}
Then, we have the following cb-isometric inclusion
\begin{eqnarray} \label{Eq-6terms}
\mathrm{X}_1 \otimes_h \mathrm{X}_2 & \subset \ \, & \mathcal{Z}_1
\otimes_h \mathcal{Z}_2 \\ \nonumber & \oplus_2 &
\mathcal{H}_{11,c_p} \otimes_h \mathrm{X}_2 \\ \nonumber &
\oplus_2 & \mathrm{X}_1 \otimes_h \mathcal{H}_{21,r_p} \\
\nonumber & \oplus_2 & \mathcal{H}_{12,oh} \otimes_h \mathcal{Z}_2
\\ \nonumber & \oplus_2 & \mathcal{Z}_1 \otimes_h
\mathcal{H}_{22,oh} \\ \nonumber & \oplus_2 & \mathcal{H}_{12,oh}
\otimes_h \mathcal{H}_{22,oh}.
\end{eqnarray}
Our reduction argument is quite similar to that in Theorem
\ref{QSQSQS}. Indeed, according to \cite{X3} we know that
$\mathrm{OH} \in \mathcal{QS}(C_p \oplus_p R_p)$ and that any
element in $\mathcal{QS}(C_p \oplus_p R_p)$ completely embeds in
$L_p(\mathcal{A})$ for some $\mathrm{QWEP}$ type $\mathrm{III}$
factor $\mathcal{A}$. This eliminates the last term in
\eqref{Eq-6terms}. The second and third terms embed into
$S_p(\mathrm{X}_1)$ and $S_p(\mathrm{X}_2)$ completely
isometrically. On the other hand, since $\mathrm{OH} \in
\mathcal{QS}(C_p \oplus_p R_p)$ and we have by hypothesis
$$\mathrm{X}_1 \in \mathcal{QS}(C_p \oplus_2 \mathrm{OH}) \quad
\mbox{and} \quad \mathrm{X}_2 \in \mathcal{QS}(R_p \oplus_2
\mathrm{OH}),$$ both $\mathrm{X}_1$ and $\mathrm{X}_2$ are
cb-isomorphic to an element in $\mathcal{QS}(C_p \oplus_p R_p)$.
Applying Xu's theorem \cite{X3} one more time, we may eliminate
these terms. Finally, for the fourth and fifth terms on the right
of \eqref{Eq-6terms}, we apply Lemma \ref{Lemma-OH-graph} and the
self-duality of $\mathrm{OH}$ to rewrite them as particular cases
of the first term $\mathcal{Z}_1 \otimes_h \mathcal{Z}_2$. \fin

\vskip5pt

Before continuing with the proof, we need more preparation. The
following discretization result might be also well-known.
Nevertheless, since we are not aware of any reference for it, we
include the proof for the sake of completeness.

\begin{lemma} \label{Lemma-Diagonalization}
Given $1 \le p \le \infty$ and a closed densely-defined operator
$\Lambda: R_p \to \mathrm{OH}$ with dense range in $\mathrm{OH}$,
there exists a diagonal operator $\mathsf{d}_{\lambda} = \sum_k
\lambda_k e_{kk}$ on $\ell_2$ such that, when regarded as a map
$\mathsf{d}_\lambda: R_p \to \mathrm{OH}$, we obtain
$$graph(\mathsf{d}_\lambda) \simeq_{cb} graph(\Lambda).$$
Moreover, the relevant constants in the cb-isomorphism above do
not depend on $\Lambda$.
\end{lemma}

\dem Let us first assume that $\Lambda$ is positive. Then, since
$R_p$ is separable we deduce from spectral calculus \cite{KR} that
there exists a $\sigma$-finite measure space $(\Omega,
\mathcal{F}, \mu)$ for which $\Lambda$ is similar to a
multiplication operator $$M_f: L_2(\Omega) \to L_2(\Omega).$$ Thus
we may assume $\Lambda=M_f$. Now, we employ a standard procedure
to create a diagonal operator. Given $\delta>0$, we may
approximate the function $f$ by an infinite simple function
$g=\summ_k (k\delta) 1_{k\delta<f\le (k+1)\delta}$. Replacing $f$
by $g$ yields a $1+\delta$ cb-isomorphism $graph (M_f) \simeq_{cb}
graph(M_g)$. Therefore, defining the measurable sets
$$\Omega_k = \Big\{ w \in \Omega \, \big| \ k \delta< f(w) \le
(k+1)\delta \Big\},$$ we have that $L_2(\Omega_k)$ is isomorphic
to $\ell_2(n_k)$ with $0 \le n_k = \dim L_2(\Omega_k) \le \infty$.
Choosing an orthonormal basis for $L_2(\Omega_k)$, we find that
$M_g$ is similar to $\mathsf{d}_{\lambda}$ where $\lambda_k = k
\delta$ with multiplicity $n_k$. This gives the assertion for
positive operators. If $\Lambda$ is not positive, we consider the
polar decomposition $\Lambda = u |\Lambda|$. By extension we may
assume that $u$ is a unitary. Thus, we get a cb-isometry $graph
(\Lambda) \simeq_{cb} graph (|\Lambda|)$. Thus, the general case
can be reduced to the case of positive operators. \fin

We now proceed as above. Given a sequence $\gamma_1, \gamma_2,
\ldots \in \mathbb{R}_+$, the diagonal map $\mathsf{d}_\gamma$ is
regarded as the density of a \emph{n.s.s.f.} weight $\psi$ on
$\mathcal{B}(\ell_2)$. We also keep the same terminology for $q_n,
\psi_n, \mathrm{k}_n, \varphi_n, \ldots$ If $1 \le p \le 2$, we
define the space $\mathcal{J}_{p',2}(\psi_n)$ as the subspace
$$\Big\{ \big( d_{\psi_n}^{\frac{1}{2p'}} z
d_{\psi_n}^{\frac{1}{2p'}}, d_{\psi_n}^{\frac{1}{2p'}} z
d_{\psi_n}^{\frac14}, d_{\psi_n}^{\frac14} z
d_{\psi_n}^{\frac{1}{2p'}}, d_{\psi_n}^{\frac14} z
d_{\psi_n}^{\frac14} \big) \, \big| \ z \in q_n
\mathcal{B}(\ell_2) q_n \Big\}$$ of the direct sum
$$\mathcal{L}_{p'}^n = \big( C_{p'}^n \otimes_h R_{p'}^n \big)
\oplus_2 \big( C_{p'}^n \otimes_h \mathrm{OH}_n \big) \oplus_2
\big( \mathrm{OH}_n \otimes_h R_{p'}^n \big) \oplus_2 \big(
\mathrm{OH}_n \otimes_h \mathrm{OH}_n \big).$$ In other words, we
may regard $\mathcal{J}_{p',2}(\psi_n)$ as an intersection of some
weighted forms of the asymmetric Schatten classes (see Remark
\ref{Remark-Asymm-Schatten}) considered above. Now we generalize
Lemmas \ref{Lemma-PredualMAMS11} and \ref{Lemma-Direct-Sum11} to
the present setting. Let us consider the dual space
$\mathcal{K}_{p,2}(\psi_n) = \mathcal{J}_{p',2}(\psi_n)^*$. We
assume as above (without lost of generality) that $\mathrm{k}_n =
\sum_{k=1}^n \gamma_k$ is an integer and define $\mathcal{A}_n$ as
in Lemma \ref{Lemma-PredualMAMS11}. If $\pi_j$ is the natural
embedding into the $j$-th component of $\mathcal{A}_n$ and we set
$x_j = \pi_j(x,-x)$, the following result is the $L_p$ version of
Lemma \ref{Lemma-PredualMAMS11}.

\begin{lemma} \label{Lemma-PredualMAMS22}
The mapping $$\omega: x \in \mathcal{K}_{p,2}(\psi_n) \mapsto
\frac{1}{\ \mathrm{k}_n} \sum_{j=1}^{\mathrm{k}_n} x_j \otimes
\delta_j \in L_p(\mathcal{A}_n; \mathrm{OH}_{\mathrm{k}_n})$$ is a
cb-embedding with cb-complemented image and constants independent
of $n$.
\end{lemma}

\dem The complete isometry $\mathcal{J}_{p',2}(\psi_n) =
\mathcal{J}_{p',2}^{\mathrm{k}_n}(q_n \mathcal{B}(\ell_2) q_n)$ is
all what we need since the argument is completed as in Lemma
\ref{Lemma-PredualMAMS11} replacing the use of Corollary
\ref{Corollary-Sigma-Infty-2} by its generalized form given in
Corollary \ref{Corollary-Sigma-pq}. On the other hand, we may
rewrite the space $\mathcal{J}_{p',2}(\psi_n)$ using the language
of conditional $L_p$ spaces. Indeed, let us consider the von
Neumann algebra $\mathcal{M}_n = q_n \mathcal{B}(\ell_2) q_n$
equipped with the state $\varphi_n$ which arises from the relation
$\psi_n = \mathrm{k}_n \varphi_n$. Then the space
$\mathcal{J}_{p',2}(\psi_n)$ has the form (see Remark
\ref{Remark-Asymm-Schatten})
$$\mathrm{k}_n^{\frac{1}{p'}} L_{p'}(\mathcal{M}_n) \cap
\mathrm{k}_n^{\frac{1}{2p'} + \frac{1}{4}}
L_{(2p',4)}(\mathcal{M}_n) \cap \mathrm{k}_n^{\frac{1}{4} +
\frac{1}{2p'}} L_{(4,2p')}(\mathcal{M}_n) \cap
\mathrm{k}_n^{\frac12} L_{(4,4)}(\mathcal{M}_n).$$ This is exactly
the definition of $\mathcal{J}_{p',2}^{\mathrm{k}_n}(q_n
\mathcal{B}(\ell_2) q_n)$ and the proof is complete. \fin

\begin{lemma} \label{Lemma-Direct-Sum}
If $\lambda_1, \lambda_2, \ldots \in \R_+$, we set
\begin{eqnarray*}
R_{p'} \cap \ell_2^{oh}(\lambda) & = & \mathrm{span} \, \Big\{
(\delta_k, \lambda_k \delta_k) \in R_{p'} \oplus_2 \mathrm{OH}
\Big\}, \\ C_{p'} \cap \ell_2^{oh}(\lambda) & = & \mathrm{span} \,
\Big\{ (\delta_k, \lambda_k \delta_k) \in C_{p'} \hskip0.5pt
\oplus_2 \mathrm{OH} \Big\}, \\ R_p + \ell_2^{oh}(\lambda) & = &
\big( R_p \oplus_2 \mathrm{OH} \big) \big/ \big( R_{p'} \cap
\ell_2^{oh}(\lambda) \big)^\perp, \\ [3pt] C_p +
\ell_2^{oh}(\lambda) & = & \big( C_p \oplus_2 \mathrm{OH} \big)
\big/ \big( C_{p'} \cap \ell_2^{oh}(\lambda) \big)^\perp.
\end{eqnarray*}
Then, there exists a n.s.s.f. weight $\psi$ on
$\mathcal{B}(\ell_2)$ such that
$$\big( C_p + \ell_2^{oh}(\lambda) \big) \otimes_h
\big( R_p + \ell_2^{oh}(\lambda) \big) = \mathcal{K}_{p,2}(\psi) =
\overline{\bigcup_{n \ge 1} \mathcal{K}_{p,2}(\psi_n)}.$$
\end{lemma}

\dem Let us define
\begin{eqnarray*}
q_n \big( R_p + \ell_2^{oh}(\lambda) \big) & = & \Big\{ (q_n(a),
q_n(b)) + \big( R_{p'} \cap \ell_2^{oh}(\lambda) \big)^\perp \,
\big| \ (a,b) \in R_p \oplus_2 \mathrm{OH} \Big\}, \\
q_n \big( C_p + \ell_2^{oh}(\lambda) \big) & = & \Big\{ (q_n(a),
q_n(b)) + \big( C_{p'} \cap \ell_2^{oh}(\lambda) \big) ^\perp \,
\big| \ (a,b) \in C_p \oplus_2 \mathrm{OH} \Big\}.
\end{eqnarray*}
Arguing as in the proof of Lemma \ref{Lemma-Direct-Sum11} ($q_n(x)
\to x$ as $n \to \infty$ in $R_p, \mathrm{OH}, C_p$), we may write
the Haagerup tensor product $\big( C_p + \ell_2^{oh}(\lambda)
\big) \otimes_h \big( R_p + \ell_2^{oh}(\lambda) \big)$ as the
direct limit below
$$\overline{\bigcup_{n \ge 1} q_n(C_p + \ell_2^{oh}(\lambda))
\otimes_h q_n(R_p + \ell_2^{oh}(\lambda))}.$$ This reduces the
problem to the finite-dimensional case. Arguing by duality, we
have to show that $q_n(C_{p'} \cap \ell_2^{oh}(\lambda)) \otimes_h
q_n(R_{p'} \cap \ell_2^{oh}(\lambda)) =
\mathcal{J}_{p',2}(\psi_n)$ for some \emph{n.s.s.f.} weight
$\psi$, where
\begin{eqnarray*}
q_n \big( C_{p'} \hskip1pt \cap \ell_2^{oh}(\lambda) \big) & = &
\mbox{span} \Big\{ (e_{i1}, \hskip1pt \lambda_i e_{i1} \hskip1pt )
\in C_{p'}^n \hskip1pt \oplus_2 \mathrm{OH}_n \Big\}, \\
q_n \big( R_{p'} \cap \ell_2^{oh}(\lambda) \big) & = & \mbox{span}
\Big\{ (e_{1j}, \lambda_j e_{1j}) \in R_{p'}^n \oplus_2
\mathrm{OH}_n \Big\}.
\end{eqnarray*}
Its Haagerup tensor product is the subspace
$$\mbox{span} \Big\{ (e_{ij}, \lambda_j e_{ij}, \lambda_i e_{ij},
\lambda_i \lambda_j e_{ij} ) \Big\} = \Big\{ (x, x
\mathsf{d}_{\lambda}, \mathsf{d}_\lambda x, \mathsf{d}_\lambda x
\mathsf{d}_\lambda) \, \big| \ x \in q_n \mathcal{B}(\ell_2) q_n
\Big\}$$ of the space $\mathcal{L}_{p'}^n$ defined above. Then, we
define $\gamma_k \in \R_+$ by the relation
$$\lambda_k = \gamma_k^{\frac14 - \frac{1}{2p'}} \quad \mbox{and}
\quad z = d_{\psi_n}^{-\frac{1}{2p'}} x
d_{\psi_n}^{-\frac{1}{2p'}},$$ where $\psi$ is the \emph{n.s.s.f.}
weight induced by $\mathsf{d}_\gamma$. This gives the space
$$q_n(\mathcal{G}_{p'}^c(\lambda)) \otimes_h
q_n(\mathcal{G}_{p'}^r(\lambda)) = \Big\{ \big(
d_{\psi_n}^{\frac{1}{2p'}} z d_{\psi_n}^{\frac{1}{2p'}},
d_{\psi_n}^{\frac{1}{2p'}} z d_{\psi_n}^{\frac{1}{4}},
d_{\psi_n}^{\frac{1}{4}} z d_{\psi_n}^{\frac{1}{2p'}},
d_{\psi_n}^{\frac{1}{4}} z d_{\psi_n}^{\frac{1}{4}} \big)
\Big\}.$$ The space on the right is by definition
$\mathcal{J}_{p',2}(\psi_n)$. This completes the proof. \fin

\demII By Lemma \ref{Lemma-Diagonalization} we may assume that the
graphs appearing in the terms $\mathcal{Z}_1$ and $\mathcal{Z}_2$
are graphs of diagonal operators $\mathsf{d}_{\lambda_1}$ and
$\mathsf{d}_{\lambda_2}$. By polar decomposition, perturbation and
complementation (as in the proof of Theorem \ref{QSQSQS} via
Lemmas \ref{Lemma-OH-graph} and \ref{Lemma-Diagonalization}), we
may assume that
\begin{eqnarray*}
\mathcal{Z}_1 & = & C_p + \ell_2^{oh}(\lambda), \\
\mathcal{Z}_2 & = & R_p + \ell_2^{oh}(\lambda),
\end{eqnarray*}
with $\lambda_1, \lambda_2, \ldots \in \R_+$ strictly positive.
According to Lemma \ref{Lemma-Direct-Sum}, we conclude that
$\mathcal{Z}_1 \otimes_h \mathcal{Z}_2$ can be identified with
$\mathcal{K}_{p,2}(\psi)$ for some \emph{n.s.s.f.} weight $\psi$
on $\mathcal{B}(\ell_2)$. It remains to construct a complete
embedding of $\mathcal{K}_{p,2}(\psi)$ into $L_p(\mathcal{A})$ for
some $\mathrm{QWEP}$ algebra $\mathcal{A}$. To that aim, we assume
without loss of generality that the $\mathrm{k}_n$'s are integers.
This allows us to proceed in the usual way. Namely, we first embed
$\mathcal{K}_{p,2}(\psi)$ into an ultraproduct
$$\mathcal{K}_{p,2}(\psi) = \overline{\bigcup_{n \ge 1}
\mathcal{K}_{p,2}(\psi_n)}^{\null} \to \prodd_{n, \mathcal{U}}
\mathcal{K}_{p,2}(\psi_n).$$ According to \cite{Ra}, this reduces
the problem to the finite-dimensional case, which follows from
Xu's cb-embedding \cite{X3} of $\mathrm{OH}$ into
$L_p(\mathcal{B})$ for some $\mathrm{QWEP}$ type $\mathrm{III}$
factor $\mathcal{B}$ and from Lemma \ref{Lemma-PredualMAMS22}
$$\mathcal{K}_{p,2}(\psi_n) \to L_p(\mathcal{A}_n;
\mathrm{OH}_{\mathrm{k}_n}) \to L_p(\mathcal{A}_n \bar\otimes
\mathcal{B}) = L_p(\mathcal{A}_n').$$ We have therefore
constructed a cb-embedding
$$\mathcal{Z}_1 \ten_h \mathcal{Z}_2 \to L_p(\mathcal{A}) \quad
\mbox{with} \quad \mathcal{A} = \Big( \prodd_{n, \mathcal{U}}
{\mathcal{A}_n'}_* \Big)^*.$$ The fact that $\mathcal{A}$ is
$\mathrm{QWEP}$ is justified as in the proof of Theorem
\ref{QSQSQS}. \fin

\begin{corollary} \label{SqQWEP}
$S_q$ cb-embeds into $L_p(\mathcal{A})$ for some $\mathrm{QWEP}$
algebra $\mathcal{A}$.
\end{corollary}

\dem Since $S_q = C_q \otimes_h R_q$, it follows from Lemma
\ref{Lemma-Motivation} and Theorem \ref{Theorem-QS}. \fin

\subsection{Embedding into the hyperfinite factor}
\label{SSS3.2}

Now we want to show that the cb-embedding $S_q \to
L_p(\mathcal{A})$ can be constructed with $\mathcal{A}$ being a
hyperfinite type III factor. Moreover, we shall prove some more
general results to be used in the next paragraph, where the
general cb-embedding will be constructed. We first establish a
transference argument, based on a noncommutative form of
Rosenthal's inequality for identically distributed random
variables in $L_1$ from \cite{J3}, which enables us to replace
freeness by some sort of independence.

\vskip5pt

Let $\mathcal{N}$ be a $\sigma$-finite von Neumann subalgebra of
some algebra $\mathcal{A}$ and let us consider a family
$\mathsf{A}_1, \mathsf{A}_2, \ldots$ of von Neumann algebras with
$\mathcal{N} \subset \mathsf{A}_k \subset \mathcal{A}$. As usual,
we require the existence of a \emph{n.f.} conditional expectation
$\mathsf{E}_\mathcal{N}: \mathcal{A} \to \mathcal{N}$. We recall
that $(\mathsf{A}_k)_{k \ge 1}$ is a system of \emph{indiscernible
independent copies over $\mathcal{N}$} (\emph{i.i.c.} in short)
when

\vskip3pt

\begin{itemize}
\item[\textbf{(i)}] If $a \in \big\langle \mathsf{A}_1,
\mathsf{A}_2, \ldots, \mathsf{A}_{k-1} \big\rangle$ and $b \in
\mathsf{A}_k$, we have $$\mathsf{E}_\mathcal{N}(ab) =
\mathsf{E}_\mathcal{N}(a) \mathsf{E}_\mathcal{N}(b).$$

\item[\textbf{(ii)}] There exist a von Neumann algebra
$\mathsf{A}$ containing $\mathcal{N}$, a normal faithful
conditional expectation $\mathsf{E}_0: \mathsf{A} \to \mathcal{N}$
and homomorphisms $\pi_k: \mathsf{A} \to \mathsf{A}_k$ such that
$$\mathsf{E}_\mathcal{N} \circ \pi_k = \mathsf{E}_0$$ and the
following holds for every strictly increasing function $\alpha:\N
\to \N$ $$\mathsf{E}_{\mathcal{N}} \big( \pi_{j_1}(a_1) \cdots
\pi_{j_m}(a_m) \big) = \mathsf{E}_\mathcal{N} \big(
\pi_{\alpha(j_1)}(a_1) \cdots \pi_{\alpha(j_m)}(a_m) \big).$$

\item[\textbf{(iii)}] There exist \emph{n.f.} conditional
expectations $\mathcal{E}_k: \mathcal{A} \to \mathsf{A}_k$ such
that
$$\mathsf{E}_\mathcal{N} = \mathsf{E}_0 \pi_k^{-1} \mathcal{E}_k
\quad \mbox{for all} \quad k \ge 1.$$
\end{itemize}
We shall say that $(\mathsf{A}_k)_{k \ge 1}$ are
\emph{symmetrically independent copies over $\mathcal{N}$}
(\emph{s.i.c.} in short) when the first condition above also holds
for $a$ in the algebra generated by $\mathsf{A}_1, \ldots,
\mathsf{A}_{k-1}, \mathsf{A}_{k+1}, \ldots$ and the second
condition holds for any permutation $\alpha$ of the integers. In
what follows, given a probability space $(\Omega, \mu)$, we shall
write $\varepsilon_1, \varepsilon_2, \ldots$ to denote an
independent family of Bernoulli random variables on $\Omega$
equidistributed on $\pm 1$. We now present the key inequality in
\cite{J3}.

\begin{lemma} \label{Lemma-Araki-Woods}
The following inequalities hold for $x \in L_1(\mathsf{A}):$
\begin{itemize}
\item[\textbf{a)}] If $(\mathsf{A}_k)_{k \ge 1}$ are i.i.c. over
$\mathcal{N}$, we have
\begin{eqnarray*}
\lefteqn{\qquad \int_\Omega \Big\| \sum_{k=1}^n \varepsilon_k
\pi_k(x) \Big\|_{L_1(\mathcal{A})} d\mu} \\ & \sim & \inf_{x =
a+b+c} n \|a\|_{L_1(\mathsf{A})} + \sqrt{n} \big\| \mathsf{E}_0
(bb^*)^{\frac12} \big\|_{L_1(\mathcal{N})} + \sqrt{n} \big\|
\mathsf{E}_0(c^*c)^{\frac12} \big\|_{L_1(\mathcal{N})}.
\end{eqnarray*}

\item[\textbf{b)}] If moreover, $\mathsf{E}_0(x)=0$ and
$(\mathsf{A}_k)_{k \ge 1}$ are s.i.c. over $\mathcal{N}$, then
\begin{eqnarray*}
\lefteqn{\qquad \Big\| \sum_{k=1}^n \pi_k(x)
\Big\|_{L_1(\mathcal{A})}} \\ & \sim & \inf_{x = a+b+c} n
\|a\|_{L_1(\mathsf{A})} + \sqrt{n} \big\| \mathsf{E}_0
(bb^*)^{\frac12} \big\|_{L_1(\mathcal{N})} + \sqrt{n} \big\|
\mathsf{E}_0(c^*c)^{\frac12} \big\|_{L_1(\mathcal{N})}.
\end{eqnarray*}
\end{itemize}
\end{lemma}

\dem We claim that $$\frac12 \Big\| \sum_{k=1}^n \pi_k(x)
\Big\|_{L_1(\mathcal{A})} \le \Big\| \sum_{k=1}^n \varepsilon_k
\pi_k(x) \Big\|_{L_1(\mathcal{A})} \le 2 \Big\| \sum_{k=1}^n
\pi_k(x) \Big\|_{L_1(\mathcal{A})}$$ for any choice of signs
$\varepsilon_1, \varepsilon_2, \ldots, \varepsilon_n$ whenever
$\mathsf{A}_1, \mathsf{A}_2, \ldots$ are symmetric independent
copies of $\mathsf{A}$ over $\mathcal{N}$ and $\mathsf{E}_0(x)=0$.
This establishes $\mathrm{(a)} \Rightarrow \mathrm{(b)}$ and so,
since the first assertion is proved in \cite{J3}, it suffices to
prove our claim. Such result will follow from the more general
statement $$\Big\| \sum_{k=1}^n \varepsilon_k \pi_k(x_k)
\Big\|_{L_1(\mathcal{A})} \le 2 \Big\| \sum_{k=1}^n \pi_k(x_k)
\Big\|_{L_1(\mathcal{A})},$$ for any family $x_1, x_2, \ldots,
x_n$ in $L_1(\mathsf{A})$ with $\mathsf{E}_0(x_k) = 0$ for $1 \le
k \le n$. Since we assume that $\mathcal{N}$ is $\sigma$-finite we
may fix a \emph{n.f.} state $\varphi$. We define $\phi = \varphi
\circ \mathsf{E}_\mathcal{N}$ and $\phi_0= \varphi \circ
\mathsf{E}_0$. According to \cite{Co} we have
$$\sigma_t^{\varphi} \circ \mathsf{E}_0 = \mathsf{E}_0 \circ
\sigma_t^{\phi_0} \quad \mbox{and} \quad \sigma_t^{\varphi} \circ
\mathsf{E}_\mathcal{N} = \mathsf{E}_\mathcal{N} \circ
\sigma_t^{\phi}.$$ Moreover, since $\mathsf{E}_\mathcal{N} =
\mathsf{E}_\mathcal{N} \circ \mathcal{E}_k$ we find $\phi = \phi
\circ \mathcal{E}_k$ which implies $$\sigma_t^{\phi} \circ
\mathcal{E}_k = \mathcal{E}_k \circ \sigma_t^{\phi}.$$ In
particular, $\sigma_t^{\phi}(\mathsf{A}_k) \subset \mathsf{A}_k$
for $k \ge 1$. Therefore, given any subset $\mathsf{S}$ of $\{1,2,
\ldots, n\}$ we find a $\phi$-invariant conditional expectation
$\mathsf{E}_\mathsf{S}: \mathcal{A} \to \mathcal{A}_\mathsf{S}$
where the von Neumann algebra $\mathcal{A}_\mathsf{S} = \langle
\mathsf{A}_k \, | \ k \in \mathsf{S} \rangle$. We claim that
$$\mathsf{E_S}(\pi_j(a))=0 \quad \mbox{whenever} \quad
\mathsf{E}_0(a)=0 \quad \mbox{and} \quad j \notin \mathsf{S}.$$
Indeed, let $b \in \mathcal{A}_\mathsf{S}$ and $a$ as above. Then
we deduce from symmetric independence
$$\phi \big( \mathsf{E_S}(\pi_j(a))b \big) = \phi \big( \pi_j(a)b
\big) = \phi \big( \mathsf{E}_\mathcal{N}(\pi_j(a)b) \big) = \phi
\big( \mathsf{E}_0(a) \mathsf{E}_\mathcal{N}(b) \big) = 0.$$ Thus
we may apply Doob's trick $$\Big\| \sum_{k \in \mathsf{S}}
\pi_k(x_k) \Big\|_{L_1(\mathcal{A})} = \Big\| \mathsf{E_S} \Big(
\sum_{k=1}^n \pi_k(x_k) \Big) \Big\|_{L_1(\mathcal{A})} \le \Big\|
\sum_{k=1}^n \pi_k(x_k) \Big\|_{L_1(\mathcal{A})}.$$ Then, the
claim follows taking $\{1,2,\ldots,n\} = \mathsf{S}_1 \cup
\mathsf{S}_{-1}$ with $\mathsf{S}_\alpha = \{ k: \, \varepsilon_k
= \alpha\}$. \fin

\begin{remark} \label{Remark-BR}
\emph{The inequalities in Lemma \ref{Lemma-Araki-Woods} generalize
the noncommutative Rosenthal inequality \cite{JX4} to the case
$p=1$ for identically distributed variables and under such notions
of noncommutative independence. Of course, in the case $1 < p < 2$
we have much stronger results from \cite{JX,JX4} and there is no
need of proving any preliminary result for our aims in this case.}
\end{remark}

Let us now generalize our previous definition of the space
$\mathcal{K}_{p,2}(\psi)$ to general von Neumann algebras. Let
$\mathcal{M}$ be a given von Neumann algebra, which we assume
$\sigma$-finite for the sake of clarity. Let $\mathcal{M}$ be
equipped with a \emph{n.s.s.f.} weight $\psi$. In other words,
$\psi$ is given by an increasing sequence (a net in the general
case) of pairs $(\psi_n,q_n)$ such that the $q_n$'s are increasing
finite projections in $\mathcal{M}$ with $\lim_n q_n = 1$ in the
strong operator topology and $\sigma_t^\psi(q_n) = q_n$. Moreover,
the $\psi_n$'s are normal positive functionals on $\mathcal{M}$
with support $q_n$ and satisfying the compatibility condition
$\psi_{n+1}(q_n x q_n) = \psi_n(x)$. As above, we shall write
$\mathrm{k}_n$ for the number $\psi_n(q_n) \in (0,\infty)$ and
(again as above) we may and will assume that the $\mathrm{k}_n$'s
are nondecreasing positive integers. In what follows we shall
write $d_{\psi_n}$ for the density on $q_n \mathcal{M} q_n$
associated to the \emph{n.f.} finite weight $\psi_n$. If $1 \le p
\le 2$, we define the space $\mathcal{J}_{p',2}(\psi_n)$ as the
closure of
$$\Big\{ \big( d_{\psi_n}^{\frac{1}{2p'}} z
d_{\psi_n}^{\frac{1}{2p'}}, d_{\psi_n}^{\frac{1}{2p'}} z
d_{\psi_n}^{\frac14}, d_{\psi_n}^{\frac14} z
d_{\psi_n}^{\frac{1}{2p'}}, d_{\psi_n}^{\frac14} z
d_{\psi_n}^{\frac14} \big) \, \big| \ z \in q_n \mathcal{M} q_n
\Big\}$$ in the direct sum
$$\mathcal{L}_{p'}^n = L_{p'}(q_n \mathcal{M} q_n)
\oplus_2 L_{(2p',4)}(q_n \mathcal{M} q_n) \oplus_2 L_{(4,2p')}(q_n
\mathcal{M} q_n) \oplus_2 L_2(q_n \mathcal{M} q_n).$$ In other
words, after considering the \emph{n.f.} state $\varphi_n$ on $q_n
\mathcal{M} q_n$ determined by the relation $\psi_n = \mathrm{k}_n
\varphi_n$ and recalling the definition of the spaces
$\mathcal{J}_{p,q}^n(\mathcal{M})$ from Section \ref{NewSect3}, we
may regard $\mathcal{J}_{p',2}(\psi_n)$ as the $4$-term
intersection space
$$\mathcal{J}_{p',2}(\psi_n) = \bigcap_{u,v \in \{2p',4\}}^{\null}
{\mathrm{k}_n}^{\frac1u + \frac1v} \, L_{(u,v)}(q_n \mathcal{M}
q_n) = \mathcal{J}_{p',2}^{\mathrm{k}_n}(q_n \mathcal{M} q_n).$$
Now we take direct limits and define
$$\mathcal{J}_{p',2} (\psi) = \overline{\bigcup_{n \ge 1}
\mathcal{J}_{p',2}(\psi_n)}^{\null},$$ where the closure is taken
with respect to the norm of the space $$\mathcal{L}_{p'} =
L_{p'}(\mathcal{M}) \oplus_2 L_{(2p',4)}(\mathcal{M}) \oplus_2
L_{(4,2p')}(\mathcal{M}) \oplus_2 L_2(\mathcal{M}).$$ To define
the space $\mathcal{K}_{p,2}(\psi)$ we also proceed as above and
consider $$\Psi_n: \mathcal{L}_p^n \to L_1(q_n \mathcal{M} q_n)$$
given by $$\Psi_n (x_1, x_2, x_3,x_4) = d_{\psi_n}^{\frac{1}{2p'}}
x_1 d_{\psi_n}^{\frac{1}{2p'}} + d_{\psi_n}^{\frac{1}{2p'}} x_2
d_{\psi_n}^{\frac{1}{4}} + d_{\psi_n}^{\frac{1}{4}} x_3
d_{\psi_n}^{\frac{1}{2p'}} + d_{\psi_n}^{\frac{1}{4}} x_4
d_{\psi_n}^{\frac{1}{4}}.$$ This gives $\ker \Psi_n =
\mathcal{J}_{p',2}(\psi_n)^{\perp}$ and we define
$$\mathcal{K}_{p,2}(\psi_n) = \mathcal{L}_p^n / \ker \Psi_n \quad
\mbox{and} \quad \mathcal{K}_{p,2}(\psi) = \overline{\bigcup_{n
\ge 1} \mathcal{K}_{p,2}(\psi_n)}^{\null},$$ where the latter is
understood as a quotient of $\mathcal{L}_p$. In other words, we
may regard the space $\mathcal{K}_{p,2}(\psi_n)$ as the sum of the
corresponding dual weighted asymmetric $L_p$ spaces considered in
the definition of $\mathcal{J}_{p',2}(\psi_n)$
\begin{equation} \label{KRescalado}
\mathcal{K}_{p,2}(\psi_n) = \sum_{u,v \in \{2p,4\}}^{\null}
{\mathrm{k}_n}^{\frac1u + \frac1v} \, L_{(u,v)}(q_n \mathcal{M}
q_n).
\end{equation}
Thus, using $\psi_n = \mathrm{k}_n \varphi_n$ backwards and taking
direct limits
$$\mathcal{K}_{p,2}(\psi) = L_p(\mathcal{M}) +
L_{(2p,4)}(\mathcal{M}) + L_{(4,2p)}(\mathcal{M}) +
L_2(\mathcal{M}),$$ where the sum is taken in $L_p(\mathcal{M})$
and the embeddings are given by
\begin{eqnarray*}
j_c(x): x \in L_{(2p,4)}(\mathcal{M}) & \mapsto & x
d_{\psi}^{\beta} \in L_p(\mathcal{M}), \\ j_r(x): x \in
L_{(4,2p)}(\mathcal{M}) & \mapsto & d_{\psi}^{\beta} x \in
L_p(\mathcal{M}),
\end{eqnarray*}
with $\beta = 1/2p - 1/4$, while the embedding of
$L_2(\mathcal{M})$ into $L_p(\mathcal{M})$ is given by $$j_2(x) =
d_{\psi}^{\beta} x d_{\psi}^{\beta}.$$

\begin{remark} \label{Remark-Wierd-Dual}
\emph{It will be important below to observe that our definition of
$\mathcal{K}_{p,2}(\psi_n)$ is slightly different to the one given
in the previous paragraph. Indeed, according to the usual duality
bracket $\langle x,y \rangle = \mbox{tr}(x^*y)$, we should have
defined
$$\mathcal{K}_{p,2}(\psi_n) = \sum_{u,v \in \{2p,4\}}^{\null}
{\mathrm{k}_n}^{-\gamma(u,v)} \, L_{(u,v)}(q_n \mathcal{M} q_n)
\quad \mbox{with} \quad \gamma(u,v) = \frac{1}{2(u/2)'} +
\frac{1}{2(v/2)'}.$$ This would give $\mathcal{K}_{p,2}(\psi_n) =
\mathcal{J}_{p',2}(\psi_n)^*$ and
$$\mathcal{K}_{p,2}(\psi_n) = \frac{1}{\mathrm{k}_n} \sum_{u,v
\in \{2p,4\}}^{\null} {\mathrm{k}_n}^{\frac1u + \frac1v} \,
L_{(u,v)}(q_n \mathcal{M} q_n).$$ However, we prefer to use
\eqref{KRescalado} in what follows for notational convenience.}
\end{remark}

Now we set some notation to distinguish between independent and
free random variables. If we fix a positive integer $n$, the von
Neumann algebra $\mathcal{A}_{ind}^n$ will denote the
$\mathrm{k}_n$-fold tensor product of $q_n \mathcal{M} q_n$ while
$\mathcal{A}_{free}^n$ will be (as usual) the $\mathrm{k}_n$-fold
free product of $q_n \mathcal{M} q_n \oplus q_n \mathcal{M} q_n$.
In other words, if we set $$\widetilde{\mathsf{A}}_{n,j} = q_n
\mathcal{M} q_n \quad \mbox{and} \quad \mathsf{A}_{n,j} = q_n
\mathcal{M} q_n \oplus q_n \mathcal{M} q_n$$ for $1 \le j \le
\mathrm{k}_n$, we define the following von Neumann algebras
\begin{eqnarray*}
\mathcal{A}_{ind}^n & = & \otimes_j \, \widetilde{\mathsf{A}}_{n,j}, \\
\mathcal{A}_{free}^n & = & \, *_j \ \mathsf{A}_{n,j}.
\end{eqnarray*}
We also consider the natural embeddings $$\pi_{ind}^j:
\widetilde{\mathsf{A}}_{n,j} \to \mathcal{A}_{ind}^n \quad
\mbox{and} \quad \pi_{free}^j: \mathsf{A}_{n,j} \to
\mathcal{A}_{free}^n.$$

\vskip5pt

We need some further information on $\mathrm{OH}$. Given $1 < p <
2$, Xu constructed in \cite{X3} a complete embedding of
$\mathrm{OH}$ into $L_p(\mathcal{A})$ with $\mathcal{A}$
hyperfinite, while for $p=1$ the corresponding cb-embedding was
recently constructed in \cite{J3}. The argument can be sketched
with the following chain
$$\mathrm{OH} \hookrightarrow \big( C_p \oplus_p R_p \big) /
graph(\mathsf{d}_\lambda)^{\perp} \simeq_{cb} C_p + R_p(\lambda)
\hookrightarrow L_p(\mathcal{A}).$$ Indeed, arguing as in Lemma
\ref{Lemma-Motivation}/Remark \ref{Remark-Graph_Sq} and applying
Lemma \ref{Lemma-Diagonalization}, we see how to regard
$\mathrm{OH}$ as a subspace of a quotient of $C_p \oplus_p R_p$ by
the annihilator of some diagonal map $\mathsf{d}_\lambda: C_{p'}
\to R_{p'}$. By the action of $\mathsf{d}_\lambda$, the
annihilator of its graph is the span of elements of the form
$(\delta_k, - \delta_k/\lambda_k)$. This suggest to regard the
quotient above as the sum of $C_p$ with a weighted form of $R_p$.
This establishes the cb-isomorphism in the middle. Then, it is
natural to guess that the complete embedding into
$L_p(\mathcal{A})$ should follow from a \emph{weighted} form of
the noncommutative Khintchine inequality. The first inequality of
this kind was given by Pisier and Shlyakhtenko in \cite{PS} for
generalized circular variables and further investigated in
\cite{JPX,X4}. However, if we want to end up with a hyperfinite
von Neumann algebra $\mathcal{A}$, we must replace generalized
circulars by their Fermionic analogues. More precisely, given a
complex Hilbert space $\mathcal{H}$, we consider its antisymmetric
Fock space $\mathcal{F}_{-1}(\mathcal{H})$. Let $c(e)$ and $a(e)$
denote the creation and annihilation operators associated with a
vector $e \in \mathcal{H}$. Given an orthonormal basis $(e_{\pm
k})_{k \ge 1}$ of $\mathcal{H}$ and a family $(\mu_k)_{k \ge 1}$
of positive numbers, we set $f_k = c(e_k) + \mu_k \, a(e_{-k})$.
The sequence $(f_k)_{k \ge 1}$ satisfies the canonical
anticommutation relations and we take $\mathcal{A}$ to be the von
Neumann algebra generated by the $f_k$'s. Taking suitable
$\mu_k$'s depending only on $p$ and the eigenvalues of
$\mathsf{d}_{\lambda}$, the Khintchine inequality associated to
the system of $f_k$'s provides the desired cb-embedding. Namely,
let $\phi$ be the quasi-free state on $\mathcal{A}$ determined by
the vacuum and let $d_\phi$ be the associated density. Then, if
$(\delta_k)_{k \ge 1}$ denotes the unit vector basis of
$\mathrm{OH}$, the cb-embedding has the form
$$w(\delta_k) = \xi_k \, d_\phi^{\frac{1}{2p}} f_k \,
d_\phi^{\frac{1}{2p}} = \xi_k \, f_{p,k}$$ for some scaling
factors $(\xi_k)_{k \ge 1}$. The necessary Khintchine type
inequalities for $1 < p < 2$ follow from the noncommutative
Burkholder inequality \cite{JX}. In the $L_1$ case, the key
inequalities follow from Lemma \ref{Lemma-Araki-Woods}, see
\cite{J3} for details. With this construction, the von Neumann
algebra $\mathcal{A}$ turns out to be the Araki-Woods factor
arising from the GNS construction applied to the CAR algebra with
respect to the quasi-free state $\phi$. In fact, using a
conditional expectation, we can replace the $\mu_k$'s by a
sequence $(\mu_k')_{k \ge 1}$ such that for every rational $0 <
\lambda < 1$ there are infinitely many $\mu_k'$'s with $\mu_k'=
\lambda/ (1 + \lambda)$. According to the results in \cite{AW}, we
then obtain the hyperfinite type $\mathrm{III}_1$ factor
$\mathcal{R}$.

\vskip5pt

On the other hand, there exists a slight modification of this
construction which will be used below. Indeed, using the
terminology introduced above and following \cite{J3} there exists
a mean-zero $\gamma_p \in L_p(\mathcal{R})$ given by a linear
combination of the $f_{p,k}$'s such that
$$w(\delta_j) = \pi_{ind}^j(\gamma_p)$$ defines a completely
isomorphic embedding $$w: \mathrm{OH}_{\mathrm{k}_n} \to
L_p(\mathcal{R}_{\otimes \mathrm{k}_n})$$ with constants
independent of $n$. Here $\mathcal{R}_{\otimes \mathrm{k}_n}$
denotes the $\mathrm{k}_n$-fold tensor product of $\mathcal{R}$.
Moreover, given any von Neumann algebra $\mathcal{A}$,
$id_{L_p(\mathcal{A})} \otimes w$ also defines an isomorphism
$$id_{L_p(\mathcal{A})} \otimes w: L_p
\big( \mathcal{A}; \mathrm{OH}_{\mathrm{k}_n} \big) \to
L_p(\mathcal{A} \bar\otimes \mathcal{R}_{\otimes \mathrm{k}_n}).$$

\begin{proposition} \label{trick}
If $1 \le p \le 2$, the map $$\xi_{ind}^n: x \in
\mathcal{K}_{p,2}(\psi_n) \mapsto \sum_{j=1}^{\mathrm{k}_n}
\pi_{ind}^j(x) \otimes \delta_j \in L_p \big( \mathcal{A}_{ind}^n;
\mathrm{OH}_{\mathrm{k}_n} \big)$$ is a completely isomorphic
embedding with relevant constants independent of $n$.
\end{proposition}

\dem According to Corollary \ref{Corollary-Sigma-pq}, this is true
for $$\xi_{free}^n: x \in \mathcal{K}_{p,2}(\psi_n) \mapsto
\sum_{j=1}^{\mathrm{k}_n} \pi_{free}^j(x,-x) \otimes \delta_j \in
L_p \big( \mathcal{A}_{free}^n; \mathrm{OH}_{\mathrm{k}_n}
\big).$$ Indeed, it follows from a simple duality argument (see
\cite{JP} or Remark 7.4 \cite{JP2}) taking Remark
\ref{Remark-Wierd-Dual} into account. According to the preceding
discussion on $\mathrm{OH}$, we deduce that $$w \circ
\xi_{free}^n: x \in \mathcal{K}_{p,2}(\psi_n) \mapsto \sum_{j=1}^n
\pi_{free}^j (x,-x) \otimes \pi_{ind}^j(\gamma_p) \in L_p \big(
\mathcal{A}_{free}^n \bar\otimes \mathcal{R}_{\otimes
\mathrm{k}_n} \big)$$ also provides a cb-isomorphism. Now, we
consider
\begin{eqnarray*}
\widetilde{\mathsf{B}}_{n,j} & = & \pi_{ind}^j
(\widetilde{\mathsf{A}}_{n,j})
\otimes \pi_{ind}^j(\mathcal{R}), \\
\mathsf{B}_{n,j} & = & \pi_{free}^j(\mathsf{A}_{n,j}) \otimes
\pi_{ind}^j(\mathcal{R}),
\end{eqnarray*}
for $1 \le j \le \mathrm{k}_n$. It is clear from the construction
that both families of von Neumann algebras are \emph{s.i.c.} over
the complex field. Therefore, Lemma \ref{Lemma-Araki-Woods}/Remark
\ref{Remark-BR} apply in both cases (note that the mean-zero
condition for the $\widetilde{\mathsf{B}}_{n,j}$'s holds due to
the fact that $\gamma_p$ is mean-zero) and hence
$$\big\| id_{S_p^m} \otimes w \circ \xi_{free}^n(x) \big\|_p
\sim_c \big\| id_{S_p^m} \otimes w \circ \xi_{ind}^n(x) \big\|_p$$
holds for every element $x\in L_p \big( \mathrm{M}_m (q_n
\mathcal{M} q_n) \big)$. This completes the proof. \fin

\begin{remark}
\emph{Proposition \ref{trick} can be regarded as a generalization
of Lemma \ref{Lemma-PredualMAMS22} for general von Neumann
algebras, where freeness is replaced by noncommutative
independence. Indeed, the only difference between both results is
the factor $1/ \mathrm{k}_n$, which is explained as a byproduct of
Remark \ref{Remark-Wierd-Dual}.}
\end{remark}

\begin{remark}
\emph{The transference argument applied in the proof of
Proposition \ref{trick} gives a result which might be of
independent interest. Given a von Neumann algebra $\mathsf{A}$,
let us construct the tensor product $\mathcal{A}_{ind}$ of
infinitely many copies of $\mathsf{A} \oplus \mathsf{A}$.
Similarly, the free product $\mathcal{A}_{free}$ of infinitely
many copies of $\mathsf{A} \oplus \mathsf{A}$ will be considered.
Following our terminology, we have maps $$\pi_{ind}^j: \mathsf{A}
\to \mathcal{A}_{ind} \quad \mbox{and} \quad \pi_{free}^j:
\mathsf{A} \to \mathcal{A}_{free}.$$ If $1 < p < q < 2$, we claim
that
\begin{equation} \label{Eq-Equivsimcb}
\Big\| \summ_j \pi_{ind}^j(x,-x) \ten \delta_j
\Big\|_{L_p(\mathcal{A}_{ind}; \ell_q)} \sim_{cb} \Big\| \summ_j
\pi_{free}^j(x,-x) \ten \delta_j \Big\|_{L_p(\mathcal{A}_{free};
\ell_q)},
\end{equation}
where the symbol $\sim_{cb}$ is used to mean that the equivalence
also holds (with absolute constants) when taking the matrix norms
arising from the natural operator space structures of the spaces
considered. The case $q=2$ follows by using exactly the same
argument as in Proposition \ref{trick}. In fact, the same idea
works for general indices. Indeed, we just need to embed $\ell_q$
into $L_p$ completely isometrically and then use the
noncommutative Rosenthal inequality \cite{JX4}. Recall that the
cb-embedding of $\ell_q$ into $L_p$ is already known at this stage
of the paper as a consequence of Corollary \ref{SqQWEP}. At the
time of this writing, it is still open whether or not
\eqref{Eq-Equivsimcb} is still valid for other values of $(p,q)$.}
\end{remark}

Our main goal in this paragraph is to generalize the complete
embedding in Proposition \ref{trick} to the direct limit
$\mathcal{K}_{p,2}(\psi)$. Of course, this is possible using an
ultraproduct procedure. However, this would not preserve
hyperfiniteness. We will now explain how the proof of Proposition
\ref{trick} allows to factorize the cb-embedding
$\mathcal{K}_{p,2}(\psi_n) \to L_p(\mathcal{A}_{ind}^n \bar\ten
\mathcal{R}_{\ten \mathrm{k}_n})$ via a three term
$\mathrm{K}$-functional. We will combine this with the concept of
noncommutative Poisson random measure from \cite{J6} to produce a
complete embedding which preserves the direct limit mentioned
above. Let us consider the operator space $$\mathcal{K}_{r \!
c_p}^p(\psi_n) = \mathrm{k}_n^{\frac1p} L_p(q_n \mathcal{M} q_n) +
\mathrm{k}_n^{\frac12} L_2^{r_p} (q_n \mathcal{M} q_n) +
\mathrm{k}_n^{\frac12} L_2^{c_p} (q_n \mathcal{M} q_n),$$ where
the norms in the $L_p$ spaces considered above are calculated with
respect to the state $\varphi_n$ arising from the relation $\psi_n
= \mathrm{k}_n \varphi_n$. More precisely, the operator space
structure is determined by
$$\|x\|_{S_p^m(\mathcal{K}_{r \! c_p}^p(\psi_n))} = \inf
\Big\{ \mathrm{k}_n^{\frac1p} \|x_1\|_{S_p^m(L_p)} +
\mathrm{k}_n^{\frac12} \|x_2\|_{S_p^m(L_2^{r_p})} +
\mathrm{k}_n^{\frac12} \|x_3\|_{S_p^m(L_2^{c_p})} \Big\},$$ where
the infimum runs over all possible decompositions $$x = x_1 +
d_{\varphi_n}^{\alpha} x_2 + x_3 d_{\varphi_n}^{\alpha},$$ with
$d_{\varphi_n}$ standing for the density associated to $\varphi_n$
and $\alpha = 1/p - 1/2$. Note that $\mathcal{K}_{r \!
c_p}^p(\psi_n)$ coincides algebraically with $L_p(q_n \mathcal{M}
q_n)$. There exists a close relation between $\mathcal{K}_{r \!
c_p}^p(\psi_n)$ and conditional $L_p$ spaces. Indeed, let us
consider the conditional expectation $\mathsf{E}_{\varphi_n}:
\mathrm{M}_m (q_n \mathcal{M} q_n) \to \mathrm{M}_m$ given by
$$\mathsf{E}_{\varphi_n} \big( (x_{ij}) \big) = \big(
\varphi_n(x_{ij}) \big) = \Big(
\frac{\psi_n(x_{ij})}{\mathrm{k}_n} \Big).$$

\begin{lemma} \label{Lemma-IsoK3Cond}
We have isometries
\begin{eqnarray*}
S_p^m \big( L_2^{r_p}(q_n \mathcal{M} q_n) \big) & = & m^{\frac1p}
L_p^r \big( \mathrm{M}_m (q_n \mathcal{M} q_n),
\mathsf{E}_{\varphi_n} \big), \\ S_p^m \big( L_2^{c_p}(q_n
\mathcal{M} q_n) \big) & = & m^{\frac1p} L_p^c \big( \mathrm{M}_m
(q_n \mathcal{M} q_n), \mathsf{E}_{\varphi_n} \big).
\end{eqnarray*}
Moreover, these isometries have the form
\begin{eqnarray*}
\big\| d_{\varphi_n}^{\frac12} a \big\|_{S_p^m (L_2^{r_p}(q_n
\mathcal{M} q_n))} & = & m^{\frac1p} \big\|
d_{\varphi_n}^{\frac1p} a \big\|_{L_p^r (\mathrm{M}_m (q_n
\mathcal{M} q_n), \mathsf{E}_{\varphi_n})}, \\ \big\| a
d_{\varphi_n}^{\frac12} \big\|_{S_p^m (L_2^{c_p}(q_n \mathcal{M}
q_n))} & = & m^{\frac1p} \big\| a d_{\varphi_n}^{\frac1p}
\big\|_{L_p^c (\mathrm{M}_m (q_n \mathcal{M} q_n),
\mathsf{E}_{\varphi_n})}.
\end{eqnarray*}
In particular, using the relation $d_{\varphi_n}^{\frac1p} =
d_{\varphi_n}^{\frac12} d_{\varphi_n}^{\alpha}$, we conclude
\begin{eqnarray*}
\lefteqn{\|x\|_{S_p^m(\mathcal{K}_{r \! c_p}^p(\psi_n))}} \\
& = & \inf_{x = x_p + x_r + x_c} \Big\{ \mathrm{k}_n^{\frac1p}
\|x_p\|_p  + \mathrm{k}_n^{\frac12} \big\|
\mathsf{E}_{\varphi_n}(x_r x_r^*)^{\frac12} \big\|_{S_p^m} +
\mathrm{k}_n^{\frac12} \big\| \mathsf{E}_{\varphi_n}(x_c^*
x_c)^{\frac12} \big\|_{S_p^m} \Big\}.
\end{eqnarray*}
\end{lemma}

\dem We have
\begin{eqnarray*}
\big\| d_{\varphi_n}^{\frac12} a \big\|_{S_p^m (L_2^{r_p}(q_n
\mathcal{M} q_n))} & = & \Big\| \mbox{tr}_{\mathcal{M}} \big(
d_{\varphi_n}^{\frac12} a a^* d_{\varphi_n}^{\frac12}
\big)^{\frac12} \hskip0.5pt \Big\|_{S_p^m} \\ & = & \Big\| \big(
id_{\mathrm{M}_m} \ten \varphi_n \big) (a a^*)^{\frac12}
\Big\|_{S_p^m}.
\end{eqnarray*}
Then, normalizing the trace on $\mathrm{M}_m$ and recalling that
$$ \Big\| \big( id_{\mathrm{M}_m} \ten \varphi_n \big) (a a^*)^{\frac12}
\Big\|_{S_p^m} = m^{\frac1p} \Big\| \mathsf{E}_{\varphi_n} \big(
d_{\varphi_n}^{\frac1p} a a^* d_{\varphi_n}^{\frac1p}
\big)^{\frac12} \Big\|_{L_p(\mathrm{M}_m(q_n \mathcal{M} q_n))}$$
when regarding the conditional expectation as a mapping
$$\mathsf{E}_{\varphi_n}: L_p \big( \mathrm{M}_m (q_n \mathcal{M}
q_n) \big) \to L_p(\mathrm{M}_m),$$ we deduce the assertion. The
column case is proved in the same way. \fin

\begin{proposition} \label{Prop-k4k3}
Let $\mathcal{R}$ be the hyperfinite $\mathrm{III}_1$ factor and
$\phi$ the quasi-free state on $\mathcal{R}$ considered above. Let
us consider the space $\mathcal{K}_{r \! c_p}^p(\phi \ten
\psi_n)$, defined as we did above. Then, there exists a completely
isomorphic embedding $$\rho_n: \mathcal{K}_{p,2}(\psi_n) \to
\mathcal{K}_{r \! c_p}^p(\phi \ten \psi_n).$$ Moreover, the
relevant constants in $\rho_n$ are independent of $n$.
\end{proposition}

\dem We have $$\mathcal{K}_{r \! c_p}^p(\phi \ten \psi_n) =
\mathrm{k}_n^{\frac1p} L_p(\mathcal{R} \bar\ten q_n \mathcal{M}
q_n) + \mathrm{k}_n^{\frac12} L_2^{r_p} (\mathcal{R} \bar\ten q_n
\mathcal{M} q_n) + \mathrm{k}_n^{\frac12} L_2^{c_p} (\mathcal{R}
\bar\ten q_n \mathcal{M} q_n).$$ The embedding is given by
$\rho_n(x) = \gamma_p \ten x$, with $\gamma_p$ the element of
$L_p(\mathcal{R})$ introduced before Proposition \ref{trick}.
Indeed, taking $\mathsf{E}_{\phi \ten \varphi_n}: \mathrm{M}_m
(\mathcal{R} \bar\ten q_n \mathcal{M} q_n) \to \mathrm{M}_m$ and
$x \in S_p^m(L_p(q_n \mathcal{M} q_n))$, we may argue as above and
obtain
\begin{eqnarray*}
\lefteqn{\|\gamma_p \ten x\|_{S_p^m(\mathcal{K}_{r \! c_p}^p(\phi
\ten \psi_n))}} \\ \!\!\!\! & = & \!\!\!\! \inf_{\gamma_p \ten x =
x_p + x_r + x_c} \Big\{ \mathrm{k}_n^{\frac1p} \|x_p\|_p +
\mathrm{k}_n^{\frac12} \big\| \mathsf{E}_{\phi \ten \varphi_n}(x_r
x_r^*)^{\frac12} \big\|_{S_p^m} + \mathrm{k}_n^{\frac12} \big\|
\mathsf{E}_{\phi \ten \varphi_n}(x_c^* x_c)^{\frac12}
\big\|_{S_p^m} \Big\}.
\end{eqnarray*}
Therefore, Lemma \ref{Lemma-Araki-Woods} and Remark
\ref{Remark-BR} give $$\|\gamma_p \ten x\|_{S_p^m(\mathcal{K}_{r
\! c_p}^p(\phi \ten \psi_n))} \sim \Big\|
\sum_{k=1}^{\mathrm{k}_n} \pi_{ind}^j(\gamma_p \ten x) \Big\|_p
\sim \Big\| \sum_{j=1}^{\mathrm{k}_n} \pi_{ind}^j(\gamma_p) \ten
\pi_{free}^j(x,-x) \Big\|_p.$$ Hence, the assertion follows as in
Proposition \ref{trick}. The proof is complete. \fin

Let $\mathcal{M}$ be a von Neumann algebra equipped with a
\emph{n.s.s.f.} weight $\psi$ and let us write $(\psi_n,q_n)_{n
\ge 1}$ for the associated sequence of $q_n$-supported weights.
Then we define the following direct limit $$\mathcal{K}_{r \!
c_p}^p(\psi) = \overline{\bigcup_{n \ge 1} \mathcal{K}_{r \!
c_p}^p(\psi_n)}.$$ We are interested in a cb-embedding
$\mathcal{K}_{r \! c_p}^p(\psi) \to L_p(\mathcal{A})$ preserving
hyperfiniteness. In the construction, we shall use a
noncommutative Poisson random measure. Let us briefly review the
main properties of this notion from \cite{J6} before stating our
result. Let $\mathcal{M}_{sa}^{f}$ stand for the subspace of
self-adjoint elements in $\mathcal{M}$ which are $\psi$-finitely
supported. Let $\mathcal{M}_\pi$ denote the projection lattice of
$\mathcal{M}$. We write $e \perp f$ for orthogonal projections. A
\emph{noncommutative Poisson random measure} is a map $\lambda:
(\mathcal{M}, \psi) \to L_1(\mathcal{A}, \Phi_{\psi})$, where
$(\mathcal{A}, \Phi_{\psi})$ is a noncommutative probability space
and the following conditions hold
\begin{itemize}
\item[\textbf{(i)}] $\lambda: \mathcal{M}_{sa}^f \to
L_1(\mathcal{A})$ is linear.

\item[\textbf{(ii)}] $\Phi_{\psi}(e^{i \lambda(x)}) =
\exp(\psi(e^{ix} - 1))$ for $x \in \mathcal{M}_{sa}^f$.

\item[\textbf{(iii)}] If $e,f \in \mathcal{M}_\pi$ and $e \perp
f$, $\lambda(e \mathcal{M} e)''$ and $\lambda(f \mathcal{M} f)''$
are strongly independent.
\end{itemize}
These properties are not yet enough to characterize $\lambda$, see
below. Let us recall that two von Neumann subalgebras
$\mathcal{A}_1, \mathcal{A}_2$ of $\mathcal{A}$ are called
strongly independent if $a_1 a_2 = a_2 a_1$ and $\Phi_{\psi}(a_1
a_2) = \Phi_{\psi}(a_1) \Phi_{\psi}(a_2)$ for any pair $(a_1,a_2)$
in $\mathcal{A}_1 \times \mathcal{A}_2$. The construction of
$\lambda$ follows by a direct limit argument. Indeed, let us show
how to produce $\lambda: (q_n \mathcal{M} q_n, \psi_n) \to
L_1(\mathcal{A}_n, \Phi_{\psi_n})$. We define
$$\mathcal{A}_n = M_s(q_n \mathcal{M} q_n) = \prod_{k=0}^\infty
\ten_s^k q_n \mathcal{M} q_n,$$ where $\ten _s^k q_n \mathcal{M}
q_n$ denotes the subspace of symmetric tensors in the $k$-fold
tensor product $(q_n \mathcal{M} q_n)_{\ten k}$. In other words,
if $\mathcal{S}_k$ is the symmetric group of permutations of $k$
elements, the space $\ten _s^k q_n \mathcal{M} q_n$ is the range
of the conditional expectation $$\mathcal{E}_k (x_1 \ten x_2 \ten
\cdots \ten x_k) = \frac{1}{k!} \sum_{\pi \in \mathcal{S}_k}
x_{\pi(1)} \ten x_{\pi(2)} \ten \cdots \ten x_{\pi(k)}.$$ Then we
set $$\lambda(x) = (\lambda_k(x))_{k \ge 0} \in M_s(q_n
\mathcal{M} q_n) \quad \mbox{with} \quad \lambda_k(x) =
\sum_{j=1}^k \pi_{ind}^j(x),$$ and properties (i), (ii) and (iii)
hold when working with the state $$\Phi_{\psi_n} \big( (z_k)_{k
\ge 1} \big) = \sum_{k=0}^\infty \frac{\exp(- \psi_n(1))}{k!} \,
\underbrace{\psi_n \ten \cdots \ten \psi_n}_{k \ \mathrm{times}}
(z_k).$$ In the following, it will be important to know the
moments with respect to this state. Given $m \ge 1$, $\Pi(m)$ will
be the set of partitions of $\{1,2, \ldots, m\}$. On the other
hand, given an ordered family $(x_\alpha)_{\alpha \in \Lambda}$ in
$\mathcal{M}$, we shall write $$\prod_{\alpha \in
\Lambda}^{\rightarrow} x_\alpha$$ for the directed product of the
$x_\alpha$'s. Then, the moments are given by the formula
\begin{itemize}
\item[\textbf{(iv)}] $\displaystyle \Phi_{\psi_n} \big(
\lambda(x_1) \lambda(x_2) \cdots \lambda(x_m) \big) =
\sum_{\begin{subarray}{c} \sigma \in \Pi(m) \\ \sigma =
\{\sigma_1, \ldots, \sigma_r\}
\end{subarray}} \prod_{k=1}^{r} \psi_n \Big( \prod_{j \in
\sigma_k}^{\rightarrow} x_j \Big)$.
\end{itemize}
Now we can say that properties (i)-(iv) determine the Poisson
random measure $\lambda$ for any given \emph{n.s.s.f.} weight
$\psi$ in $\mathcal{M}$. According to a uniqueness result from
\cite{J3} which provides a noncommutative form of the Hamburger
moment problem, it turns out that there exists a state preserving
embedding of $M_s(q_{n_1} \mathcal{M} q_{n_1})$ into $M_s(q_{n_2}
\mathcal{M} q_{n_2})$ for $n_1 \le n_2$ and such that the map
$\lambda = \lambda_{n_1}$ constructed for $q_{n_1}$ may be
obtained as a restriction of $\lambda_{n_2}$. This allows to take
direct limits. More precisely, let us define $M_s(\mathcal{M})$ to
be the ultra-weak closure of the direct limit of the $M_s(q_n
\mathcal{M} q_n)$'s. Then, there exists a \emph{n.f.} state
$\Phi_\psi$ on $M_s(\mathcal{M})$ and a map $\lambda$ which
assigns to every self-adjoint operator $x$ (with $\mbox{supp}(x)
\le e$ for some $\psi$-finite projection $e$ in $\mathcal{M}$) a
self-adjoint unbounded operator $\lambda(x)$ affiliated to
$M_s(\mathcal{M})$ and such that $$\Phi_\psi (e^{i \lambda(x)}) =
\exp (\psi(e^{ix}-1)).$$

\begin{theorem} \label{Kp} Let $1\le p\le 2$ and $\psi$ be a n.s.s.f.
weight on $\mathcal{M}$. Then there exists a von Neumann algebra
$\mathcal{A}$, which is hyperfinite when $\mathcal{M}$ is
hyperfinite, and a completely isomorphic embedding
$$\mathcal{K}_{r \! c_p}^p(\psi) \to L_p(\mathcal{A}).$$
\end{theorem}

\dem Let us set the $s$-fold tensor product $$\mathcal{B}_{n,s} =
\Big( L_\infty[0,1] \bar\ten \big[ q_n \mathcal{M} q_n \oplus q_n
\mathcal{M} q_n \big] \Big)_{\otimes s}.$$ Given $s \ge
\mathrm{k}_n$, we define the mapping $\Lambda_{n,s}:
\mathcal{K}_{r \! c_p}^p(\psi_n) \to L_p(\mathcal{B}_{n,s})$ by
\begin{eqnarray*}
\Lambda_{n,s} \big( d_{\varphi_n}^{\frac{1}{2p}} x
d_{\varphi_n}^{\frac{1}{2p}} \big) & = & \sum_{j=1}^s \pi_{ind}^j
\Big( 1_{[0,\mathrm{k}_n/s]} \ten d_{\varphi_n}^\frac{1}{2p} \,
(x,-x) \, d_{\varphi_n}^\frac{1}{2p} \Big) \\ & = & \sum_{j=1}^s
d_{n,s}^\frac{1}{2p} \, \pi_{ind}^j \big( 1_{[0,\mathrm{k}_n/s]}
\ten (x,-x) \big) \, d_{n,s}^\frac{1}{2p},
\end{eqnarray*}
where $d_{n,s}$ is the density associated to the $s$-fold tensor
product state
$$\phi_{n,s} = \Big[ \int_0^1 \cdot \, dt \ten \frac12
(\varphi_n \oplus \varphi_n) \Big]_{\ten s} = \underbrace{\phi_n
\ten \phi_n \ten \cdots \ten \phi_n}_{s \ \mathrm{times}}.$$ If we
tensor $\Lambda_{n,s}$ with the identity map on $\mathrm{M}_m$,
the resulting mapping gives a sum of symmetrically independent
mean-zero random variables over $\mathrm{M}_m$. Therefore, taking
$x \in S_p^m(\mathcal{K}_{r \! c_p}^p(\psi_n))$ and applying Lemma
\ref{Lemma-Araki-Woods}/Remark \ref{Remark-BR}
\begin{eqnarray*}
\big\| \Lambda_{n,s} \big( d_{\varphi_n}^{\frac{1}{2p}} x
d_{\varphi_n}^{\frac{1}{2p}} \big) \big\|_p & \sim & \inf \Big\{
s^{\frac1p} \|a\|_p + s^{\frac12} \big\| \mathsf{E}_{\phi_n} (b
b^*)^{\frac12} \big\|_{S_p^m} + s^\frac12 \big\|
\mathsf{E}_{\phi_n} (c^* c)^{\frac12} \big\|_{S_p^m} \Big\},
\end{eqnarray*}
where $\mathsf{E}_{\phi_n}$ denotes the conditional expectation
$$\mathsf{E}_{\phi_n}: \mathrm{M}_m \Big( L_\infty[0,1] \bar\ten
\big[ q_n \mathcal{M} q_n \oplus q_n \mathcal{M} q_n \big] \Big)
\to \mathrm{M}_m$$ and the infimum runs over all possible
decompositions
\begin{equation} \label{Eq-Decomp-Restriction}
1_{[0, \mathrm{k}_n/s]} \ten d_{\varphi_n}^\frac{1}{2p} \, (x,-x)
\, d_{\varphi_n}^\frac{1}{2p} = a+b+c.
\end{equation}
Multiplying at both sizes by $1_{[0, \mathrm{k}_n/s]} \ten 1$, we
obtain a new decomposition which vanishes over $(\mathrm{k}_n/s,
1]$. Thus, since this clearly improves the infimum above, we may
assume this property in all decompositions considered. Moreover,
we claim that we can also restrict the infimum above to those
decompositions $a+b+c$ which are constant on $[0,
\mathrm{k}_n/s]$. Indeed, given any decomposition of the form
\eqref{Eq-Decomp-Restriction} we take averages at both sizes and
produce another decomposition $a_0 + b_0 + c_0$ given by the
relations $$\big( a_0, b_0, c_0 \big) = 1_{[0, \mathrm{k}_n/s]}
\ten \frac{s}{\mathrm{k}_n} \int_0^{\frac{\mathrm{k}_n}{s}} \big(
a(t), b(t), c(t) \big) \, dt.$$ Then, our claim is a consequence
of the inequalities
\begin{eqnarray*}
\|a_0\|_p & \le & \|a\|_p, \\ \big\| \mathsf{E}_{\phi_n} (b_0
b_0^*)^{\frac12} \big\|_{S_p^m} & \le & \big\| \mathsf{E}_{\phi_n}
(b b^*)^{\frac12} \big\|_{S_p^m},
\\ \big\| \mathsf{E}_{\phi_n} (c_0^* c_0)^{\frac12} \big\|_{S_p^m}
& \le & \big\| \mathsf{E}_{\phi_n} (c^* c)^{\frac12}
\big\|_{S_p^m}.
\end{eqnarray*}
The fist one is justified by means of the inequality
$$\lambda^{\frac1p} \Big\| \frac{1}{\lambda} \int_0^\lambda a(t)
\, dt \Big\|_{L_p(\mathcal{M})} \le \|a\|_{L_p([0,\lambda]
\bar\ten \mathcal{M})},$$ which follows easily  by complex
interpolation. The two other inequalities arise as a consequence
of Kadison's inequality $\mathsf{E}(x) \mathsf{E}(x^*) \le
\mathsf{E}(xx^*)$ applied to the conditional expectation
$$\mathsf{E} = \frac{1}{\lambda} \int_0^\lambda \cdot \, dt.$$
Our considerations allow us to assume
$$\big( a,b,c \big) = 1_{[0, \mathrm{k}_n/s]} \ten \big( x_p,
x_r, x_c \big)$$ for some $x_p, x_r, x_c \in S_p^m(L_p(q_n
\mathcal{M} q_n))$. This gives rise to
$$\big\| \Lambda_{n,s} \big( d_{\varphi_n}^{\frac{1}{2p}} x
d_{\varphi_n}^{\frac{1}{2p}} \big) \big\|_p \sim \inf \Big\{
\mathrm{k}_n^{\frac1p} \|x_p\|_p + \mathrm{k}_n^\frac12 \big\|
\mathsf{E}_{\varphi_n} (x_r x_r^*)^{\frac12} \big\|_{S_p^m} +
\mathrm{k}_n^{\frac12} \big\| \mathsf{E}_{\varphi_n} (x_c^*
x_c)^{\frac12} \big\|_{S_p^m} \Big\},$$ where the infimum runs
over all possible decompositions
$$d_{\varphi_n}^{\frac{1}{2p}} x d_{\varphi_n}^{\frac{1}{2p}} =
x_p + x_r + x_c.$$ This shows that $\Lambda_{n,s}: \mathcal{K}_{r
\! c_p}^p(\psi_n) \to L_p(\mathcal{B}_{n,s})$ is a completely
isomorphic embedding with constants independent of $n$ or $s$. We
are not ready yet to take direct limits. Before that, we use the
algebraic central limit theorem to identify the moments in the
limit as $s \to \infty$. To calculate the joint moments we set
$$\zeta_n = \frac12 (\varphi_n \oplus \varphi_n) \quad \mbox{and}
\quad \zeta_{n,s} = \underbrace{\zeta_n \ten \zeta_n \ten \cdots
\ten \zeta_n}_{s \ \mathrm{times}}$$ and recall that the map
$\Lambda_{n,s}$ corresponds to
$$u_{n,s} (x) = \sum_{j=1}^s \pi_{ind}^j \Big(
1_{[0,\mathrm{k}_n/s]} \ten (x,-x) \Big).$$ Then, the joint
moments are given by
\begin{eqnarray*}
\lefteqn{\phi_{n,s} \big( u_{n,s}(x_1) \cdots u_{n,s}(x_m) \big)}
\\ & = & \sum_{j_1, j_2, \ldots, j_m = 1}^s \int_{[0,1]^s}
\prod_{i=1}^m \pi_{ind}^{j_i}(1_{[0,\mathrm{k}_n/s]})(t) \, dt \
\zeta_{n,s} \Big( \prod_{1 \le i \le m}^{\rightarrow}
\pi_{ind}^{j_i}(x_i,-x_i) \Big) \\ & = & \sum_{\begin{subarray}{c}
\sigma \in \Pi(m) \\ \sigma = \{\sigma_1, \ldots, \sigma_r \}
\end{subarray}} \sum_{(j_1, \ldots, j_m) \sim \sigma} \Big(
\frac{\mathrm{k}_n}{s} \Big)^r \prod_{k=1}^r \zeta_n \Big[ \Big(
\prod_{i \in \sigma_k}^{\rightarrow} x_i, (-1)^{|\sigma_k|}
\prod_{i \in \sigma_k}^{\rightarrow} x_i \Big) \Big],
\end{eqnarray*}
where $|\sigma_k|$ denotes the cardinality of $\sigma_k$ and we
write $(j_1, \ldots, j_m) \sim \sigma$ when $j_a = j_b$ if and
only if there exits $1 \le k \le r$ such that $j_a, j_b \in
\sigma_k$. Therefore, recalling that $\zeta_n = \frac12 (\varphi_n
\oplus \varphi_n)$, the only partitions which contribute to the
sum above are the even partitions satisfying $|\sigma_k| \in 2 \N$
for $1 \le k \le r$. Let us write $\Pi_e(m)$ for the set of even
partitions. Then, using $\psi_n = \mathrm{k}_n \varphi_n$ we
deduce
\begin{eqnarray*}
\phi_{n,s} \big( u_{n,s}(x_1) \cdots u_{n,s}(x_m) \big) & = &
\sum_{\begin{subarray}{c} \sigma \in \Pi_e(m) \\ \sigma =
\{\sigma_1, \ldots, \sigma_r \} \end{subarray}} \frac{|\{(j_1,
\ldots, j_m) \sim \sigma\}|}{s^r} \prod_{k=1}^r \psi_n \Big(
\prod_{i \in \sigma_k}^{\rightarrow} x_i \Big) \\ & = &
\sum_{\begin{subarray}{c} \sigma \in \Pi_{e}(m)
\\ \sigma = \{\sigma_1, \ldots , \sigma_r\} \end{subarray}}
\frac{s!}{s^r (s-r)!} \prod_{k=1}^r \psi_n \Big( \prod_{i \in
\sigma_k}^{\rightarrow} x_i \Big).
\end{eqnarray*}
Therefore, taking limits $$\lim_{s \to \infty} \phi_{n,s} \big(
u_{n,s}(x_1) \cdots u_{n,s}(x_m) \big) = \sum_{\begin{subarray}{c}
\sigma \in \Pi_{e}(m)
\\ \sigma = \{\sigma_1, \ldots , \sigma_r\} \end{subarray}}
\prod_{k=1}^r \psi_n \Big( \prod_{i \in \sigma_k}^{\rightarrow}
x_i \Big).$$ These moments coincide with the moments of the
Poisson random process $$\lambda: (q_n \mathcal{M} q_n, \psi_n)
\to L_1(M_s(q_n \mathcal{M} q_n), \Phi_{\psi_n}).$$ Hence, the
noncommutative version of the Hamburger moment problem from
\cite{J3} provides a state preserving homomorphism between the von
Neumann algebra which generate the operators $$\Big\{ e^{i
u_{n,s}(x)} \, \big| \ x \in q_n \mathcal{M} q_n, \, s \ge 1
\Big\}$$ and the von Neumann subalgebra of $M_s(q_n \mathcal{M}
q_n)$ generated by
$$\Big\{ e^{i \lambda(x)} \, \big| \ x \in (q_n
\mathcal{M} q_n)_{sa}^f \Big\}.$$ In particular, taking
$\mathcal{A}_n = M_s(q_n \mathcal{M} q_n)$
\begin{eqnarray*}
\lefteqn{\big\| d_{\psi_n}^{\frac{1}{2p}} x
d_{\psi_n}^{\frac{1}{2p}}
\big\|_{S_p^m(\mathcal{K}_{r \! c_p}^p(\psi_n))}} \\
& \sim & \lim_{s \to \infty} \big\| \Lambda_{n,s} \big(
d_{\psi_n}^\frac{1}{2p} x d_{\psi_n}^\frac{1}{2p} \big)
\big\|_{S_p^m(L_p(\mathcal{B}_{n,s}))} = \big\|
d_{\Phi_{\psi_n}}^{\frac{1}{2p}} \lambda(x)
d_{\Phi_{\psi_n}}^{\frac{1}{2p}}
\big\|_{S_p^m(L_p(\mathcal{A}_n))}.
\end{eqnarray*}
Now, we use from \cite{J6} that $$\big(
M_s(\mathcal{M}),\Phi_{\psi} \big) = \overline{\bigcup_{n \ge 1}
\big( M_s(q_n \mathcal{M} q_n), \Phi_{\psi_n} \big)}$$ exists.
Therefore, the map $$\Lambda \big( d_{\psi}^{\frac{1}{2p}} x
d_{\psi}^{\frac{1}{2p}} \big) = d_{\Phi_{\psi}}^{\frac{1}{2p}}
\lambda(x) d_{\Phi_{\psi}}^{\frac{1}{2p}}$$ extends to a complete
embedding
$$\mathcal{K}_{r \! c_p}^p(\psi) = \limm_n \mathcal{K}_{r \!
c_p}^p (\psi_n) \to L_p(M_s(\mathcal{M})).$$ Moreover, if
$\mathcal{M}$ is hyperfinite so is $M_s(q_n \mathcal{M} q_n)$ and
hence the limit $M_s(\mathcal{M})$. \fin

\begin{corollary} \label{Corollary-Kp2-Lp}
Let $1\le p\le 2$ and $\psi$ be a n.s.s.f. weight on
$\mathcal{M}$. Then there exists a von Neumann algebra
$\mathcal{A}$, which is hyperfinite when $\mathcal{M}$ is
hyperfinite, and a completely isomorphic embedding
$$\mathcal{K}_{p,2}(\psi) \to L_p(\mathcal{A}).$$
\end{corollary}

\dem Let us set $$\mathcal{RB}_{n,s} = \Big( L_\infty[0,1]
\bar\ten \big[ (\mathcal{R} \bar\ten q_n \mathcal{M} q_n) \oplus
(\mathcal{R} \bar\ten q_n \mathcal{M} q_n) \big] \Big)_{\ten s}.$$
By Proposition \ref{Prop-k4k3} and Theorem \ref{Kp}, the map
$$\Lambda_{n,s} \circ \rho_n: \mathcal{K}_{p,2}(\psi_n) \to
\mathcal{K}_{r \! c_p}^p(\phi \ten \psi_n) \to
L_p(\mathcal{RB}_{n,s})$$ provides a complete isomorphism with
constant independent of $n$ and $s$. Using the algebraic central
limit theorem to take limits in $s$ and the noncommutative version
of the Hamburger moment problem one more time, we obtain a
complete embedding $$\Lambda_n \circ \rho_n:
\mathcal{K}_{p,2}(\psi_n) \to \mathcal{K}_{r \! c_p}^p(\phi \ten
\psi_n) \to L_p \big( M_s(\mathcal{R} \bar\ten q_n \mathcal{M}
q_n) \big).$$ Taking direct limits we obtain a cb-embedding which
preserves hyperfiniteness. \fin

\begin{corollary} \label{caoarao}
$S_q$ cb-embeds into $L_p(\mathcal{A})$ for some hyperfinite
factor $\mathcal{A}$.
\end{corollary}

\dem According to the complete isometry $S_q = C_q \otimes_h R_q$
and Remark \ref{Remark-Graph_Sq}, it suffices to embed the first
term $\mathcal{Z}_1 \otimes \mathcal{Z}_2$ on the right of
\eqref{Eq-6terms} into $L_p(\mathcal{A})$ for some hyperfinite von
Neumann algebra $\mathcal{A}$. However, following Part II of the
proof of Theorem \ref{Theorem-QS}, we know that $\mathcal{Z}_1
\otimes \mathcal{Z}_2$ embeds completely isomorphically into
$\mathcal{K}_{p,2}(\psi)$, where $\psi$ denotes some
\emph{n.s.s.f.} weight on $\mathcal{B}(\ell_2)$. Therefore, the
assertion follows from Corollary \ref{Corollary-Kp2-Lp}. This
completes the proof. \fin

\begin{remark}
\emph{Theorem \ref{Theorem-QS} easily generalizes to the context
of Corollary \ref{caoarao}. More precisely, given operator spaces
$\mathrm{X}_1 \in \mathcal{QS}(C_p \oplus_2 \mathrm{OH})$ and
$\mathrm{X}_2 \in \mathcal{QS}(R_p \oplus_2 \mathrm{OH})$ and
combining the techniques applied so far, it is rather easy to find
a hyperfinite type $\mathrm{III}$ factor $\mathcal{A}$ and a
completely isomorphic embedding
$$\mathrm{X}_1 \ten_h \mathrm{X}_2 \to L_p(\mathcal{A}).$$ }
\end{remark}

\begin{remark} \label{estfindim}
\emph{In contrast with Corollary \ref{SqQWEP}, where free products
are used, the complete embedding of $S_q$ into $L_p$ given in
Corollary \ref{caoarao} provides estimates on the dimension of
$\mathcal{A}$ in the cb-embedding
$$S_q^m \to L_p(\mathcal{A}).$$ Indeed, a quick look at our
construction shows that $$S_q^m = C_q^m \ten_h R_q^m \to
\mathcal{K}_{p,2}(\psi_n) \to L_p \big( \mathcal{A}_{ind}^n;
\mathrm{OH}_{\mathrm{k}_n} \big) = L_p \big( \mathrm{M}_n^{\ten
\mathrm{k}_n}; \mathrm{OH}_{\mathrm{k}_n} \big),$$ with $n \sim m
\log m$, see \cite{J3} for this last assertion. This chain
essentially follows from Remark \ref{Remark-Graph_Sq}, Lemma
\ref{Lemma-Direct-Sum} and Proposition \ref{trick}. On the other
hand, given any parameter $\gamma > 0$ and according once more to
to \cite{J3}, we know that $\mathrm{OH}_{\mathrm{k}_n}$ embeds
completely isomorphically into $S_p^{\mathrm{w}_n}$ for
$\mathrm{w}_n = \mathrm{k}_n^{\gamma \mathrm{k}_n}$ with absolute
constants depending only on $\gamma$ and that $\mathrm{k}_n \sim
n^{\alpha_p}$. Combining the embeddings mentioned so far, we have
found a complete embedding
$$S_q^m \to S_p^M \quad \mbox{with} \quad M \sim m^{\beta_p
m^{\alpha_p}}.$$}
\end{remark}

\subsection{Embedding for general von Neumann algebras}

We conclude this paper with the proof of our main embedding result
in full generality. We shall need to extend our definition
\eqref{Eq-RCLp} to the case $1 \le q \le 2$. This is easily done
as follows
\begin{equation}
\begin{array}{rcl}
L_q^r(\mathcal{M}) & = & \big[ L_1(\mathcal{M}),
L_2^r(\mathcal{M})\big]_{\frac{2}{q'}}, \\ [5pt]
L_q^c(\mathcal{M}) & = & \big[ L_1(\mathcal{M}),
L_2^c(\mathcal{M})\big]_{\frac{2}{q'}}.
\end{array}
\end{equation}
The space $L_q$ is given by complex interpolation. Therefore we
will encode complex interpolation in a suitable graph. This
follows Pisier's approach \cite{P4} to the main result in
\cite{J2}. Indeed, consider a fixed $0<\theta<1$ and let
$\mu_\theta = (1-\theta) \mu_0 + \theta \mu_1$ be the harmonic
measure on the boundary of the strip $\mathcal{S}$ associated to
the point $z= \theta$, as defined at the beginning of Section
\ref{NewSect1}. We consider the space
$$\mathcal{H}_2 = \Big\{ \big( f_{|_{\partial_0}},
f_{|_{\partial_1}} \big) \, \big| \ f: \mathcal{S} \to \C \ \mbox{
analytic} \Big\} \subset L_2(\partial_0) \oplus L_2(\partial_1).$$
We need operator-valued versions of this space given by subspaces
$$\mathcal{H}_{2p',2}^r(\mathcal{M},\theta) \subset \Big(
L_2^{c_{p'}}(\partial_0) \ten_h L_{2p'}^r(\mathcal{M}) \Big)
\oplus \Big( L_2^{oh}(\partial_1) \ten_h L_4^r(\mathcal{M}) \Big)
= \mathcal{O}_{p',0}^r \oplus \mathcal{O}_{p',1}^r,$$
$$\mathcal{H}_{2p',2}^c(\mathcal{M},\theta) \subset \Big(
L_{2p'}^c(\mathcal{M}) \ten_h L_2^{r_{p'}}(\partial_0) \Big)
\oplus \Big( L_4^c(\mathcal{M}) \ten_h L_2^{oh}(\partial_1) \Big)
= \mathcal{O}_{p',0}^c \oplus \mathcal{O}_{p',1}^c.$$ More
precisely, if $\mathcal{M}$ comes equipped with a \emph{n.s.s.f.}
weight $\psi$ and $d_{\psi}$ denotes the associated density,
$\mathcal{H}_{2p',2}^r(\mathcal{M},\theta)$ is the subspace of all
pairs $(f_0,f_1)$ of functions in $\mathcal{O}_{p',0}^r \oplus
\mathcal{O}_{p',1}^r$ such that for every scalar-valued analytic
function $g: \mathcal{S} \to \C$ (extended non-tangentially to the
boundary) with $g(\theta)=0$, we have
$$(1-\theta) \int_{\partial_0} g(z) \, d_\psi^{\frac14 -
\frac{1}{2p'}} f_0(z) \, d\mu_0(z) + \theta \int_{\partial_1} g(z)
f_1(z) \, d\mu_1(z) = 0.$$ Similarly, the condition on
$\mathcal{H}_{2p',2}^c(\mathcal{M},\theta)$ is $$(1-\theta)
\int_{\partial_0} g(z) f_0(z) \, d_\psi^{\frac14 - \frac{1}{2p'}}
\, d\mu_0(z) + \theta \int_{\partial_1} g(z) f_1(z) \, d\mu_1(z) =
0.$$ We shall also need to consider the subspaces
\begin{eqnarray*}
\mathcal{H}_{r,0} & = & \Big\{ (f_0,f_1) \in
\mathcal{H}_{2p',2}^r(\mathcal{M},\theta) \, \big| \ (1-\theta)
\int_{\partial_0} d_{\psi}^{\frac14 - \frac{1}{2p'}} f_0 d\mu_0 +
\theta \int_{\partial_1} f_1 d\mu_1 = 0 \Big\},
\\ \mathcal{H}_{c,0} & = & \Big\{ (f_0,f_1) \in
\mathcal{H}_{2p',2}^c(\mathcal{M},\theta) \, \big| \ (1-\theta)
\int_{\partial_0} f_0 d_{\psi}^{\frac14 - \frac{1}{2p'}} d\mu_0 +
\theta \int_{\partial_1} f_1 d\mu_1 = 0 \Big\}.
\end{eqnarray*}

\begin{remark} \label{Remark-Bar-Ten}
\emph{In order to make all the forthcoming duality arguments work,
we need to introduce a slight modification of these spaces for
$p=1$. Indeed, in that case the spaces defined above must be
regarded as subspaces of
\begin{eqnarray*}
\mathcal{H}_{\infty,2}^r(\mathcal{M},\theta) & \subset & \Big(
L_2^{c}(\partial_0) \bar\ten \mathcal{M} \Big) \oplus \Big(
L_2^{oh}(\partial_1) \ten_h L_4^r(\mathcal{M}) \Big),
\\ \mathcal{H}_{\infty,2}^c(\mathcal{M},\theta) & \subset &
\Big( \mathcal{M} \bar\ten L_2^{r}(\partial_0) \Big) \oplus \Big(
L_4^c(\mathcal{M}) \ten_h L_2^{oh}(\partial_1) \Big).
\end{eqnarray*}
The von Neumann algebra tensor product used above is the weak
closure of the minimal tensor product, which in this particular
case coincides with the Haagerup tensor product since we have
either a column space on the left or a row space on the right. In
particular, we just take the closure in the weak operator
topology.}
\end{remark}

\begin{lemma} \label{step1}
Let $\mathcal{M}$ be a finite von Neumann algebra equipped with a
n.f. state $\varphi$ and let $d_{\varphi}$ be the associated
density. If $2 \le q' < p'$ and $\frac{1}{2q'} =
\frac{1-\theta}{2p'} + \frac{\theta}{4}$, we have complete
contractions
$$u_r: d_{\varphi}^{\frac{1}{2q'}} x
\in L_{2q'}^r(\mathcal{M}) \mapsto \Big( 1 \ten
d_{\varphi}^{\frac{1}{2p'}} x, 1 \ten d_{\varphi}^{\frac14} x
\Big) + \mathcal{H}_{r,0} \in
\mathcal{H}_{2p',2}^r(\mathcal{M},\theta) / \mathcal{H}_{r,0},$$
$$u_c: x d_{\varphi}^{\frac{1}{2q'}} \in
L_{2q'}^c(\mathcal{M}) \mapsto \Big( x d_{\varphi}^{\frac{1}{2p'}}
\ten 1, x d_{\varphi}^{\frac14} \ten 1 \Big) + \mathcal{H}_{c,0}
\in \mathcal{H}_{2p',2}^c(\mathcal{M},\theta) /
\mathcal{H}_{c,0}.$$
\end{lemma}

\dem By symmetry, it suffices to consider the column case. Let $x$
be an element in $\mathrm{M}_m (L_{2q'}^c(\mathcal{M}))$ of norm
less than $1$. According to our choice of $0 < \theta < 1$, we
find that $$\mathrm{M}_m (L_{2q'}^c(\mathcal{M})) = \big[
\mathrm{M}_m (L_{2p'}^c(\mathcal{M})), \mathrm{M}_m
(L_4^c(\mathcal{M})) \big]_\theta.$$ Thus, there exists $f:
\mathcal{S} \to \mathrm{M}_m(\mathcal{M})$ analytic such that
$f(\theta) = x$ and
$$\max \Big\{ \sup_{z \in \partial_0}
\big\| f(z) d_{\varphi}^{\frac{1}{2p'}} \big\|_{\mathrm{M}_m
(L_{2p'}^c(\mathcal{M}))}, \sup_{z \in
\partial_1} \big\| f(z) d_{\varphi}^\frac14 \big\|_{\mathrm{M}_m
(L_4^c(\mathcal{M}))} \Big\} \le 1.$$ If $1 \le s \le \infty$ and
$j \in \{0,1\}$, we claim that
\begin{equation} \label{Eq-Cla1}
\big\| f_{\mid_{\partial_j}} d_{\varphi}^{\frac{1}{2s}}
\big\|_{\mathrm{M}_m(L_{2s}^c(\mathcal{M}) \otimes_h
L_2^{r_{s}}(\partial_j))} \le \sup_{z \in
\partial_j} \big\| f(z) d_{\varphi}^{\frac{1}{2s}}
\big\|_{\mathrm{M}_m(L_{2s}^c(\mathcal{M}))}.
\end{equation}
Before proving our claim, let us finish the proof. Taking $f_j =
f_{\mid_{\partial_j}}$, we have
$$\big( f_0 d_{\varphi}^{\frac{1}{2p'}}, f_1 d_{\varphi}^{\frac{1}{4}}
\big) - \Big( x d_{\varphi}^{\frac{1}{2p'}} \ten 1, x
d_{\varphi}^{\frac14} \ten 1 \Big) \in \mathcal{H}_{c,0}.$$ Indeed
by analyticity, we have
$$(1-\theta) \int_{\partial_0} f_0 \, d_{\varphi}^{\frac{1}{2p'}}
d_{\varphi}^{\frac14 - \frac{1}{2p'}} \, d\mu_0 + \theta
\int_{\partial_1} f_1 \, d_{\varphi}^{\frac14} \, d\mu_1 =
\int_{\partial \mathcal{S}} f d_{\varphi}^{\frac14} d\mu_\theta =
f(\theta) d_{\varphi}^\frac14 = x d_{\varphi}^\frac14.$$ This
implies from \eqref{Eq-Cla1} applied to $(s,j)=(p',0)$ and
$(s,j)=(2,1)$ that
$$\Big\| u_c \big( x d_{\varphi}^{\frac{1}{2q'}} \big)
\Big\|_{\mathrm{M}_m(\mathcal{H}_{2p',2}^c/ \mathcal{H}_{c,0})}
\le \Big\| \big( f_0 d_{\varphi}^{\frac{1}{2p'}}, f_1
d_{\varphi}^{\frac{1}{4}} \big)
\Big\|_{\mathrm{M}_m(\mathcal{H}_{2p',2}^c)} \le 1.$$ Hence, it
remains to prove our claim \eqref{Eq-Cla1}. We must show that the
identity map $L_{\infty} \big( \partial_j; L_{2s}^c(\mathcal{M})
\big) \to L_{2s}^c(\mathcal{M}) \ten_h L_2^{r_{s}}(\partial_j)$ is
a complete contraction. By complex interpolation, we have
\begin{eqnarray*}
L_{\infty} \big( \partial_j; L_{2s}^c(\mathcal{M}) \big) & = &
\big[ L_{\infty} \big( \partial_j; \mathcal{M} \big), L_{\infty}
\big( \partial_j; L_{2}^c(\mathcal{M}) \big) \big]_{\frac{1}{s}},
\\ L_{2s}^c(\mathcal{M}) \ten_h L_2^{r_{s}}(\partial_j) & = &
\big[ \mathcal{M} \ten_h L_2^r(\partial_j), L_2^c(\mathcal{M})
\ten_h L_2^c(\partial_j) \big]_{\frac{1}{s}} \\ & = & \big[
\mathcal{M} \ten_{\mathrm{min}} L_2^r(\partial_j),
L_2^c(\mathcal{M}) \ten_{\mathrm{min}} L_2^c(\partial_j)
\big]_{\frac{1}{s}}.
\end{eqnarray*}
In other words, we must study the identity mappings
\begin{eqnarray*}
\mathcal{M} \ten_{\mathrm{min}} L_{\infty} (\partial_j) & \to &
\mathcal{M} \ten_{\mathrm{min}} L_2^r (\partial_j), \\
L_{2}^c(\mathcal{M}) \ten_{\mathrm{min}} L_{\infty} (\partial_j) &
\to & L_2^c(\mathcal{M}) \ten_{\mathrm{min}} L_2^c(\partial_j).
\end{eqnarray*}
However, this automatically reduces to see that we have complete
contractions
\begin{eqnarray*}
L_{\infty} (\partial_j) & \to & L_2^r (\partial_j), \\
L_{\infty} (\partial_j) & \to & L_2^c(\partial_j).
\end{eqnarray*}
Therefore it suffices  to observe that
\begin{eqnarray*}
\|f\|_{\mathrm{M}_m(L_2^r(\partial_j))}^2 = \Big\|
\int_{\partial_j} f f^* d\mu_j \Big\|_{\mathrm{M}_m} & \le &
\mu_j(\partial_j) \, \sup_{z \in
\partial_j} \|f(z)\|_{\mathrm{M}_m}^2 = \|f\|_{\mathrm{M}_m
(L_\infty(\partial_j))}^2, \\
\|f\|_{\mathrm{M}_m(L_2^c(\partial_j))}^2 = \Big\|
\int_{\partial_j} f^* f d\mu_j \Big\|_{ \mathrm{M}_m} & \le &
\mu_j(\partial_j) \, \sup_{z \in
\partial_j} \|f(z)\|_{\mathrm{M}_m}^2 = \|f\|_{\mathrm{M}_m
(L_\infty(\partial_j))}^2.
\end{eqnarray*}
This completes the proof for $1 < p \le 2$. In the case $p = 1$,
we have overlooked the fact that the definition of
$\mathcal{H}_{\infty,2}^c(\mathcal{M},\theta)$ (see Remark
\ref{Remark-Bar-Ten} above) is slightly different. The only
consequence of this point is that we also need the inequality
$$\big\| f_{\mid_{\partial_0}} \big\|_{\mathrm{M}_m(\mathcal{M}
\bar\otimes L_2^{r}(\partial_0))} \le \sup_{z \in
\partial_0} \|f(z)\|_{\mathrm{M}_m(\mathcal{M})}.$$ However,
this is proved once more as above. The proof is complete. \fin

Lemma \ref{step1} is closely related to Lemma
\ref{Lemma-Motivation} and a similar result holds on the preduals.
More precisely, we begin by defining the operator-valued Hardy
spaces which arise as subspaces
\begin{eqnarray*}
\mathcal{H}_{2p,2}^c(\mathcal{M},\theta) & \subset & \Big(
L_2^{c_{p}}(\partial_0) \ten_h L_{\frac{2p}{p+1}}^c(\mathcal{M})
\Big) \oplus \Big( L_2^{oh}(\partial_1) \ten_h
L_{\frac43}^c(\mathcal{M}) \Big), \\
\mathcal{H}_{2p,2}^r(\mathcal{M},\theta) & \subset & \Big(
L_{\frac{2p}{p+1}}^r(\mathcal{M}) \ten_h L_2^{r_{p}}(\partial_0)
\Big) \oplus \Big( L_{\frac43}^r(\mathcal{M}) \ten_h
L_2^{oh}(\partial_1) \Big),
\end{eqnarray*}
formed by pairs $(f_0,f_1)$ respectively satisfying
\begin{eqnarray*}
(1-\theta) \int_{\partial_0} g(z) f_0(z) \, d\mu_0(z) + \theta
\int_{\partial_1} g(z) d_{\varphi}^{\frac{p+1}{2p} - \frac{3}{4}}
f_1(z) \, d\mu_1(z) & = & 0, \\ (1-\theta) \int_{\partial_0} g(z)
f_0(z) \, d\mu_0(z) + \theta \int_{\partial_1} g(z) f_1(z)
d_{\varphi}^{\frac{p+1}{2p} - \frac{3}{4}} \, d\mu_1(z) & = & 0,
\end{eqnarray*}
for all scalar-valued analytic function $g: \mathcal{S} \to \C$
(extended non-tangentially to the boundary) with $g(\theta)=0$.
The subspaces $\mathcal{H}_{r,0}'$ and $\mathcal{H}_{c,0}'$ are
defined accordingly. In other words, we have
\begin{eqnarray*}
\mathcal{H}_{c,0}' & = & \Big\{ (f_0,f_1) \in
\mathcal{H}_{2p,2}^c(\mathcal{M},\theta) \, \big| \ (1-\theta)
\int_{\partial_0} f_0 \, d\mu_0 + \theta \int_{\partial_1}
d_{\varphi}^{\frac{p+1}{2p} - \frac{3}{4}} f_1 \, d\mu_1 = 0
\Big\}, \\ \mathcal{H}_{r,0}' & = & \Big\{ (f_0,f_1) \in
\mathcal{H}_{2p,2}^r(\mathcal{M},\theta) \, \big| \ (1-\theta)
\int_{\partial_0} f_0 \, d\mu_0 + \theta \int_{\partial_1} f_1
d_{\varphi}^{\frac{p+1}{2p} - \frac{3}{4}} \, d\mu_1 = 0 \Big\}.
\end{eqnarray*}

\begin{lemma} \label{LemmadeDualidad}
Let $\mathcal{M}$ be a finite von Neumann algebra equipped with a
n.f. state $\varphi$ and let $d_{\varphi}$ be the associated
density. Taking the same values for $p,q$ and $\theta$ as above,
we have complete contractions $$w_r: d_{\varphi}^{\frac{q+1}{2q}}
x \in L_{\frac{2q}{q+1}}^c(\mathcal{M}) \mapsto \Big( 1 \ten
d_{\varphi}^{\frac{p+1}{2p}} x, 1 \ten d_{\varphi}^{\frac34} x
\Big) + \mathcal{H}_{c,0}' \in
\mathcal{H}_{2p,2}^c(\mathcal{M},\theta) / \mathcal{H}_{c,0}',$$
$$w_c: x d_{\varphi}^{\frac{q+1}{2q}} \in
L_{\frac{2q}{q+1}}^r(\mathcal{M}) \mapsto \Big( x
d_{\varphi}^{\frac{p+1}{2p}} \ten 1, x d_{\varphi}^{\frac34} \ten
1 \Big) + \mathcal{H}_{r,0}' \in
\mathcal{H}_{2p,2}^r(\mathcal{M},\theta) / \mathcal{H}_{r,0}'.$$
\end{lemma}

\dem If $1 \le s \le \infty$ and $\mathcal{M}_m =
\mathrm{M}_m(\mathcal{M})$, we have
\begin{eqnarray}
\label{EqAmalgRow} S_1^m \big( L_{\frac{2s}{s+1}}^r(\mathcal{M})
\big) & = & S_1^m \Big( \big[ L_1(\mathcal{M}), L_2^r(\mathcal{M})
\big]_{\frac{1}{s'}} \Big) \\ \nonumber & = & \Big[ S_1^m \big(
L_1(\mathcal{M}) \big), S_2^m \big( L_2(\mathcal{M}) \big) S_2^m
\Big]_{\frac{1}{s'}} \\ \nonumber & = & S_{\frac{2s}{s+1}}^m \big(
L_{\frac{2s}{s+1}}(\mathcal{M}) \big) S_{2s'}^m =
L_{\frac{2s}{s+1}}(\mathcal{M}_m) S_{2s'}^m.
\end{eqnarray}
Indeed, the second identity follows from \eqref{Eq-Schatten-RC}
and duality (alternatively, one may argue directly as we did in
Lemma \ref{Lemma-IsoK3Cond}), while the third one follows from
Theorem \ref{TheoApp1}. Now, let us consider an element of norm
less than $1$ $$x d_{\varphi}^{\frac{q+1}{2q}} \in S_1^m \big(
L_{\frac{2q}{q+1}}^r(\mathcal{M}) \big).$$ The isometry above
provides a factorization $$x d_{\varphi}^{\frac{q+1}{2q}} = \alpha
\beta \gamma \delta \in L_{\frac43}(\mathcal{M}_m)
L_{\rho_1}(\mathcal{M}_m) S_{2p'}^m S_{\rho_2}^m$$ with $\alpha,
\beta, \gamma, \delta$ in the unit balls of their respective
spaces with $$\frac{1}{\rho_1} = \frac{q+1}{2q} - \frac34 \quad
\mbox{and} \quad \frac{1}{\rho_2} = \frac{1}{2q'} -
\frac{1}{2p'}.$$ Moreover, by polar decomposition and
approximation, we may assume that $\beta$ and $\delta$ are
strictly positive elements. In particular, motivated by the
complex interpolation isometry $$L_{\frac{2q}{q+1}}(\mathcal{M}_m)
S_{2q'}^m = \Big[ L_{\frac{2p}{p+1}}(\mathcal{M}_m) S_{2p'}^m,
L_{\frac{4}{3}}(\mathcal{M}_m) S_{4}^m \Big]_\theta = \big[
\mathrm{X}_0, \mathrm{X}_1 \big]_\theta,$$ we take $(\beta_\theta,
\delta_\theta) = (\beta^{1/(1-\theta)},\delta^{1/\theta})$ and
define $$f: z \in \mathcal{S} \mapsto \alpha \beta_{\theta}^{1-z}
\gamma \delta_\theta^z \in \mathrm{X}_0 + \mathrm{X}_1.$$ Since
$f$ is analytic and $f(\theta) = x d_\varphi^{\frac{q+1}{2q}}$, we
conclude $$\big( f_{\mid_{\partial_0}}, f_{\mid_{\partial_1}}
\big) \in \big( x d_\varphi^{\frac{p+1}{2p}} \ten 1, x
d_\varphi^{\frac{3}{4}} \ten 1 \big) + \mathcal{H}_{r,0}'.$$
Therefore, taking $f_{\mid_{\partial_j}} = f_j$ we obtain the
following estimate
\begin{eqnarray*}
\lefteqn{\Big\| w_c \big( x d_{\varphi}^{\frac{q+1}{2q}} \big)
\Big\|_{S_1^m(\mathcal{H}_{2p,2}^r/ \mathcal{H}_{r,0}')}} \\ & \le
& \max \Big\{ \big\| f_0 \big\|_{S_1^m
(L_{\frac{2p}{p+1}}^r(\mathcal{M}) \ten_h L_2^{r_p}(\partial_0))},
\big\| f_1 \big\|_{S_1^m (L_{\frac{4}{3}}^r(\mathcal{M}) \ten_h
L_2^{oh}(\partial_1))} \Big\} \\ & = & \max \Big\{ \big\| \alpha
\beta_\theta \gamma \big\|_{S_1^m
(L_{\frac{2p}{p+1}}^r(\mathcal{M}))}, \big\| \alpha \gamma
\delta_\theta \big\|_{S_1^m (L_{\frac{4}{3}}^r(\mathcal{M}))}
\Big\} \le 1,
\end{eqnarray*}
where the last identity follows from \eqref{EqAmalgRow} and the
fact that $(\alpha \beta_\theta, \gamma \delta_\theta)$ are in the
unit balls of $L_{2p/p+1}(\mathcal{M}_m)$ and $S_4^m$
respectively. The assertion for the mapping $w_r$ is proved
similarly. The proof is complete. \fin

\begin{proposition} \label{Proposition-CBIsomorphisms-Quo}
Let $\mathcal{M}$ be a finite von Neumann algebra equipped with a
n.f. state $\varphi$ and let $d_{\varphi}$ be the associated
density. If $2 \le q' < p'$ and $\frac{1}{2q'} =
\frac{1-\theta}{2p'} + \frac{\theta}{4}$, we have complete
isomorphisms
$$u_r: d_{\varphi}^{\frac{1}{2q'}} x
\in L_{2q'}^r(\mathcal{M}) \mapsto \Big( 1 \ten
d_{\varphi}^{\frac{1}{2p'}} x, 1 \ten d_{\varphi}^{\frac14} x
\Big) + \mathcal{H}_{r,0} \in
\mathcal{H}_{2p',2}^r(\mathcal{M},\theta) / \mathcal{H}_{r,0},$$
$$u_c: x d_{\varphi}^{\frac{1}{2q'}} \in
L_{2q'}^c(\mathcal{M}) \mapsto \Big( x d_{\varphi}^{\frac{1}{2p'}}
\ten 1, x d_{\varphi}^{\frac14} \ten 1 \Big) + \mathcal{H}_{c,0}
\in \mathcal{H}_{2p',2}^c(\mathcal{M},\theta) /
\mathcal{H}_{c,0}.$$
\end{proposition}

\dem This follows easily from the identity
$$\mbox{tr}_\mathcal{M} (x^* y) = \int_{\partial \mathcal{S}}
\mbox{tr}_\mathcal{M} \big( f(\overline{z})^* g(z) \big) d
\mu_\theta(z),$$ valid for any pair of analytic functions $f$ and
$g$ such that $(f(\theta), g(\theta)) = (x,y)$. Indeed, according
to the definition of the mappings $u_r, u_c, w_r, w_c$, this means
that we have
\begin{eqnarray*}
\big\langle u_r(x_1), w_r(y_1) \big\rangle & = & \langle x_1,y_1
\rangle \\ \big\langle u_c(x_2), w_c(y_2) \big\rangle & = &
\langle x_2, y_2 \rangle
\end{eqnarray*}
for any $$(x_1,x_2,y_1,y_2) \in L_{2q'}^r(\mathcal{M}) \times
L_{2q'}^c(\mathcal{M}) \times L_{\frac{2q}{q+1}}^c(\mathcal{M})
\times L_{\frac{2q}{q+1}}^r(\mathcal{M}).$$ In particular, we
deduce $$w_r^* u_r = id_{L_{2q'}^r(\mathcal{M})} \quad \mbox{and}
\quad w_c^* u_c = id_{L_{2q'}^c(\mathcal{M})}.$$ Therefore, the
assertion follows from Lemma \ref{step1} and Lemma
\ref{LemmadeDualidad}. \fin

Let $\mathcal{M}$ be a von Neumann algebra equipped with a
\emph{n.f.} state $\varphi$ and let $\mathcal{M}_m$ be the tensor
product $\mathrm{M}_m \ten \mathcal{M}$. Then, if $\mathsf{E}_m =
id_{\mathrm{M}_m} \ten \varphi: \mathcal{M}_m \to \mathrm{M}_m$
denotes the associated conditional expectation, the following
generalizes Lemma \ref{Lem-oss-Jpq}.

\begin{lemma} \label{Lemma-BYint}
If $1/r = 1/q - 1/p$, we have isometries
\begin{eqnarray*}
L_{(2r,\infty)}^{2p}(\mathcal{M}_m,\mathsf{E}_m) & = & C_{p}^m
\ten_h L_{2q}^r(\mathcal{M}) \ten_h R_m, \\
L_{(\infty,2r)}^{2p}(\mathcal{M}_m,\mathsf{E}_m) & = & C_m \ten_h
L_{2q}^c(\mathcal{M}) \ten_h R_p^m.
\end{eqnarray*}
\end{lemma}

\dem By Kouba's theorem, the spaces on the right form complex
interpolation families with respect to both indices $p$ and $q$.
Let us see that the same happens for the conditional $L_p$ spaces
on the left. Indeed, if we fix $p$ and move $1 \le q \le p$ so
that $1/q = 1-\theta + \theta/p$, we have to justify the
isometries
\begin{eqnarray*}
L_{(2r,\infty)}^{2p}(\mathcal{M}_m,\mathsf{E}_m) & = & \Big[
L_{(2p',\infty)}^{2p}(\mathcal{M}_m,\mathsf{E}_m),
L_{(\infty,\infty)}^{2p}(\mathcal{M}_m,\mathsf{E}_m)\Big]_\theta,
\\ L_{(\infty,2r)}^{2p}(\mathcal{M}_m,\mathsf{E}_m) & = & \Big[
L_{(\infty,2p')}^{2p}(\mathcal{M}_m,\mathsf{E}_m),
L_{(\infty,\infty)}^{2p}(\mathcal{M}_m,\mathsf{E}_m)\Big]_\theta.
\end{eqnarray*}
As far as $p$ is finite this is part of Theorem \ref{TheoApp1},
while the remaining case follows from Lemma \ref{Lemintinfty}.
This means that it suffices to prove the assertion for $q=1$ since
the case $q=p$ follows from the trivial isometries
$$\begin{array}{rclcl} L_{2p}(\mathcal{M}_m) & = &
L_{2p}^r(\mathcal{M}_m) & = & C_p^m \ten_h L_{2p}^r(\mathcal{M})
\ten_h R_m, \\ [3pt] L_{2p}(\mathcal{M}_m) & = &
L_{2p}^c(\mathcal{M}_m) & = & C_m \ten_h L_{2p}^c(\mathcal{M})
\ten_h R_p^m.
\end{array}$$
Now we claim that we also have
\begin{eqnarray*}
L_{(2p',\infty)}^{2p}(\mathcal{M}_m,\mathsf{E}_m) & = & \Big[
L_{(\infty,\infty)}^{2}(\mathcal{M}_m,\mathsf{E}_m),
L_{(2,\infty)}^{\infty}(\mathcal{M}_m,\mathsf{E}_m)\Big]_{\frac{1}{p'}},
\\ L_{(\infty,2p')}^{2p}(\mathcal{M}_m,\mathsf{E}_m) & = & \Big[
L_{(\infty,\infty)}^{2}(\mathcal{M}_m,\mathsf{E}_m),
L_{(\infty,2)}^{\infty}(\mathcal{M}_m,\mathsf{E}_m)\Big]_{\frac{1}{p'}}.
\end{eqnarray*}
Indeed, recalling that
$$L_{(\infty,\infty)}^{2}(\mathcal{M}_m,\mathsf{E}_m) =
L_{2}(\mathcal{M}_m),$$ we deduce that both spaces above are
reflexive. Therefore, using Theorem \ref{TheoApp1} (b) for
amalgamated $L_p$ spaces and duality, we deduce our claim and we
are reduced to show that
$$L_{(2,\infty)}^{\infty}(\mathcal{M}_m,\mathsf{E}_m) = \mathrm{M}_m
\big( L_{2}^r(\mathcal{M}) \big) \quad \mbox{and} \quad
L_{(\infty,2)}^{\infty}(\mathcal{M}_m,\mathsf{E}_m) = \mathrm{M}_m
\big( L_{2}^c(\mathcal{M}) \big).$$ However, these isometries are
exactly \eqref{Eq-Schatten-RC}. The proof is complete. \fin

Let $\mathcal{M}$ be a von Neumann algebra equipped with a
\emph{n.f.} state $\varphi$ and let $\mathcal{N}$ be a von Neumann
subalgebra of $\mathcal{M}$. Let $\mathsf{E}: \mathcal{M} \to
\mathcal{N}$ denote the corresponding conditional expectation. In
Section \ref{NewSect3}, we defined the spaces
\begin{eqnarray*} \mathcal{R}_{2p,q}^n(\mathcal{M}, \mathsf{E}) &
= & n^{\frac{1}{2p}} \, L_{2p}^r(\mathcal{M}) \, \cap \,
n^{\frac{1}{2q}} \, L_{(2r,\infty)}^{2p} (\mathcal{M},
\mathsf{E}), \\ \mathcal{C}_{2p,q}^n \, (\mathcal{M}, \mathsf{E})
& = & n^{\frac{1}{2p}} \, L_{2p}^c(\mathcal{M}) \, \cap \,
n^{\frac{1}{2q}} \, L_{(\infty,2r)}^{2p} (\mathcal{M},
\mathsf{E}),
\end{eqnarray*}
and mentioned the isomorphism from \cite{JP2}
\begin{equation} \label{Eq-Amalgam-Factor}
\mathcal{J}_{p,q}^n(\mathcal{M}, \mathsf{E}) \simeq
\mathcal{R}_{2p,q}^n(\mathcal{M}, \mathsf{E}) \oM
\mathcal{C}_{2p,q}^n(\mathcal{M}, \mathsf{E}).
\end{equation}
In fact, to be completely fair we should say that we have slightly
modified the definition of $\mathcal{R}_{2p,q}^n(\mathcal{M},
\mathsf{E})$ and $\mathcal{C}_{2p,q}^n(\mathcal{M}, \mathsf{E})$
by considering the row/column o.s.s. of $L_{2p}(\mathcal{M})$.
However, the new definition coincides with the former one in the
Banach space level. Hence, since we do not even have an operator
space structure for these spaces, our modification is only
motivated for notational convenience below. Namely, inspired by
Lemma \ref{Lemma-BYint}, we introduce the operator spaces
\begin{eqnarray*}
\mathcal{R}_{2p,q}^n(\mathcal{M}) & = & n^{\frac{1}{2p}} \,
L_{2p}^r(\mathcal{M}) \, \cap \, n^{\frac{1}{2q}} \,
L_{2q}^r (\mathcal{M}), \\
\mathcal{C}_{2p,q}^n \, (\mathcal{M}) & = & n^{\frac{1}{2p}} \,
L_{2p}^c(\mathcal{M}) \, \cap \, n^{\frac{1}{2q}} \, L_{2q}^c
(\mathcal{M}).
\end{eqnarray*}
These spaces give rise to the complete isomorphism
\begin{equation} \label{Eq-CB-Factor}
\mathcal{J}_{p,q}^n(\mathcal{M}) \simeq_{cb}
\mathcal{R}_{2p,q}^n(\mathcal{M}) \ten_{\mathcal{M},h}
\mathcal{C}_{2p,q}^n(\mathcal{M}).
\end{equation}
Indeed, taking $(\mathcal{M}, \mathcal{N}, \mathsf{E}) =
(\mathcal{M}_m, \mathrm{M}_m, \mathsf{E}_m)$ in
\eqref{Eq-Amalgam-Factor} we have
\begin{eqnarray*}
S_p^m \big( \mathcal{J}_{p,q}^n(\mathcal{M}) \big) & = &
\mathcal{J}_{p,q}^n(\mathcal{M}_m, \mathsf{E}_m) \\ [3pt] & \simeq
& \mathcal{R}_{2p,q}^n(\mathcal{M}_m, \mathsf{E}_m)
\ten_{\mathcal{M}_m} \mathcal{C}_{2p,q}^n(\mathcal{M}_m,
\mathsf{E}_m) \\ [2pt] & = & S_{(2p,\infty)}^m \big(
\mathcal{R}_{2p,q}^n(\mathcal{M}) \big) \ten_{\mathcal{M}_m}
S_{(\infty,2p)}^m \big( \mathcal{C}_{2p,q}^n (\mathcal{M}) \big)
\\ & = & C_p^m \ten_h \Big( \mathcal{R}_{2p,q}^n(\mathcal{M})
\ten_{\mathcal{M},h} \mathcal{C}_{2p,q}^n (\mathcal{M}) \Big)
\ten_h R_p^m \\ & = & S_p^m \big(
\mathcal{R}_{2p,q}^n(\mathcal{M}) \ten_{\mathcal{M},h}
\mathcal{C}_{2p,q}^n(\mathcal{M}) \big),
\end{eqnarray*}
where the third identity follows from Lemma \ref{Lemma-BYint}
after taking in consideration our \emph{new} definition of
$\mathcal{R}_{2p,q}^n(\mathcal{M}, \mathsf{E})$ and
$\mathcal{C}_{2p,q}^n(\mathcal{M}, \mathsf{E})$. In other words,
\eqref{Eq-Amalgam-Factor} and \eqref{Eq-CB-Factor} are the
amalgamated and operator space versions of the same factorization
isomorphism. Now assume that $\mathcal{M}$ is equipped with a
\emph{n.s.s.f.} weight $\psi$, given by the increasing sequence
$(\psi_n, q_n)_{n \ge 1}$. Then, we may generalize the
factorization result above in the usual way. Namely, assuming by
approximation that $\mathrm{k}_n = \psi_n(1)$ are positive
integers, we define
\begin{eqnarray*}
\mathcal{R}_{2p,q}(\psi_n) & = & \mathrm{k}_n^{\frac{1}{2p}}
L_{2p}^r (q_n \mathcal{M} q_n) \cap \mathrm{k}_n^{\frac{1}{2q}}
L_{2q}^r (q_n \mathcal{M} q_n), \\ \mathcal{C}_{2p,q}(\psi_n) & =
& \mathrm{k}_n^{\frac{1}{2p}} L_{2p}^c (q_n \mathcal{M} q_n) \cap
\mathrm{k}_n^{\frac{1}{2q}} L_{2q}^c (q_n \mathcal{M} q_n).
\end{eqnarray*}
This gives the complete isomorphism
$$\mathcal{J}_{p,q}(\psi_n) \simeq_{cb} \mathcal{R}_{2p,q}(\psi_n)
\ten_{\mathcal{M},h} \mathcal{C}_{2p,q}(\psi_n) = \bigcap_{u,v \in
\{2p,2q\}} \mathrm{k}_n^{\frac1u + \frac1v} L_{(u,v)}(q_n
\mathcal{M} q_n).$$ Then, taking direct limits we obtain the space
$$\mathcal{J}_{p,q}(\psi) = \mathcal{R}_{2p,q}(\psi) \ten_{\mathcal{M},h}
\mathcal{C}_{2p,q}(\psi) = \bigcap_{u,v \in \{2p,2q\}}
L_{(u,v)}(\mathcal{M}).$$

\begin{lemma} \label{psi} Let $\mathcal{M}$ be a von Neumann
algebra equipped with a n.s.s.f. weight $\psi$. Then, there exists
a n.s.s.f. weight $\xi$ on $\mathcal{B}(\ell_2)$ such that the
following complete isomorphisms hold
\begin{eqnarray*}
\mathcal{H}_{2p',2}^r(\mathcal{M},\theta) \ten_h R & \simeq_{cb} &
\mathcal{R}_{2p',2}(\psi \ten \xi), \\ C \ten_h
\mathcal{H}_{2p',2}^c(\mathcal{M},\theta) & \simeq_{cb} & \,
\mathcal{C}_{2p',2} \, (\psi \ten \xi).
\end{eqnarray*}
\end{lemma}

\dem By symmetry, we only consider the column case. Let us first
observe that $\mathcal{H}_2$ is indeed the graph of an injective
closed densely-defined (unbounded) operator with dense range. This
is quite similar to Remark \ref{Remark-Graph_Sq}. It follows from
the three lines lemma that for $z=a+ib$
$$|f(z)| \le \|f_{|_{\partial_0}}\|_{L_2(\partial_0,\mu_a)}^{1-a}
\|f_{\mid_{\partial_1}}\|_{L_2(\partial_1,\mu_a)}^a.$$ Since
$\mu_a$ and $\mu_{\theta}$ have the same null sets, we deduce that
$$\pi_j(f) = f_{\mid_{\partial_j}} \in L_2(\partial_j, \mu_\theta)$$
are injective for $j=0,1$ when restricted to analytic functions.
Thus, the mapping $\Lambda(\pi_0(f)) = \pi_1(f)$ is an injective
closed densely-defined operator with dense range and
$\mathcal{H}_2$ is its graph. Let $\Lambda = u |\Lambda|$ be the
polar decomposition. Since $M_u = 1 \ten u^*$ defines a complete
isometry (recall that $\Lambda$ has dense range)
$$L_4^c(\mathcal{M}) \ten_h L_2^{oh}(\partial_1) \to
L_4^c(\mathcal{M}) \ten_h L_2^{oh}(\partial_0),$$ we may replace
$\Lambda$ by $|\Lambda|$ in the definition of
$\mathcal{H}_{2p',2}^c(\mathcal{M},\theta)$. Using the
discretization Lemma \ref{Lemma-Diagonalization}, we may also
replace $L_2(\partial_0)$ by $\ell_2$ and the operator $|\Lambda|$
by a diagonal operator $\mathsf{d}_{\lambda}$. These
considerations provide a cb-isomorphism
$$C \ten_h \mathcal{H}_{2p',2}^c(\mathcal{M},\theta) \simeq_{cb}
\Big( C \ten_h L_{2p'}^c (\mathcal{M}) \ten_h R_{p'} \Big) \cap
\Big( C \ten_h L_4^c(\mathcal{M}) \ten_h \ell_2^{oh}(\lambda)
\Big),$$ where $\ell_2^{oh}(\lambda)$ is the weighted form of
$\mathrm{OH}$ which arises from the action of
$\mathsf{d}_\lambda$. The assertion follows by a direct limit
argument. Indeed, the \emph{n.s.s.f.} weight $\psi$ on
$\mathcal{M}$ is given by the sequence $(\psi_n, q_n)_{n \ge 1}$.
On the other hand, we may consider the \emph{n.s.s.f.} weight
$\xi$ on $\mathcal{B}(\ell_2)$ determined by the sequence
$(\xi_n,\pi_n)_{n \ge 1}$, where $\pi_n$ is the projection onto
the first $n$ coordinates and $\xi_n$ is the finite weight on
$\pi_n \mathcal{B}(\ell_2) \pi_n$ given by $$\xi_n \Big( \pi_n
\big( \summ_{ij} x_{ij} e_{ij} \big) \pi_n \Big) = \sum_{k=1}^n
\gamma_k x_{kk} \quad \mbox{with} \quad \gamma_k^{\frac14 -
\frac{1}{2p'}} = \lambda_k.$$ Let us define the parameters
$\mathrm{k}_n' = \xi_n(1)$ and $\mathrm{w}_n = \mathrm{k}_n
\mathrm{k}_n'$. Then, arguing as we did in Lemma
\ref{Lemma-Direct-Sum}, it turns out that the intersection space
above is the direct limit of the following sequence of spaces
$$\mathrm{w}_n^{\frac{1}{2p'}}
\Big( L_{2p'}^c \big( q_n \mathcal{M} q_n \bar\ten \pi_n
\mathcal{B}(\ell_2) \pi_n \big) \Big) \cap
\mathrm{w}_n^{\frac{1}{4}} \Big( L_{4}^c \big( q_n \mathcal{M} q_n
\bar\ten \pi_n \mathcal{B}(\ell_2) \pi_n \big) \Big).$$ However,
the latter space is $\mathcal{C}_{2p',2}(\psi_n \ten \xi_n)$. This
completes the proof. \fin

\begin{proposition} \label{Proposition-phipsixi}
The predual space of $$\mathcal{H}_{2p',2}^r (\mathcal{M}, \theta)
\ten_{\mathcal{M},h} \mathcal{H}_{2p',2}^c(\mathcal{M},\theta)$$
embeds completely isomorphically into $\mathcal{K}_{r \! c_p}^p
(\phi \ten
 \psi \ten \xi)$ for some n.s.s.f. weight $\xi$ on
$\mathcal{B}(\ell_2)$ and where $\phi$ denotes the quasi-free
state over the hyperfinite $\mathrm{III}_1$ factor $\mathcal{R}$
considered in Proposition $\ref{Prop-k4k3}$.
\end{proposition}

\dem According to Lemma \ref{psi}
\begin{eqnarray*}
\lefteqn{\mathcal{H}_{2p',2}^r (\mathcal{M}, \theta)
\ten_{\mathcal{M},h} \mathcal{H}_{2p',2}^c(\mathcal{M},\theta)} \\
& = \ \, \null & \Big( \mathcal{H}_{2p',2}^r (\mathcal{M}, \theta)
\ten_h R \Big) \ten_{\mathcal{M} \bar{\ten} \mathcal{B}(\ell_2),h}
\Big( C \ten_h \mathcal{H}_{2p',2}^c(\mathcal{M},\theta) \Big) \\
& \simeq_{cb} & \mathcal{R}_{2p',2}(\psi \ten \xi)
\ten_{\mathcal{M} \bar{\ten} \mathcal{B}(\ell_2),h}
\mathcal{C}_{2p',2}(\psi \ten \xi) \ \simeq_{cb} \
\mathcal{J}_{p',2}(\psi \ten \xi).
\end{eqnarray*}
However, $\mathcal{J}_{p',2}(\psi \ten \xi)$ is a direct limit of
spaces
$$\mathcal{J}_{p',2} (\psi_n \ten \xi_n) =
\mathcal{J}_{p',2}^{\mathrm{w}_n} \Big( q_n \mathcal{M} q_n
\bar\ten \pi_n \mathcal{B}(\ell_2) \pi_n \Big).$$ According to
Proposition \ref{Prop-k4k3}, the direct limit
$$\mathcal{K}_{p,2}(\psi \ten \xi) = \limm_n
\mathcal{K}_{p,2}^n(\psi_n \ten \xi_n)$$ of the corresponding
predual spaces cb-embeds into $\mathcal{K}_{r \! c_p}^p(\phi \ten
\psi \ten \xi)$. \fin

Now we are ready to prove our main result. In the proof we shall
need to work with certain quotient of
$\mathcal{H}_{2p',2}^r(\mathcal{M},\theta) \ten_{\mathcal{M},h}
\mathcal{H}_{2p',2}^c(\mathcal{M},\theta)$. Namely, recalling the
subspaces $\mathcal{H}_{r,0}$ and $\mathcal{H}_{c,0}$, we set
$$\mathcal{Q}_{2p',2}(\mathcal{M},\theta) = \big(
\mathcal{H}_{2p',2}^r(\mathcal{M},\theta) / \mathcal{H}_{r,0}
\big) \ten_{\mathcal{M},h} \big(
\mathcal{H}_{2p',2}^c(\mathcal{M},\theta) / \mathcal{H}_{c,0}
\big).$$ We claim that $\mathcal{Q}_{2p',2}(\mathcal{M},\theta)$
is a quotient of $\mathcal{H}_{2p',2}^r(\mathcal{M},\theta)
\ten_{\mathcal{M},h} \mathcal{H}_{2p',2}^c(\mathcal{M},\theta)$.
Indeed, according to the definition of the
$\mathcal{M}$-amalgamated Haagerup tensor product of two operator
spaces (see Remark \ref{Re-oM}), we may write
$\mathcal{Q}_{2p',2}(\mathcal{M},\theta)$ as a quotient of the
Haagerup tensor product
$$\Lambda_{2p',2}(\mathcal{M},\theta) = \big(
\mathcal{H}_{2p',2}^r(\mathcal{M},\theta) / \mathcal{H}_{r,0}
\big) \ten_h \big( \mathcal{H}_{2p',2}^c(\mathcal{M},\theta) /
\mathcal{H}_{c,0} \big)$$ by the closed subspace spanned by the
differences $x_1 \gamma \ten x_2 - x_1 \ten \gamma x_2$, with
$\gamma \in \mathcal{M}$. Therefore, it suffices to see that the
space $\Lambda_{2p',2}(\mathcal{M},\theta)$ is a quotient of
$\mathcal{H}_{2p',2}^r(\mathcal{M},\theta) \ten_h
\mathcal{H}_{2p',2}^c(\mathcal{M},\theta)$. However, this follows
from the projectivity of the Haagerup tensor product and our claim
follows.


\begin{theorem}
Let $1 \le p < q \le 2$ and let $\mathcal{M}$ be a von Neumann
algebra. Then, there exists a sufficiently large von Neumann
algebra $\mathcal{A}$ and a completely isomorphic embedding of
$L_q(\mathcal{M})$ into $L_p(\mathcal{A})$, where both spaces are
equipped with their respective natural operator space structures.
Moreover, we have
\begin{itemize}
\item[\textbf{(a)}] If $\mathcal{M}$ is $\mathrm{QWEP}$, we can
choose $\mathcal{A}$ to be $\mathrm{QWEP}$.

\item[\textbf{(b)}] If $\mathcal{M}$ is hyperfinite, we can choose
$\mathcal{A}$ to be hyperfinite.
\end{itemize}
\end{theorem}

\dem Let us first assume that $\mathcal{M}$ is finite. According
to Theorem \ref{Kp} and Proposition \ref{Proposition-phipsixi}, it
suffices to prove that $L_{q'}(\mathcal{M})$ is completely
isomorphic to a quotient of
$\mathcal{H}_{2p',2}^r(\mathcal{M},\theta) \ten_{\mathcal{M}}
\mathcal{H}_{2p',2}^c(\mathcal{M},\theta)$. This follows from
Proposition \ref{Proposition-CBIsomorphisms-Quo} since
$$L_{q'}(\mathcal{M}) \simeq_{cb} L_{2q'}^r(\mathcal{M})
\ten_{\mathcal{M},h} L_{2q'}^c(\mathcal{M}) \simeq_{cb}
\mathcal{Q}_{2p',2}(\mathcal{M},\theta).$$ The construction of the
cb-embedding for a general von Neumann algebra $\mathcal{M}$ can
be obtained by using Haagerup's approximation theorem \cite{H2}
and the fact that direct limits are stable in our construction.
Indeed, Haagerup theorem shows that for every $\sigma$-finite von
Neumann algebra $\mathcal{M}$, the space $L_q(\mathcal{M})$ is
complemented in a direct limit of $L_q$ spaces over finite von
Neumann algebras. Finally, if $\mathcal{M}$ is any von Neumann
algebra, we observe that $L_q(\mathcal{M})$ can always be written
as a direct limit of $L_q$ spaces associated to $\sigma$-finite
von Neumann algebras. On the other hand, the stability of
hyperfiniteness follows directly from our construction. Indeed,
our construction goes as follows
$$L_q(\mathcal{M}) \to \Big( \mathcal{H}_{2p',2}^r (\mathcal{M},
\theta) \ten_{\mathcal{M},h}
\mathcal{H}_{2p',2}^c(\mathcal{M},\theta) \Big)_* \to
\mathcal{K}_{r \! c_p}^p(\phi \ten \psi \ten \xi) \to
L_p(\mathcal{A})$$ where the first embedding follows as above, the
second from Proposition \ref{Proposition-phipsixi} and the last
one from Theorem \ref{Kp}. In particular, it turns out that the
von Neumann algebra $\mathcal{A}$ is of the form
$$\mathcal{A} = M_s \big( \mathcal{R} \bar\ten \mathcal{M}
\bar\ten \mathcal{B}(\ell_2) \big),$$ which is hyperfinite when
$\mathcal{M}$ is hyperfinite and is a factor when $\mathcal{M}$ is
a factor. Finally, it remains to justify that the $\mathrm{QWEP}$
is preserved. If $\mathcal{M}$ is $\mathrm{QWEP}$, there exists a
completely isometric embedding of $L_q(\mathcal{M})$ into
$L_q(\mathcal{M}_\mathcal{U})$ with $\mathcal{M}_\mathcal{U}$ of
the form $$\mathcal{M}_\mathcal{U} = \Big( \prodd_{n,\mathcal{U}}
S_1^n \Big)^*.$$ Since we know from Corollary \ref{SqQWEP} that
the Schatten class $S_q^n$ embeds completely isomorphically into
$L_p(\mathcal{A}_n)$ with relevant constants independent of $n$
and $\mathcal{A}_n$ being $\mathrm{QWEP}$, we find a completely
isomorphic embedding
$$L_q(\mathcal{M}) \to L_p(\mathcal{A}_\mathcal{U}) \quad
\mbox{with} \quad \mathcal{A}_\mathcal{U} = \Big( \prodd_{n,
\mathcal{U}} {\mathcal{A}_n}_* \Big)^*.$$ This proves the
assertion since $\mathcal{A}_\mathcal{U}$ is $\mathrm{QWEP}$. The
proof is complete. \fin


\bibliographystyle{amsplain}

\begin{thebibliography}{10}
\bibitem {AW} H. Araki and E.J. Woods, A classification of factors.
Publ. Res. Inst. Math. Sci. Ser. A \textbf{4} (1968/69), 51-130.
\bibitem {BDK} J. Bretagnolle, D. Dacunha-Castelle and J.L.
Krivine, Lois stable et space $L^p$. Ann. Inst. H. Poincar{\'e}
\textbf{2} (1966), 231-259.
\bibitem {Co} A. Connes, Une classification des facteurs
de type III. Ann. Sci. \'{E}cole Norm. Sup. \textbf{6} (1973),
133-252.
\bibitem {ER} E.G. Effros and Z.J. Ruan, Operator Spaces. London
Math. Soc. Monogr. \textbf{23}, Oxford University Press, 2000.
\bibitem {H2} U. Haagerup, Non-commutative integration theory.
Unpublished manuscript (1978). See also Haagerup's Lecture given
at the Symposium in Pure Mathematics of the Amer. Math. Soc.
Queens University, Kingston, Ontario, 1980.
\bibitem {H} U. Haagerup, $L_p$ spaces associated with an
arbitrary von Neumann algebra. Alg\`{e}bres d'op\'{e}rateurs et
leurs applications en physique math\'{e}matique, CNRS (1979),
175-184.
\bibitem {HRS} U. Haagerup, H.P. Rosenthal and F.A. Sukochev,
Banach Embedding Properties of Non-commutative $L_p$-Spaces. Mem.
Amer. Math. Soc. \textbf{163}, 2003.
\bibitem {J00} M. Junge, Factorization theory for spaces of
operators. Habilitation thesis. Kiel university, 1996.
\bibitem {J0} M. Junge, Embeddings of non-commutative $L_p$-spaces
into non-commutative $L_1$-spaces, $1 < p < 2$. Geom. Funct. Anal.
\textbf{10} (2000), 389-406.
\bibitem {J1} M. Junge, Doob's inequality for non-commutative
martingales. J. reine angew. Math. \textbf{549} (2002), 149-190.
\bibitem {J2} M. Junge, Embedding of the operator space
$\mathrm{OH}$ and the logarithmic \lq little Grothendieck
inequality\rq. Invent. Math. \textbf{161} (2005), 225-286.
\bibitem {J5} M. Junge, Fubini's theorem for ultraproducts of
noncommutative $L_p$-spaces. Canad. J. Math. \textbf{56} (2004),
983-1021.
\bibitem {J3} M. Junge, Operator spaces and Araki-Woods factors.
Preprint. ArXiv: math.OA/0504255.
\bibitem {J6} M. Junge, Noncommutative Poisson random measure.
In progress.
\bibitem {JP} M. Junge and J. Parcet, The norm of sums of
independent noncommutative random variables in $L_p(\ell_1)$. J.
Funct. Anal. \textbf{221} (2005), 366-406.
\bibitem {JP2} M. Junge and J. Parcet, Theory of Amalgamated $L_p$
Spaces in Noncommutative Probability. Preprint. ArXiv:
math.OA/0511406.
\bibitem {JP3} M. Junge and J. Parcet, Rosenthal's theorem for
subspaces of noncommutatuve $L_p$. Preprint. ArXiv:
math.FA/0604510.
\bibitem {JPX} M. Junge, J. Parcet and Q. Xu, Rosenthal type
inequalities for free chaos. Preprint. ArXiv: math.OA/0511732.
\bibitem {JX} M. Junge and Q. Xu, Noncommutative
Burkholder/Rosenthal inequalities. Ann. Probab. \textbf{31}
(2003), 948-995.
\bibitem {JX2} M. Junge and Q. Xu, Noncommutative maximal ergodic
theorems. To appear in J. Amer. Math. Soc.
\bibitem {JX4} M. Junge and Q. Xu, Noncommutative
Burkholder/Rosenthal inequalities II: Applications. Preprint,
2005.
\bibitem {KR} R.V. Kadison and J.R. Ringrose, Fundamentals of the
Theory of Operator Algebras I and II. Grad. Stud. Math.,
\textbf{15} \& \textbf{16}, American Mathematical Society, 1997.
\bibitem {Ko} H. Kosaki, Applications of the complex
interpolation method to a von Neumann algebra. J. Funct. Anal.
\textbf{56} (1984), 29-78.
\bibitem {LuP} F. Lust-Piquard and G. Pisier, Non-commutative
Khintchine and Paley inequalities. Ark. Mat. \textbf{29} (1991),
241-260.
\bibitem {MaP} B. Maurey and G. Pisier, S\'{e}ries de variables
al\'{e}atoires vectorielles ind\'{e}pendantes et
propri\'{e}t\'{e}s g\'{e}om\'{e}triques des espaces de Banach.
Studia Math. \textbf{58} (1976), 45-90.
\bibitem {Pa} J. Parcet, $\mathrm{B}$-convex operator spaces.
Proc. Edinburgh Math. Soc. \textbf{46} (2003), 649-668.
\bibitem {PP} J. Parcet and G. Pisier, Non-commutative Khintchine
type inequalities associated with free groups. Indiana Univ. Math.
J. \textbf{54} (2005), 531-556.
\bibitem {PR} J. Parcet and N. Randrianantoanina, Gundy's
decomposition  for non-commutative martingales and applications.
Proc. London Math. Soc. \textbf{93} (2006), 227-252.
\bibitem {Pau} V. Paulsen, Completely bounded maps and operator
algebras. Cambridge Studies in Advanced Mathematics \textbf{78}.
Cambridge University Press, 2002.
\bibitem {Pis} G. Pisier, Factorization of operators
through $L_{p\infty}$ or $L_{p1}$ and non-commutative
generalizations. Math. Ann. \textbf{276} (1986), 105-136.
\bibitem {P0} G. Pisier, Projections from a von Neumann algebra
onto a subalgebra. Bull. Soc. Math. France \textbf{123} (1995),
139-153.
\bibitem {P1} G. Pisier, The Operator Hilbert Space OH, Complex
Interpolation and Tensor Norms. Mem. Amer. Math. Soc. \textbf{122}
(1996).
\bibitem {P2} G. Pisier, Non-Commutative Vector Valued
$L_p$-Spaces and Completely $p$-Summing Maps. Ast\'{e}risque
\textbf{247} (1998).
\bibitem {P3} G. Pisier, Introduction to Operator Space Theory.
Cambridge University Press, 2003.
\bibitem {P5} G. Pisier, The operator Hilbert space OH and type
III von Neumann algebras. Bull. London Math. Soc. \textbf{36}
(2004), 455-459.
\bibitem {P4} G. Pisier, Completely bounded maps into certain
Hilbertian operator spaces. Internat. Math. Res. Notices
\textbf{74} (2004), 3983-4018.
\bibitem {PS} G. Pisier and D. Shlyakhtenko, Grothendieck's
theorem for operator spaces. Invent. Math. \textbf{150} (2002),
185-217.
\bibitem {PX1} G. Pisier and Q. Xu, Non-commutative martingale
inequalities. Comm. Math. Phys. \textbf{189} (1997), 667-698.
\bibitem {PX2} G. Pisier and Q. Xu, Non-commutative $L_p$-spaces.
Handbook of the Geometry of Banach Spaces II (Ed. W.B. Johnson and
J. Lindenstrauss) North-Holland (2003), 1459-1517.
\bibitem {Narcisse} N. Randrianantoanina,
Personal communication.
\bibitem {Ra} Y. Raynaud, On ultrapowers of non-commutative $L_p$
spaces. J. Operator Theory \textbf{48} (2002), 41-68.
\bibitem {Ro0} H.P. Rosenthal, On the subspaces of $L^p$ $(p>2)$
spanned by sequences of independent random variables. Israel J.
Math. \textbf{8} (1970), 273-303.
\bibitem {Ro} H.P. Rosenthal, On subspaces of $L_p$. Ann. of Math.
\textbf{97} (1973), 344-373.
\bibitem {R} Z.J. Ruan, Subspaces of $C^*$-algebras. J. Funct.
Anal. \textbf{76} (1988), 217-230.
\bibitem {S} D. Shlyakhtenko, Free quasi-free states. Pacific J.
Math. \textbf{177} (1997), 329-368.
\bibitem {Ta2} M. Takesaki, Conditional expectations in von
Neumann algebras. J. Func. Anal. \textbf{9} (1972), 306-321.
\bibitem {T1} M. Terp, $L_p$ spaces associated with von Neumann
algebras. Math. Institute Copenhagen University, 1981.
\bibitem {V2} D.V. Voiculescu, A strengthened asymptotic
freeness result for random matrices with applications to free
entropy. Internat. Math. Res. Notices \textbf{1} (1998), 41-63.
\bibitem {VDN} D.V. Voiculescu, K. Dykema and A. Nica, Free
random variables. CRM Monograph Series \textbf{1}, American
Mathematical Society, 1992.
\bibitem {X4} Q. Xu, Operator-space Grothendieck inequalities
for noncommutative $L_p$-spaces. Duke Math. J. \textbf{131}
(2006), 525-574.
\bibitem {X3} Q. Xu, Embedding of $C_q$ and $R_q$ into
noncommutative $L_p$-spaces, $1 \le p < q \le 2$. Preprint. ArXiv:
math.FA/0505307.
\bibitem {X2} Q. Xu, Real interpolation approach to Junge's
works on embedding of $\mathrm{OH}$ and the little Grothendieck
inequality. Preprint, 2004.
\bibitem {X} Q. Xu, A description of $(C_p[L_p(M)],
R_p[L_p(M)])_{\theta}$. Proc. Roy. Soc. Edinburgh Sect. A
\textbf{135} (2005), 1073-1083.
\end{thebibliography}

\end{document}